\documentclass[onefignum,onetabnum]{siamart251216}

\usepackage{algorithmic}
\Crefname{ALC@unique}{Line}{Lines}
\usepackage{booktabs}
\usepackage{array}
\usepackage{makecell}
\usepackage{jabbrv}
\DefineSpuriousJournalWord{on}
\DefineSpuriousJournalWord{Its}
\usepackage{rotating}
\usepackage{amsopn}
\usepackage{lipsum}
\usepackage{epstopdf}
\usepackage{graphicx}
\usepackage[caption=false]{subfig}
\usepackage{amsfonts}
\usepackage{amssymb}
\usepackage{amsmath}
\usepackage{tabu}
\usepackage[utf8]{inputenc}
\usepackage[english]{babel}
\usepackage{mathtools}
\usepackage{listings}
\usepackage{enumitem}
\usepackage{hyperref}
\usepackage{bm}
\usepackage{microtype}
\usepackage{overpic}
\usepackage{cite}
\setlist[enumerate]{nosep}
\setlist[enumerate]{topsep=0pt}

\usepackage{verbatimbox}

\newcommand*{\dd}{{\,\mathrm{d}}}
\newcommand{\surf}{M}
\renewcommand{\vec}[1]{{\bm{#1}}}

\DeclareMathOperator{\im}{Im}

\DeclareMathOperator{\area}{area}
\DeclareMathOperator{\length}{length}
\DeclareMathOperator{\dist}{dist}

\DeclareMathOperator{\spec}{Sp}
\DeclareMathOperator{\res}{res}

\usepackage{matlab-prettifier}
\lstdefinestyle{matlab-custom}{
	language=Matlab,
	basicstyle=\footnotesize\ttfamily,
	keywordstyle=\bfseries\color{green!40!black},
	commentstyle=\itshape\color{purple!40!black},
	identifierstyle=\color{blue},
	stringstyle=\color{orange}
}

\begin{document}

\headers{Computing Spectra of Operators on Surfaces}{G. Conradie, M. J. Colbrook, and D. Fortunato}

\title{Computing Spectral Properties of Differential Operators on Surfaces}

\author{Gustav Conradie\thanks{Department of Applied Mathematics and Theoretical Physics, University of Cambridge, Cambridge, CB3 0WA, UK 
  (\email{gjc51@cam.ac.uk}, \email{m.colbrook@damtp.cam.ac.uk}).}
  \and Matthew J. Colbrook\footnotemark[1]
\and Daniel Fortunato\thanks{Center for Computational Mathematics and Center for Computational Biology, Flatiron Institute, New York, NY 10010, USA (\email{dfortunato@flatironinstitute.org}).}}

\maketitle

\begin{abstract}
Computing spectra of differential operators on surfaces is crucial in many applications, ranging from medical image analysis to nanoscale quantum structures. However, discretization of the infinite-dimensional operator and the underlying surface can lead to inaccurate and misleading results. This paper introduces the \texttt{SurfSpec} suite of algorithms to compute eigenvalues and spectral measures of operators on surfaces with error control and convergence guarantees. We combine methods based on complex contour integration and residuals to compute discrete spectra and high-frequency eigenmodes, including of nonlinear eigenvalue problems. We then compute convolutions of spectral measures with rational kernels to deal with continuous spectra. In both cases, the key tool is the operator resolvent applied to functions, which is computed using high-order methods that solve surface PDEs. We illustrate our algorithms on a variety of differential operators and surfaces, including a triceratops.
\end{abstract}


\vspace{-2mm}


\vspace{-2mm}

\linespread{0.94}
\section{Introduction}

Given a surface $M \subset \mathbb{R}^3$ and a Hilbert space $\mathcal{H}$ of functions on $M$, we study the spectral properties of a linear differential operator $\mathcal{L}$ with domain $\mathcal{D}(\mathcal{L}) \subset \mathcal{H}$. The associated eigenvalue problem is
\begin{equation}\setlength\abovedisplayskip{6pt}\setlength\belowdisplayskip{6pt}\setlength\abovedisplayskip{6pt}\setlength\belowdisplayskip{6pt}
    \label{eqn:intro_diff_eig_prob}
    \mathcal{L}f=\lambda f,\quad f\in\mathcal{D}(\mathcal{L}),\quad \lambda\in\mathbb{C}.
\end{equation}
For surfaces with boundary, boundary conditions such as Dirichlet or Neumann may be imposed. Since $\mathcal{H}$ is infinite-dimensional, the spectrum $\spec(\mathcal{L})$, defined formally as
\begin{equation}\setlength\abovedisplayskip{6pt}\setlength\belowdisplayskip{6pt}\setlength\abovedisplayskip{6pt}\setlength\belowdisplayskip{6pt}
    \spec(\mathcal{L})=\{z\in\mathbb{C}:z\mathcal{I}-\mathcal{L}\text{ is not boundedly invertible}\},
\end{equation}
may include continuous components in addition to eigenvalues, giving rise to a more general and difficult spectral problem than \cref{eqn:intro_diff_eig_prob}. Accurately computing the spectra of operators is crucial, as surface spectral problems (SSPs) play a key role in various applications. Eigenvalues and eigenfunctions of the Laplace--Beltrami operator are central to spectral shape analysis, serving as surface fingerprints or ``Shape-DNA'' \cite{reuter2005laplace,reuter2006laplace}, and are used for segmentation \cite{reuter2009discrete}, local registration \cite{hamidian2019surface}, heat kernel smoothing \cite{seo2010heat}, and constructing geometry-adapted function bases \cite{levy2006laplace}. These methods are widely applied in medical image analysis \cite{seo2011laplace,niethammer2007global,tan2014spectral} and computer graphics \cite{caissard2019laplace,zhang2007spectral}. Beyond shape analysis, Laplace--Beltrami spectra are used for dimension reduction and data representation \cite{belkin2003laplacian,coifman2006diffusion}, and have been applied in modeling cell behavior \cite{nwogbaga2025protein}. Spectral analysis of other differential operators is also important. For instance, the spectrum of a Schrödinger operator with curvature-dependent potential arises in the stability analysis of interfaces in reaction–diffusion systems \cite{alikakos1996critical,harrell1996second,harrell1998laplace}, and in modeling constrained particles \cite{da1981quantum} or nanoscale quantum structures \cite{duclos1995curvature,exner1999optimal}. Computational tools for computing spectra can also inform conjectures in spectral geometry or, when error control is provided, be used in computer-assisted proofs \cite{nigam2025intersection,jakobson2005large,nigam2020proof}. 

SSPs are typically approached via a ``discretize-then-solve'' framework, where the operator is approximated by a finite matrix and a (potentially generalized) eigenvalue problem is solved. Discretization methods include embedding techniques \cite{macdonald2011solving,brandman2008level,lee2024simple}, meshless methods \cite{wu2022surface,venn2024meshfree}, finite elements \cite{lu2021geometrically,reuter2006laplace}, and finite differences \cite{lee2024simple,glowinski2008computing}. However, discretization can introduce issues like spectral pollution (spurious eigenvalues), spectral invisibility (missing true eigenvalues), and failure to capture continuous spectra  \cite{colbrook2025avoiding,lewin2010spectral,boulton2012generalized,boulton2016spectral}. These problems often persist even with finer discretizations and are further exacerbated in SSPs due to the added complexity of discretizing the surface.

In this work, we adopt a ``solve-then-discretize'' framework for SSPs, aiming to delay discretization as much as possible.
Rather than discretizing the operator itself, we discretize only the functions it acts on. This is advantageous, as controlling errors in function discretization is generally more tractable than in operator discretization. This methodology has been successfully applied to infinite-dimensional eigenvalue problems (not on surfaces) \cite{contFEAST,colbrook2025avoiding}, and also underpins work on semigroups \cite{colbrook2022computing}, spectral measures \cite{webb2021spectra,specsolve}, spectra \cite{colbrook2025avoiding,colbrook2024computation,sci_l2N}, continuous Krylov methods \cite{gilles2019continuous,olver2009gmres,vainikko2004gmres}, and infinite-dimensional QL/QR algorithms \cite{colbrook2019infinite,webb2017isospectral}. 
In this paper, we extend three such algorithms to SSPs: \texttt{contFEAST} \cite{contFEAST}, which computes eigenvalues of differential operators, \texttt{infBEYN} \cite{colbrook2025avoiding}, which solves nonlinear eigenvalue problems, and \texttt{SpecSolve} \cite{specsolve}, which computes spectral measures of self-adjoint operators.

Our first key tool is the operator resolvent \cite[Sec.~III.6]{kato_1980}:
\begin{equation*}\setlength\abovedisplayskip{0pt}\setlength\belowdisplayskip{6pt}\setlength\abovedisplayskip{6pt}\setlength\belowdisplayskip{6pt}R(z,\mathcal{L})=(z\mathcal{I}-\mathcal{L})^{-1},\quad z\notin\spec(\mathcal{L}).\end{equation*}
Contour integrals involving the resolvent yield projections onto eigenspaces, enabling accurate eigenvalue computation, and convolution with high-order rational kernels allows computation of spectral measures. We only use applications of the resolvent to functions, so we can avoid the woes of discretization described above and choose the discretization size adaptively depending on $z$ to maintain high accuracy. Additionally, the resolvent norm provides a posteriori error bounds, essential for verifying results.

Our second key tool is the high-order fast direct surface PDE solver in \cite{fortunato2024high}, which we use to compute the action of the resolvent by solving shifted linear systems:
\begin{equation}\setlength\abovedisplayskip{6pt}\setlength\belowdisplayskip{6pt}\setlength\abovedisplayskip{6pt}\setlength\belowdisplayskip{6pt}\label{eqn:shifted_lin}(z\mathcal{I}-\mathcal{L})u=f,\quad u\in\mathcal{D}(\mathcal{L}).\end{equation}
This is the only step where discretization occurs. The solver achieves high-order accuracy, is structured to handle multiple right-hand sides efficiently, and supports a broad class of PDEs and surfaces. Moreover, in all of our algorithms the shifted systems for different $z$ can be solved in parallel, enabling effective use of modern computing architectures.

Combining these key tools, we introduce \texttt{SurfSpec}, a suite of routines for computing spectral properties of differential operators on surfaces. To the best of our knowledge, this package is the first to enable error-controlled computation of high-frequency eigenmodes for general PDEs on generic smooth surfaces in $\mathbb{R}^3$, as well as spectral measures and solutions to nonlinear eigenvalue problems. All algorithms come with convergence guarantees obtained by combining the results in the papers \cite{contFEAST,colbrook2025avoiding,specsolve} with the convergence of the high-order solver \cite{fortunato2024high}. Code for \texttt{SurfSpec} and all examples in this paper is available at \url{https://github.com/GustavConradie1/SurfSpec}.

\subsection{Organization of the paper}

In \cref{surfacefun}, we introduce the \texttt{surfacefun} PDE solver and outline its key features. \Cref{sec:diff_eval_prob} covers surface differential eigenvalue problems, detailing how resolvent-based contour integral methods can be used to compute eigenvalues with error control, and how this framework extends to nonlinear eigenvalue problems. In \cref{sec:surfspecmeas}, we define spectral measures of linear operators and show how to compute them using resolvent-based convolution methods, including on unbounded surfaces. Each section includes numerical examples demonstrating the accuracy and practicality of the proposed methods.

\section{Solving surface PDEs}\label{surfacefun}

Let $\surf$ be a smooth manifold embedded in $\mathbb{R}^3$ that is locally parametrized by a smooth map $\vec{x}(\xi,\eta):\mathbb{R}^2\rightarrow\mathbb{R}^3$. We begin by briefly discussing how to define derivatives, and hence differential equations, on $\surf$. The metric tensor $g$ on $\surf$ encapsulates its geometry and is given by
\begin{equation*}\setlength\abovedisplayskip{6pt}\setlength\belowdisplayskip{6pt}
g = \begin{pmatrix}
g_{\xi\xi}  & g_{\xi\eta} \\
g_{\eta\xi} & g_{\eta\eta}
\end{pmatrix}
=
\begin{pmatrix}
\vec{x}_\xi  \cdot \vec{x}_\xi & \vec{x}_\xi  \cdot \vec{x}_\eta \\
\vec{x}_\eta \cdot \vec{x}_\xi & \vec{x}_\eta \cdot \vec{x}_\eta
\end{pmatrix},
\end{equation*}	
where $\vec{x}_\xi$ and $\vec{x}_\eta$ represent the partial derivatives of $\vec{x}$ with respect to $\xi$ and $\eta$ respectively. 
The corresponding inverse metric tensor is
\begin{equation*}\setlength\abovedisplayskip{6pt}\setlength\belowdisplayskip{6pt}
g^{-1} = \begin{pmatrix}
g^{\xi\xi}  & g^{\xi\eta} \\
g^{\eta\xi} & g^{\eta\eta}
\end{pmatrix}
=
\frac{1}{|g|}
\begin{pmatrix}
g_{\eta\eta} & -g_{\xi\eta} \\
-g_{\eta\xi}  &  g_{\xi\xi}
\end{pmatrix},
\end{equation*}
where $|g|=\det g$. The components of the inverse metric tensor are used to define tangential derivative operators on the surface. In particular, if we define
\begin{equation*}\setlength\abovedisplayskip{6pt}\setlength\belowdisplayskip{6pt}
\xi_\vec{x}  = g^{\xi\xi}  \vec{x}_\xi + g^{\eta\xi}  \vec{x}_\eta, \quad
\eta_\vec{x} = g^{\xi\eta} \vec{x}_\xi + g^{\eta\eta} \vec{x}_\eta,
\end{equation*}
then the tangential derivatives in the Cartesian coordinate directions of a scalar field $u = u(\xi,\eta) = u(\vec{x}(\xi,\eta))$ on $\surf$ are given by~\cite{Frankel2012}
\begin{equation*}\setlength\abovedisplayskip{6pt}\setlength\belowdisplayskip{6pt}
\partial_x^\surf u = \xi_x\tfrac{\partial u}{\partial\xi} + \eta_x\tfrac{\partial u}{\partial \eta}, \quad
\partial_y^\surf u = \xi_y\tfrac{\partial u}{\partial\xi} + \eta_y\tfrac{\partial u}{\partial \eta}, \quad
\partial_z^\surf u = \xi_z\tfrac{\partial u}{\partial\xi} + \eta_z\tfrac{\partial u}{\partial \eta}.
\end{equation*}
Familiar quantities can be defined in terms of these tangential derivative operators. For example, the gradient of a scalar field $u$ defined on $\surf$ is a tangential vector field $\nabla_\surf u = (\partial_x^\surf u, \partial_y^\surf u, \partial_z^\surf u)$, and the divergence of a vector field $\vec{u} = (u_1, u_2, u_3)$ tangent to $\surf$ is a scalar field $\nabla_\surf \cdot u = \partial_x^\surf u_1 + \partial_y^\surf u_2 + \partial_z^\surf u_3$. The surface Laplacian $\Delta_\surf$, the negative of which is known as the Laplace--Beltrami operator, is
\begin{equation*}\setlength\abovedisplayskip{6pt}\setlength\belowdisplayskip{6pt}
\Delta_\surf u = \nabla_\surf \cdot \nabla_\surf u = \partial_x^\surf \partial_x^\surf u + \partial_y^\surf \partial_y^\surf u + \partial_z^\surf \partial_z^\surf u.
\end{equation*}
In this work, we consider general second-order, variable-coefficient elliptic partial differential operators $\mathcal{L}_\surf$ on $\surf$, which take the form
\begin{equation*}\setlength\abovedisplayskip{6pt}\setlength\belowdisplayskip{6pt}
\mathcal{L}_\surf u(\vec{x}) = \sum_{i=1}^3 \sum_{j=i}^3 a_{ij}(\vec{x}) \, \partial_i^\surf \partial_j^\surf u(\vec{x}) + \sum_{i=1}^3 b_i(\vec{x}) \,\partial_i^\surf u(\vec{x}) + c(\vec{x}) u(\vec{x}),
\end{equation*}
for smooth coefficients \(a_{ij}\), \(b_i\), and \(c\) on \(M\), and $\partial_1^\surf$, $\partial_2^\surf$ and $\partial_3^\surf$ identified with $\partial_x^\surf$, $\partial_y^\surf$ and $\partial_z^\surf$ respectively; the methods used can be extended straightforwardly to higher-order PDEs. The corresponding surface PDE on $\surf$ then takes the form
\begin{equation}\setlength\abovedisplayskip{6pt}\setlength\belowdisplayskip{6pt}\label{eq:surface_pde}
\mathcal{L}_\surf u(\vec{x}) = f(\vec{x}), \qquad \vec{x}\in\surf,
\end{equation}
where $f(\vec{x})$ is a smooth function defined on $\surf$. If $M$ is a surface with boundary (an open, as opposed to closed, surface), for $u(\vec{x})$ to be uniquely determined \cref{eq:surface_pde} may be accompanied by boundary conditions, e.g., $u(\vec{x}) = g(\vec{x})$ for $\vec{x} \in \partial\surf$ and some function $g(\vec{x})$.
For notational convenience in the rest of the paper, we drop the subscripts and superscripts indicating the surface.

\begin{figure}[t]
	\centering
	\includegraphics[width=0.9\linewidth]{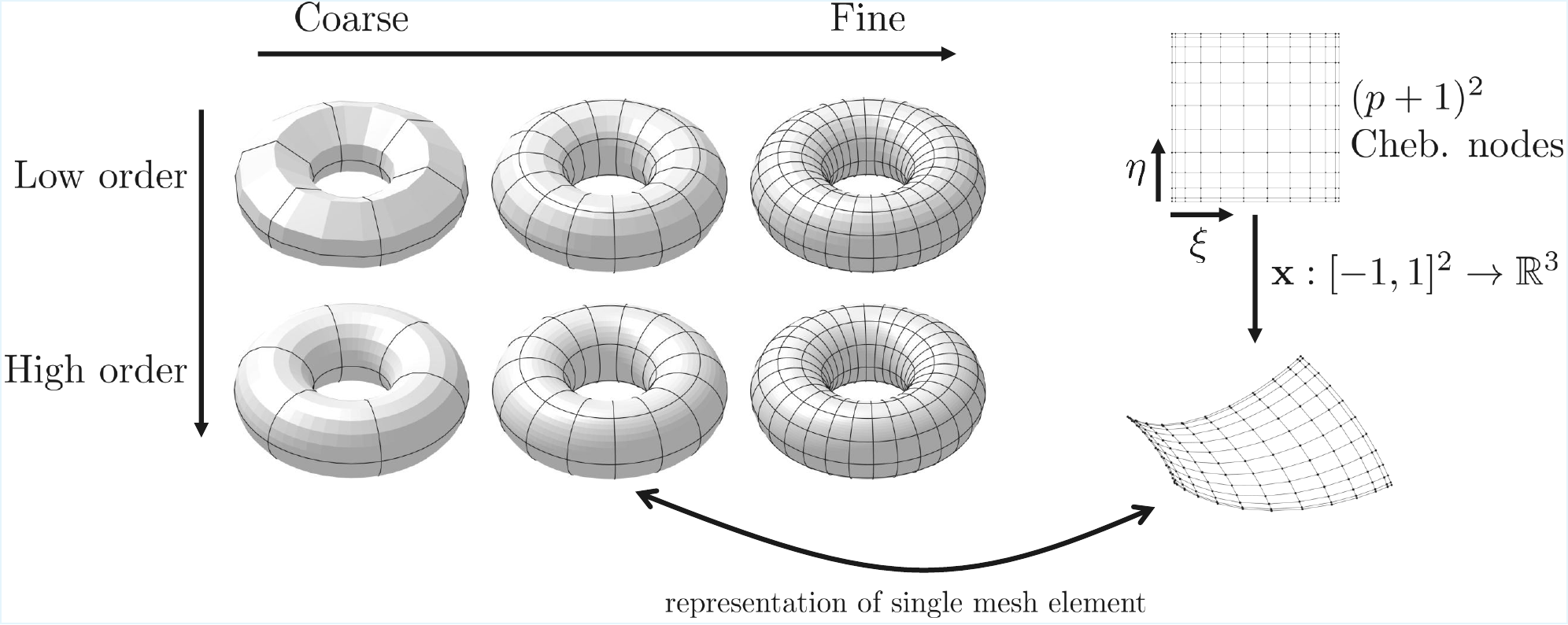}
    \caption{Left: A surface mesh may be refined along two axes. Right: On each element, the surface is sampled at tensor-product Chebyshev nodes of the second kind.}
	\label{fig:mesh_refinement}
\end{figure}

To solve problems of the form \cref{eq:surface_pde}, we use a recently developed fast direct solver for surface PDEs based on spectral collocation on unstructured quadrilateral surface meshes. We give a rough overview of the scheme; for more details, see~\cite{fortunato2024high}. Let $\{\mathcal{E}_k\}_{k=1}^N$ denote a mesh of $N$ elements each discretized to order $p$. In particular, each element $\mathcal{E}_k$ is represented by sampling $\surf$ at tensor-product second-kind Chebyshev nodes of order $p$ on $[-1,1]^2$, resulting in $(p+1)^2$ nodes per element. The function $\mathbf{x}:[-1,1]^2\rightarrow\mathbb{R}^3$ which locally parameterizes the surface on each mesh element is then constructed by Chebyshev interpolation. Discretization can be refined either by increasing $N$ (mesh refinement) or $p$ (higher local resolution) (see~\cref{fig:mesh_refinement}). Local metric tensors are computed via spectral differentiation, and discrete surface differential operators are assembled using matrix algebra to yield a spectral collocation scheme. Functions such as the variable coefficients are sampled at the same Chebyshev nodes on each element.

We split the solution to \cref{eq:surface_pde} into homogeneous and particular components, and discretize these problems separately on each element, resulting in $2N$ $(p+1)^2 \times (p+1)^2$ linear systems. The solution over the whole surface is then computed in three steps:%
\vspace{1em}
\begin{enumerate}
\item \textbf{Local factorization.} On each mesh element, we perform a local solve to construct solution operators that map Dirichlet boundary data to the homogeneous local solution, along with Dirichlet-to-Neumann maps that yield the corresponding Neumann data. We also compute discrete particular solutions and their associated fluxes.%
\vspace{1em}
\item \textbf{Global factorization.} In an upward pass, adjacent elements are ``glued'' to form solution operators, particular solutions, Dirichlet-to-Neumann maps and particular fluxes on the merged parent elements. Each merge enforces continuity of the solution and its binormal derivative along the interface using the Dirichlet-to-Neumann maps and particular fluxes from the child elements.%
\vspace{1em}
\item \textbf{Solve.} In a downward pass, we begin at the top-level element and use its solution operator, Dirichlet data, and particular fluxes to reconstruct the interface solutions on its child elements. This process continues recursively, reconstructing solutions on all merged interfaces. Once all boundary data is computed for each mesh element, the local solution operators and particular solutions from the factorization step are used to obtain the full solution.%
\end{enumerate}
\vspace{1em}

This approach is fast: for fixed $p$, the local factorization stage takes $\mathcal{O}(N)$ operations, the global factorization stage takes $\mathcal{O}(N^{3/2})$ operations, and the solve stage takes $\mathcal{O}(N\log N)$ operations. Additionally, any boundary data passed to the solver is used only in the solve stage, so we can solve the same PDE with different boundary data without recomputing any factorizations. Updating {the right-hand side $f$} is also efficient, requiring only $\mathcal{O}(N\log N)$ steps, as we do not need to recompute solution operators and Dirichlet-to-Neumann maps---just particular solutions and fluxes. The memory complexity also scales at rate $\mathcal{O}(N\log N)$.

We use the open-source MATLAB package \texttt{surfacefun}~\cite{surfacefun}, which implements the above solver and provides convenient high-order tools for computing with functions on surfaces. In what follows, it forms a key building block for computing spectral properties of surface operators.

\section{Differential eigenvalue problems}\label{sec:diff_eval_prob}

This section outlines our approach to computing discrete spectra of differential operators on surfaces. We begin with linear eigenvalue problems and then extend to nonlinear cases, providing error control for both regimes. The core strategy is to evaluate contour integrals involving the resolvent. 

\subsection{SurfFEAST for linear differential eigenvalue problems}\label{sec:surfFEAST_linear}

The work \cite{contFEAST} provides an analogue of the FEAST algorithm \cite{OriginalFEAST} for differential operators on $[-1,1]$. Motivated by this, we develop an algorithm for solving eigenproblems involving differential operators $\mathcal{L}$ on general surfaces $M$, which we call \texttt{SurfFEAST}. We assume that the spectrum of $\mathcal{L}$ is discrete in a domain $\Omega$, so any eigenvalue of $\mathcal{L}$ in $\Omega$ is isolated from the rest of the spectrum, and that the boundary of $\Omega$ does not intersect $\spec(\mathcal{L})$. The key theoretical tool is the Riesz projector, which is defined for the domain $\Omega$ by
\[\setlength\abovedisplayskip{6pt}\setlength\belowdisplayskip{6pt}
\mathcal{P}_\Omega=\cfrac{1}{2\pi i}\int_{\partial\Omega}(z\mathcal{I}-\mathcal{L})^{-1}\dd z. 
\]
Under our assumptions, this operator is the projection onto the finite-dimensional invariant subspace spanned by eigenfunctions with eigenvalues in $\Omega$ \cite[Sec.~III.6]{kato_1980}. By applying it to suitable test functions, we can recover the corresponding eigenfunctions and eigenvalues inside $\Omega$. Our algorithm for this is summarized in \cref{surfFEAST}; for the remainder of this section we elucidate the key features of the algorithm.

\begin{algorithm}[t]
\caption{\texttt{SurfFEAST}: Computing eigenvalues of a surface differential operator.}\label{surfFEAST}
\textbf{Input:} Operator $\mathcal{L}$, quadrature weights $\{w_j\}_{j=1}^\ell$ and nodes $\{z_j\}_{j=1}^\ell$, test functions $F=[f_1\cdots f_{m+m_0}]$, $\epsilon_{\res},\epsilon_{\mathrm{tol}}>0$.
\begin{algorithmic}[1]
\REPEAT
\STATE{Solve $(z_k\mathcal{I}-\mathcal{L})G_k=F$ for $k=1,\ldots,\ell$ with relevant boundary conditions.}
\STATE{Compute QR decomposition of  $(2\pi i)^{-1}\sum_{k=1}^\ell w_k\hat{G}_k$ for orthonormal basis $Q$.}
\STATE{Compute a truncated SVD of $Q\approx \mathcal{U}^*\Sigma\mathcal{V}$, excluding singular values $<\epsilon_{\mathrm{tol}}$.}
\STATE{Compute $L=Q^*\mathcal{L}Q$, solve $LX=X\Lambda$ for $X$ and $\Lambda$, and update $F=QX$.}
\UNTIL{$\|\mathcal{L}F-F\Lambda\|\leq\epsilon_{\res}\|\Lambda\|$.}
\end{algorithmic}
\textbf{Output:} Eigenvalues $\Lambda$, eigenfunctions $U=QX$.
\end{algorithm}

Assume that the sum of the algebraic multiplicities of eigenvalues inside the contour is $m$. We select $m + m_0$ test functions $f_1, \ldots, f_{m + m_0}$, where $m_0$ is an oversampling parameter, and let $F$ be the quasimatrix with columns $f_1,\ldots,f_{m+m_0}$. We construct $F$ using band-limited random functions\footnote{{A band-limited periodic random function on a cuboid is given by a truncated Fourier series with coefficients drawn from a Gaussian distribution; the nonperiodic analog (which we use) is constructed by restricting a band-limited periodic random function defined on a large cuboid containing $M$ to the surface $M$ \cite{randfun,contFEAST}.}} as in \cite{randfun,contFEAST}. Using more test functions than eigenvalues ensures all eigenvalues (counting multiplicity) are captured, improves convergence and helps avoid contamination from nearby eigenvalues just outside the contour \cite{nakatsukasa2020sharp,peter2014feast}. In practice, $m$ is often unknown but can be estimated via asymptotic formulas like Weyl's law (see \cref{laplace}) or stochastic trace estimation \cite{meyer2021hutch++,persson2022improved,zvonek2025conthutch++}; overestimating is preferable. Any excess rank can later be removed via truncated SVD as discussed later in this section; see \cite{quasiQR} for details on quasimatrix SVDs.

Given the quasimatrix $F$, we want to compute $\mathcal{P}_{\Omega}F$.
We approximate the resulting integral via a quadrature rule as follows:
\begin{equation*}\setlength\abovedisplayskip{6pt}\setlength\belowdisplayskip{6pt}
\label{eqn:approximate spectral projector}
\mathcal{P}_\Omega F=\left(\cfrac{1}{2\pi i}\int_{\partial\Omega}(z\mathcal{I}-\mathcal{L})^{-1}\dd z\right)F\approx \hat{\mathcal P}_{\Omega}F=\frac{1}{2\pi i}\sum_{k=1}^\ell w_k(z_k\mathcal{I}-\mathcal{L})^{-1}F.
\end{equation*}
Here $z_1,\ldots,z_\ell$ and $w_1,\ldots,w_\ell$ are a set of quadrature nodes and weights, respectively; for the experiments in this paper we use the exponentially convergent trapezoidal rule \cite{traprule}. For both empirical and theoretical reasons \cite{colbrook2025avoiding}, we use medium-sized contours to locate approximate locations of the eigenvalues before using small circular contours to ``zoom in'' on a particular eigenvalue and compute it to higher accuracy. If $\mathcal{L}$ is self-adjoint, then all of its eigenvalues are real and so we use circular or elliptical contours centered on the real axis, with quadrature nodes taken to lie off the real axis (see \eqref{eqn:erb_inverseproblem}). If \(F\) is real-valued, then\[\setlength\abovedisplayskip{6pt}\setlength\belowdisplayskip{6pt}(z\mathcal{I}-\mathcal{L})^{-1}F+(\overline{z}\mathcal{I}-\mathcal{L})^{-1}F=2\mathrm{Re}((z\mathcal{I}-\mathcal{L})^{-1})F\]
and so we can halve the number of computations required; for example, all Laplace--Beltrami eigenfunctions can be taken to be real-valued and so we may take $F$ real.

Construction of $\hat{\mathcal P}_\Omega F$ then requires the solution to $\ell$ shifted linear surface PDEs, each with $m+m_0$ right-hand sides, i.e.,
\begin{equation}\setlength\abovedisplayskip{6pt}\setlength\belowdisplayskip{6pt}
\label{eqn:shifted linear ODEs}
(z_k\mathcal{I}-\mathcal{L})g_{i,k}=f_i,\qquad  i=1,\ldots,m+m_0, \qquad k=1,\ldots,\ell.
\end{equation}
The solver discussed in \cref{surfacefun} has the advantage that when solving the PDE for $m+m_0$ different {right-hand sides $f$}, we do not need to recompute solution operators and Dirichlet-to-Neumann maps~\cite{fortunato2024high}. It also enables a high degree of parallelization: we may solve the shifted linear systems in parallel for each shift, and for each shifted linear system the local factorization step (see \cref{surfacefun}) can be carried out in parallel on each mesh element. The solver produces approximate solutions $\hat{g}_{i,k}$ to \cref{eqn:shifted linear ODEs}. If the quasimatrix with columns $\hat{g}_{1,k},\dots,\hat{g}_{m+m_0,k}$ is denoted by $\hat{G}_k$, then  
\begin{equation*}\setlength\abovedisplayskip{6pt}\setlength\belowdisplayskip{6pt}
\hat{\mathcal P}_\Omega F=\frac{1}{2\pi i}\sum_{k=1}^\ell w_k\hat{G}_k.
\end{equation*}
The columns of  $V=\mathcal{P}_\Omega F$ span the subspace of eigenfunctions whose eigenvalues lie inside $\Omega$, and $\hat V=\hat{\mathcal P}_\Omega F$ approximates this. We compute an orthonormal basis $\smash{\{q_j\}_{j=1}^m}$ of $\hat V$ with respect to the inner product of $\mathcal{H}$ via a QR factorization \cite{quasiQR}; we let $Q$ be the quasimatrix with columns $q_j$. The choice of band-limited periodic functions usually ensures that $\hat{\mathcal P}_\mathcal{V}F$ has linearly independent and well-conditioned columns when $\mathcal{L}$ is normal \cite{contFEAST}.

Next, we apply a truncated SVD to $Q$ to compensate for overestimating $m$. {If $m_0>0$, the span of the columns of $\hat{V}$ will include approximate eigenfunctions corresponding to eigenvalues which lie outside $\Omega$, as well as (rarely) spurious eigenvalues, and the truncated SVD will remove some of these depending on the tolerance $\epsilon_{\mathrm{tol}}$; typically, singular values drop sharply after those corresponding to true eigenvalues inside the contour. The larger $m_0$ is, the more eigenvalues outside the contour will be picked up, beginning with the ones closest to the contour. Eigenvalues outside the contour will typically be approximated inaccurately, with the accuracy decreasing the further outside the contour they are, so it is desirable to remove them. }

Finally, we take a Rayleigh-Ritz projection $L=Q^*\mathcal{L}Q\in\mathbb{C}^{n\times n}$. The elements of $L$ are given by the inner products $L_{ij}=\langle q_i,\mathcal{L}q_j\rangle$. This is done in \texttt{surfacefun} for the $L^2(M)$ norm using a quadrature rule based on Chebyshev points and weights and using the square root of the Jacobian. Ignoring numerical error, the eigenvalues $\Lambda$ of $L$ are exactly the {eigenvalues of $\mathcal{L}$ that lie inside the given contour} \cite{contFEAST}; hence, we have reduced the infinite-dimensional eigenvalue problem to a finite-dimensional one without directly discretizing $\mathcal{L}$. The corresponding eigenfunctions are given by the columns of $U=QX$, where $X$ is the matrix of eigenfunctions of $L$ given by $LX=X\Lambda$. 

To improve accuracy, we can iterate the algorithm, updating the test functions according to $F=QX$ and repeating the previous steps. We terminate when $\|\mathcal{L}F-F\Lambda\|\leq \epsilon_{\res}\|\Lambda\|$; motivation for this termination condition can be found in the next section. If the contour encircles only a single eigenvalue of multiplicity $1$, no benefit is obtained from iterations beyond the first; otherwise, taking two iterations is typically both necessary and sufficient for accurate computations \cite{horning2022twice}. These iterations can be interpreted in terms of rational subspace iteration \cite{peter2014feast}.

\subsection{Error control}
\label{errctrl}

Due to the difficulties of computing spectral properties of infinite-dimensional operators numerically, verification is crucial. For a normal operator $\mathcal{L}$, it is known that distance to the spectrum is exactly the inverse of the resolvent norm \cite{davies_2007}, i.e.,
\begin{equation}\label{eqn:specdist}\setlength\abovedisplayskip{6pt}\setlength\belowdisplayskip{6pt}
\dist(z,\spec(\mathcal{L}))=\|(z\mathcal{I}-\mathcal{L})^{-1}\|^{-1},
\end{equation}
where $\|(z\mathcal{I}-\mathcal{L})^{-1}\|^{-1}=0$ if $z\in\spec(\mathcal{L})$. Moreover, one can characterize the resolvent norm in terms of the injection modulus \cite{sci_l2N}
\[\setlength\abovedisplayskip{6pt}\setlength\belowdisplayskip{6pt}
\|(z\mathcal{I}-\mathcal{L})^{-1}\|^{-1}=\sigma_{\inf}(z\mathcal{I}-\mathcal{L}),
\quad
\sigma_{\inf}(\mathcal{L})=\inf\{\|\mathcal{L}f\|:f\in\mathcal{D}(\mathcal{L}),\|f\|=1\}.
\]
Combining these, we get an a posteriori estimate for the distance to the true spectrum of a proposed eigenvalue $\lambda$: if the corresponding normalized proposed eigenvector is $f$, then the distance of $\lambda$ to the spectrum is bounded above by
\begin{equation}\setlength\abovedisplayskip{6pt}\setlength\belowdisplayskip{6pt}\label{eqn:residual_defn}
\dist(\lambda,\spec(\mathcal{L}))\leq \|(\lambda\mathcal{I}-\mathcal{L})f\|.
\end{equation}
We refer to $\|(\lambda\mathcal{I}-\mathcal{L})f\|$ as the residual. Using this error bound, we can discard spurious eigenvalues and avoid spectral pollution. This also provides justification for the choice of termination condition in \cref{surfFEAST}. 

Furthermore, if $\mathcal{L}$ is self-adjoint we can use similar arguments to verify the solution to the shifted linear systems in \eqref{eqn:shifted linear ODEs}. {In particular, since $\spec(\mathcal{L})\subset \mathbb{R}$ it follows by \cref{eqn:specdist} that $\|(z\mathcal{I}-\mathcal{L})^{-1}\|=1/\dist(z,\spec(\mathcal{L}))\leq1/|\mathrm{Im}(z)|$ for all $z$ with $\mathrm{Im}(z)\neq 0$, and so}
\begin{equation}\setlength\abovedisplayskip{6pt}\setlength\belowdisplayskip{6pt}\label{eqn:erb_inverseproblem}
    \|\hat{g}_{i,k}-(z_k\mathcal{I}-\mathcal{L})^{-1}f_i\|\leq  \cfrac{\|(z_k\mathcal{I}-\mathcal{L})\hat{g}_{i,k}-f_i\|}{|\mathrm{Im}(z_k)|}.
\end{equation}
This allows for adaptive solution of the shifted linear systems to improve accuracy and decrease the eigenvalue residual.

For a nonnormal operator, obtaining bounds for the distance to the spectrum in terms of the resolvent is trickier, though possible \cite[Sec.~2.2.2]{colbrook2024computation}.
Moreover, bounds on the resolvent norm can be interpreted in terms of pseudospectral set inclusions \cite{trefethen2020spectra}, which still provide verification and help avoid spectral pollution. 
Given $\epsilon>0$, we define the $\epsilon$-pseudospectrum of $\mathcal{L}$ by
\[\label{eqn:pspec_defn}\setlength\abovedisplayskip{6pt}\setlength\belowdisplayskip{6pt}\spec_{\epsilon}(\mathcal{L})=\mathrm{Cl}\left(\{z\in\mathbb{C}:\|(z\mathcal{I}-\mathcal{L})^{-1}\|^{-1}<\epsilon\}\right),\]
where $\mathrm{Cl}$ denotes the closure of the set. The pseudospectrum captures how sensitive the spectrum is to perturbations, as shown by the following relation \cite[Thm.~4.3]{trefethen2020spectra}:
\[\setlength\abovedisplayskip{6pt}\setlength\belowdisplayskip{6pt}\spec_{\epsilon}(\mathcal{L})=\mathrm{Cl}\left(\bigcup_{\|{E}\|<\epsilon}\spec(\mathcal{L}+{E})\right).\]
The pseudospectrum converges to the spectrum as $\epsilon\rightarrow 0$ in a suitable topology which captures spectral pollution and spectral invisibility \cite{sci_l2N}. {Theoretical bounds on how numerical errors affect the results of contour methods for differential eigenvalue problems, and corresponding convergence results, are typically stated in terms of pseudospectral inclusions \cite{colbrook2025avoiding,contFEAST}.}

\subsection{Example: Laplace--Beltrami and Schr\"odinger operators}
\label{laplace}

\begin{figure}
    \centering
    \includegraphics[width=0.4\linewidth]{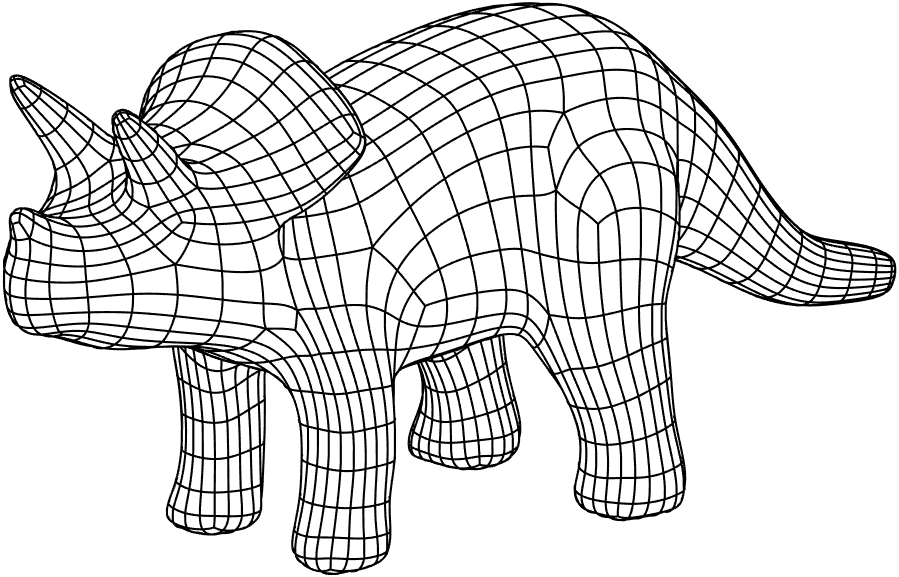}
    \caption{The triceratops mesh used in the numerical experiments in \cref{laplace}.}
    \label{fig:triceratops mesh}
\end{figure}

As a first example, we consider the Laplace--Beltrami operator $-\Delta: \mathcal{D}(-\Delta)\rightarrow L^2(\surf)$ with Dirichlet boundary conditions \cite{craioveanu_2001}. It has a discrete spectrum consisting of eigenvalues $0\leq\lambda_1\leq\lambda_2\leq\cdots\rightarrow\infty$, with each eigenvalue having finite multiplicity. For closed surfaces $\lambda_1=0$, while for Dirichlet problems on open surfaces $\lambda_1>0$ \cite{chavel_1984}. We define the eigenvalue counting function $N(\lambda)$ to be the number of eigenvalues less than or equal to $\lambda$. Weyl's law \cite{weyl} for surfaces embedded in $\mathbb{R}^3$ states that 
\begin{equation}\setlength\abovedisplayskip{6pt}\setlength\belowdisplayskip{6pt}\label{weylone}
N(\lambda)= \frac{1}{4\pi}\area(\surf)\lambda+o(\lambda)\quad\text{as}\quad\lambda\rightarrow\infty.
\end{equation}
We use this as a point of comparison for our numerical computations, and to estimate $m$. We also compare to the additional conjectured terms in Weyl's law \cite{weylmore}:
\begin{equation}\setlength\abovedisplayskip{6pt}\setlength\belowdisplayskip{6pt}\label{weyltwo}
    N(\lambda)= \frac{1}{4\pi}\area(\surf)\lambda- \cfrac{1}{4\pi}\length (\partial \surf)\lambda^{1/2}+o(\lambda^{1/2})\quad\text{as}\quad\lambda\rightarrow\infty.
\end{equation}
Note that for closed surfaces, the second term is $0$. A summary of progress towards proving this conjecture can be found in \cite{100years}.

As a second example, we compute eigenvalues for surface Schr\"odinger operators
$
\mathcal{L}=-\Delta+V(\mathbf{x})$. On compact manifolds, for potentials satisfying a local Lipschitz condition, \cref{weylone} is also satisfied \cite{schrweyl}. We choose a harmonic oscillator potential $V(\mathbf{x})=\|\mathbf{x}\|^2$. Though it looks simple, for complicated manifolds the potential is far from polynomial in the coordinate charts.

\subsubsection{A worked example: a triceratops}\label{sec:triceratops}

We begin with an in-depth example to analyze convergence of \cref{surfFEAST} as different parameters are varied. We apply \cref{surfFEAST} to the Laplace--Beltrami operator and the Schr\"odinger operator with harmonic oscillator potential on the triceratops mesh shown in \cref{fig:triceratops mesh}. The original mesh was obtained from \url{https://quaternius.com}, and has been smoothed and simplified using Rhino 8; the final mesh has 1781 faces and is of order $12$. The Python script \texttt{surfacepts} converts the SubD mesh exported by Rhino to a \texttt{.csv} file that can be imported and used in the \texttt{surfacefun} code.

We compute the first 1000 eigenvalues of both the Laplace--Beltrami operator and Schr\"odinger operator using medium-sized elliptical contours and (over)estimating $m$ based on Weyl's law. {In particular, we used elliptical contours of semi-major axis $10.1$ and semi-minor axis $1$ which are centered at points $(10(2n+1),0)$ for $n\geq 0$, and took $m_0=30>m$. The overlap between different contours is to mitigate any decrease in accuracy near the edge of the contour, as the eigenvalues we are interested in are at least $0.1$ away from the edge of the contour and we can compare results from overlapping contours to ensure that no eigenvalues are missed. Further discussion on the choice of contour can be found in \cite{colbrook2025avoiding,colbrook2023computing,hale2008new,traprule} and the references contained within. We also use a truncated SVD with $\epsilon_{\mathrm{tol}}=10^{-8}$ to minimize the effects of overestimating $m$.} Due to the exponential convergence of the quadrature formula, relatively few quadrature nodes are needed. We solve only $l=15$ shifted linear systems for each contour integral (which, after exploiting symmetry, corresponds to $30$ quadrature nodes); {subsequent analysis will show that even fewer nodes would suffice}. 
 
\begin{figure}
    \centering
        \includegraphics[width=0.43\linewidth]{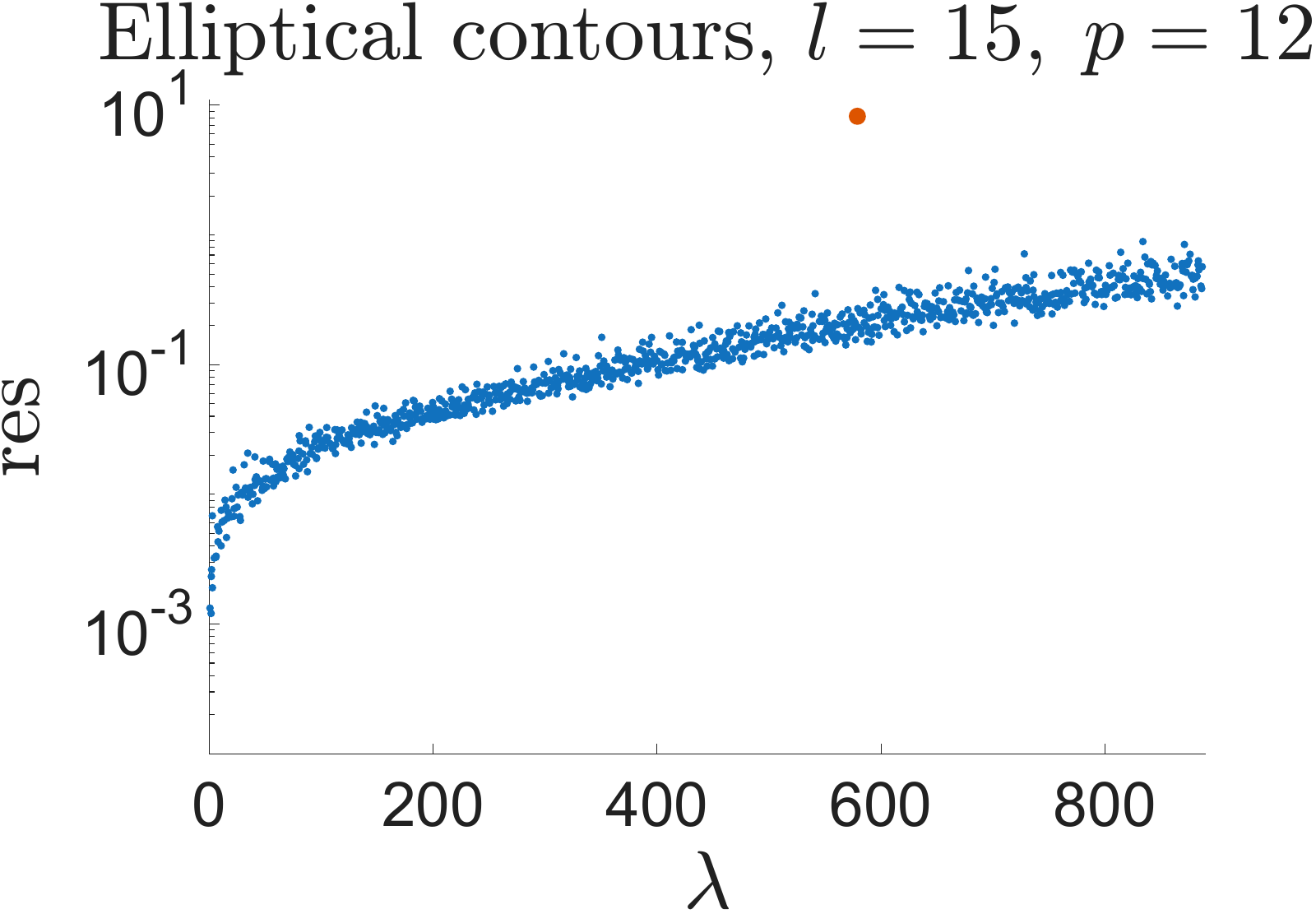}\hfill
        \includegraphics[width=0.43\linewidth]{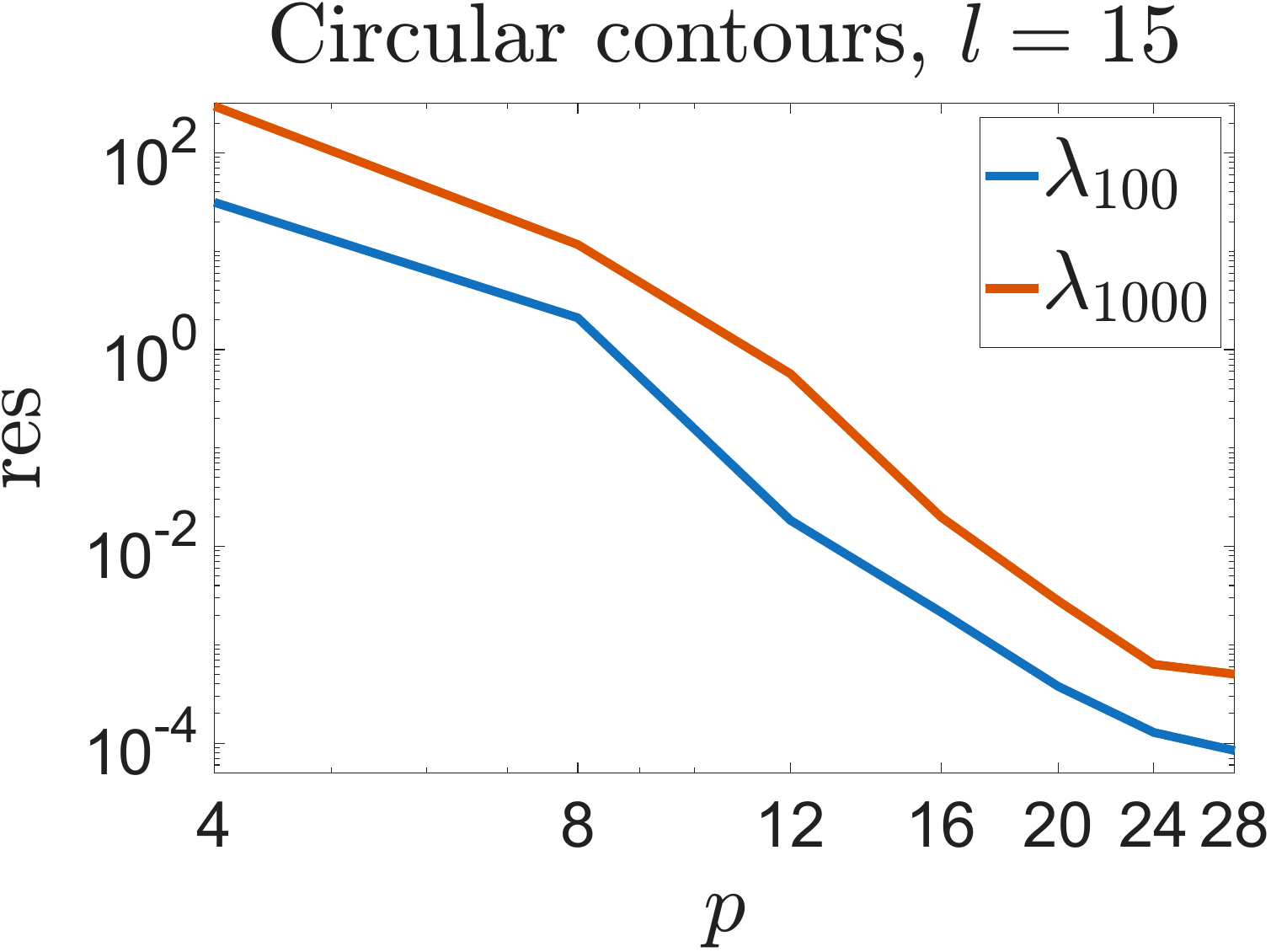}\\\vspace{0.3cm}
    \includegraphics[width=0.43\linewidth]{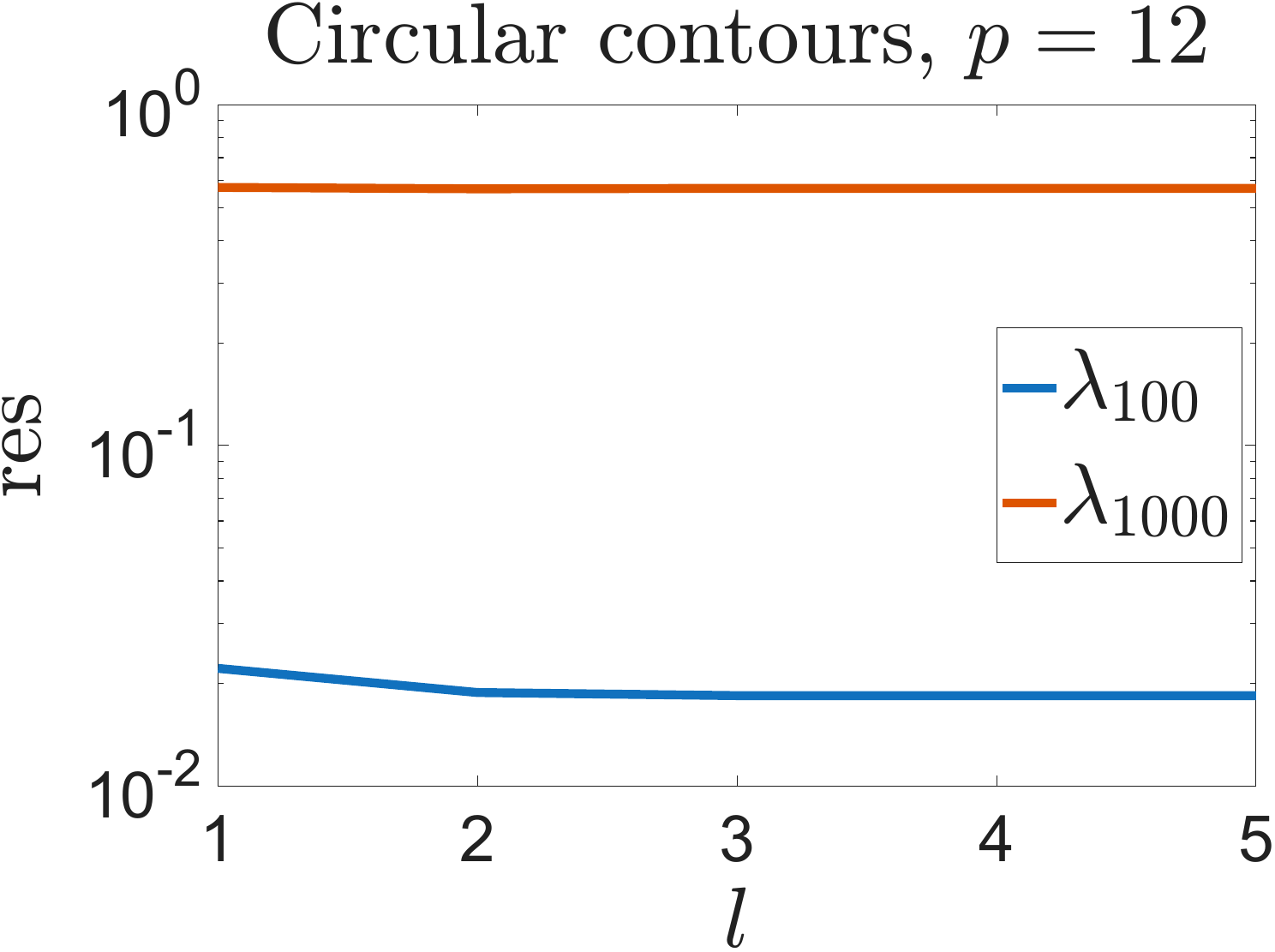}\hfill
    \includegraphics[width=0.43\linewidth]{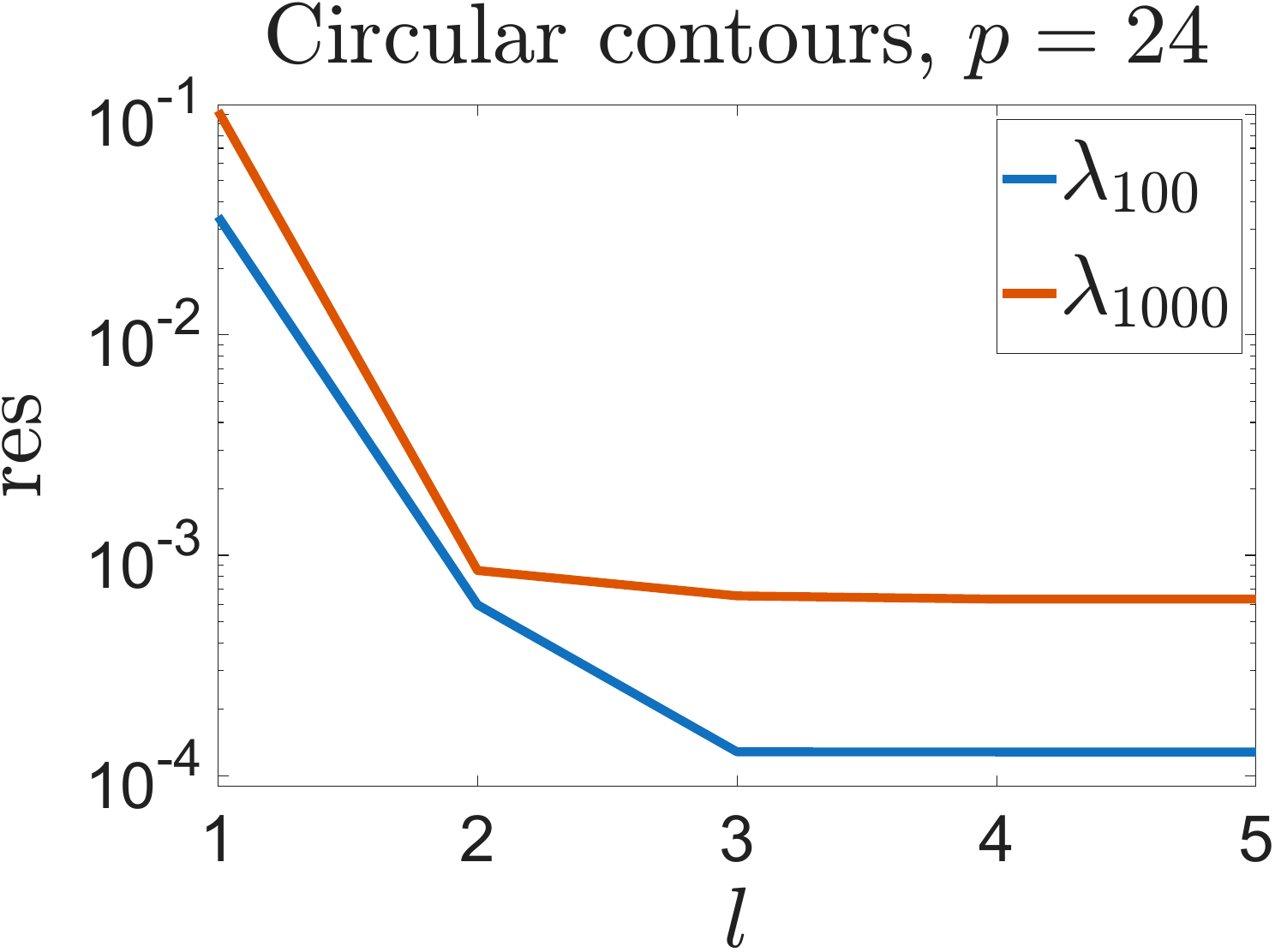}
    \caption{{Top left: The residuals of the first $1000$ eigenvalues of the Laplace-Beltrami operator on a triceratops mesh computed using elliptical contours; the only spurious eigenvalue is highlighted in red. Top right: The residuals of the $100\text{th}$ and $1000\text{th}$ Laplace-Beltrami eigenvalues on a triceratops mesh with varying order using circular contours with $l=15$. Bottom left: The residuals of the $100\text{th}$ and $1000\text{th}$ eigenvalues on a coarse triceratops mesh using circular contours with different numbers of quadrature nodes. Bottom left: The residuals of the $100\text{th}$ and $1000\text{th}$ eigenvalues on a fine triceratops mesh using circular contours with different numbers of quadrature nodes.}}
    \label{fig:schr_res_order}
\end{figure}

 We also compute the corresponding residuals following \eqref{eqn:residual_defn}. {The residual calculations are helpful to avoid any spurious eigenvalues resulting from not knowing $m$, inaccurate linear solves, the quadrature rule not having converged or contamination from eigenvalues just outside the contour: $1$ eigenvalue across the first $1000$ computed for the Laplace--Beltrami operator and $2$ eigenvalues for the Schr\"odinger operator were spurious (having residual much larger than $1$). The residuals of the eigenvalues computed for the Laplace--Beltrami operator using the described setup are shown on the top left of \cref{fig:schr_res_order}. We clearly identify the spurious eigenvalue (highlighted in red), }and also note that the higher frequency eigenmodes are harder to accurately resolve and hence have larger residuals. {If we use circular contours with the same centres and radius $10.1$ instead, the eigenvalues and residuals obtained are almost the same (the computed residuals typically agree to $8$ decimal places) but there are $3$ and $10$ spurious eigenvalues for the Laplace--Beltrami and Schr\"odinger operators respectively.}

{Unlike for many ``discrete-then-solve'' schemes, where spurious eigenvalues can persist even as large discretization sizes are taken \cite{lewin2010spectral,colbrook2025avoiding}, for our contour methods the spurious eigenvalues can be easily identified by residuals and disappear for even slight improvements of the discretization size. For example, for the spurious eigenvalue of the Laplace--Beltrami operator $\lambda_{\mathrm{spurious}}\approx 578.91$, if we increase $l$ from $15$ to $20$, reduce $m$ by the number of computed eigenvalues outside the contour (but still including the spurious eigenvalue) or equivalently decrease $\epsilon_{\mathrm{tol}}$ from $10^{-8}$ to $10^{-6}$, or adjust the centre of the contour to be closer to the spurious eigenvalue, then the algorithm no longer outputs the spurious eigenvalue. This is not the case for increases of $p$, indicating that the spurious eigenvalue arises due to inaccurate quadrature rules or overestimating $m$ rather than inaccuracy in the solves. The same observations hold for the two spurious eigenvalues outputted by the algorithm applied to the Schr\"odinger operator. By using residuals to identify spurious eigenvalues, we can adaptively tune parameters to obtain accurate results.}

We can also ``zoom in'' on the individual eigenvalues using a circular contour to compute them to higher accuracy. Restricting to the {$100\text{th}$ and $1000\text{th}$} eigenvalue, we vary the number of quadrature nodes and the refinement of the mesh to see how convergence is affected. {Throughout, residuals are smaller for the lower frequency eigenvalue.} The bottom panels of \cref{fig:schr_res_order} show convergence as we vary the number of quadrature nodes for surfaces of order $12$ and $24$ respectively; we see that the convergence is rapid, with residuals quickly stagnating even after only using $2\times 3$ quadrature nodes. {Slightly more quadrature nodes are needed for convergence for the higher-order surface; this is as error arising from coarseness of the surface no longer dominates the quadrature error.} As expected, refining the mesh by increasing the order $p$ does improve accuracy and decrease the residual; recall that the order is adjusted by varying the number of Chebyshev nodes on each patch. This is shown clearly on the right of \cref{fig:schr_res_order}{; we use $l=15$ nodes throughout to ensure that quadrature errors are not affecting the eigenvalues or residuals.} {The rate of decrease appears to stagnate as we increase $p$; this a restriction of $p$-refinement for the solver used, as the discretized geometries considered are only approximately smooth \cite{fortunato2024high}. In \cref{fig:nlep_convergence} we consider $h$-refinement and a linear (on a log-log plot) relation between the residual and number of mesh elements is observed.}

\subsubsection{Further examples}

To demonstrate the versatility of our code, we apply it to six different surfaces with both known and unknown spectra:
\begin{itemize}
    \item \textit{Sphere}: the unit 2-sphere has eigenvalues $l(l+1)$ for $l$ a non-negative integer with multiplicity $2l+1$, with eigenfunctions the spherical harmonics \cite{berger_2000};
    \item \textit{Hemisphere}: the 2-hemisphere with Dirichlet boundary conditions has eigenvalues given by $l(l+1)$ for $l$ a positive integer, with multiplicity $l$ \cite{hemisphere};
    \item \textit{Blob}: a unit sphere deformed by a smooth, random function;
    \item \textit{Torus}: the surface of revolution of $(x-2)^2+z^2=1$ around the $y$-axis;
    \item \textit{M{\"o}bius strip}: the parametrization obtained by sweeping out a line segment; this surface is non-orientable;
    \item \textit{Triceratops}: the imported mesh as discussed in the previous subsection.
\end{itemize}

\begin{figure}[t]
\centering
\setlength{\tabcolsep}{6pt}
\begin{tabular}{>{\centering\arraybackslash}m{0.25\textwidth}
                >{\centering\arraybackslash}m{0.32\textwidth}
                >{\centering\arraybackslash}m{0.32\textwidth}}
\toprule
\textbf{1000th Efun} & \textbf{Count} & \textbf{Ratio} \\
\midrule

\begin{tabular}{@{}c@{}}
$\lambda = 992.0000000$\\
$\operatorname{res}(\lambda)=2.87\times 10^{-9}$\\[3pt]
\includegraphics[height=0.10\textwidth]{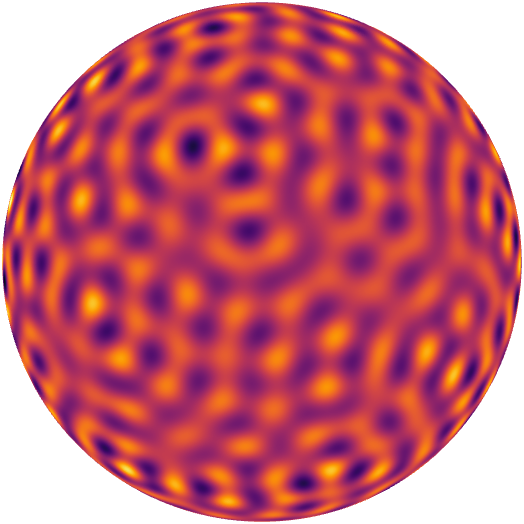}
\end{tabular}
&
\raisebox{-0.5\height}{\includegraphics[width=0.32\textwidth]{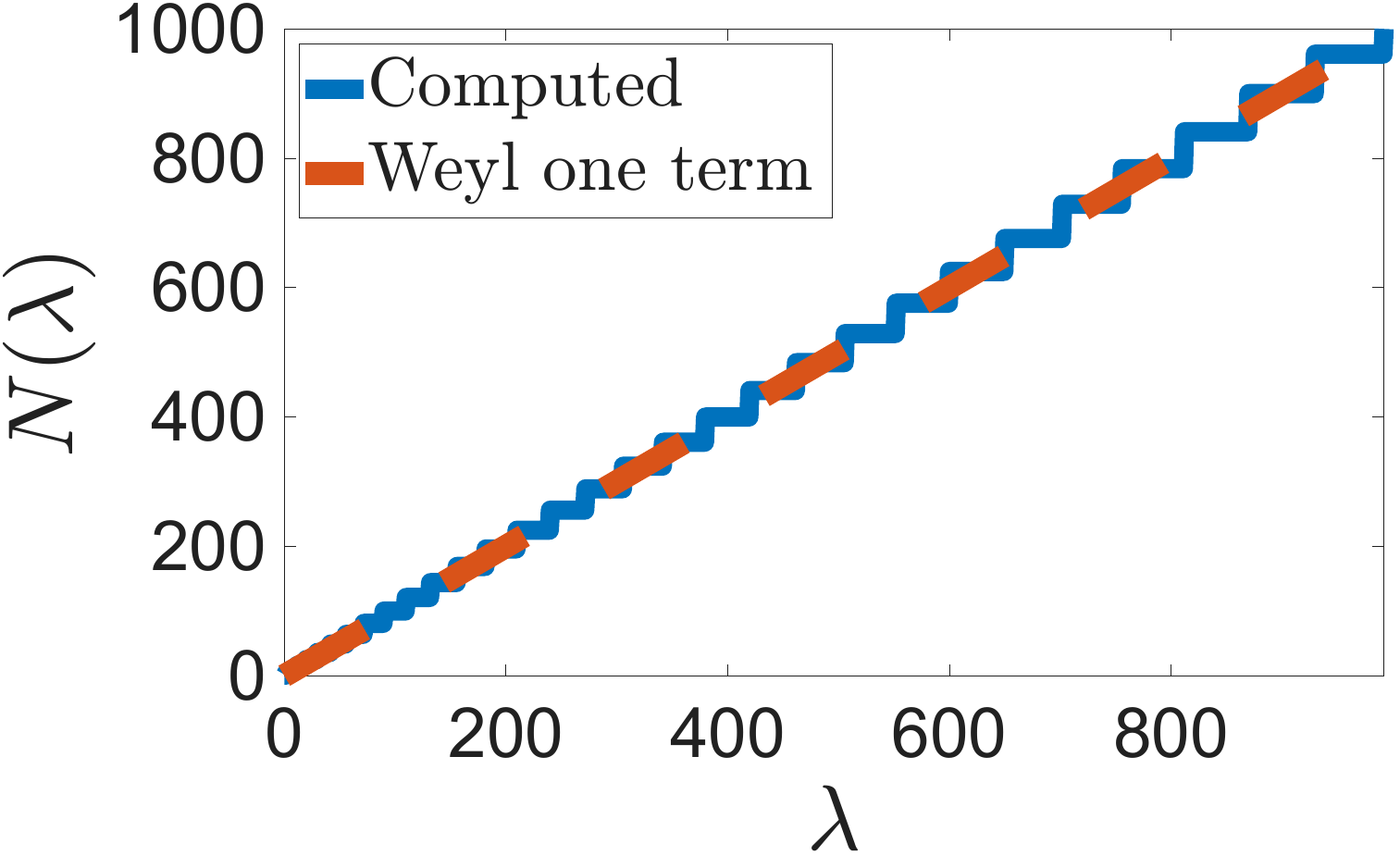}}
&
\raisebox{-0.5\height}{\includegraphics[width=0.32\textwidth]{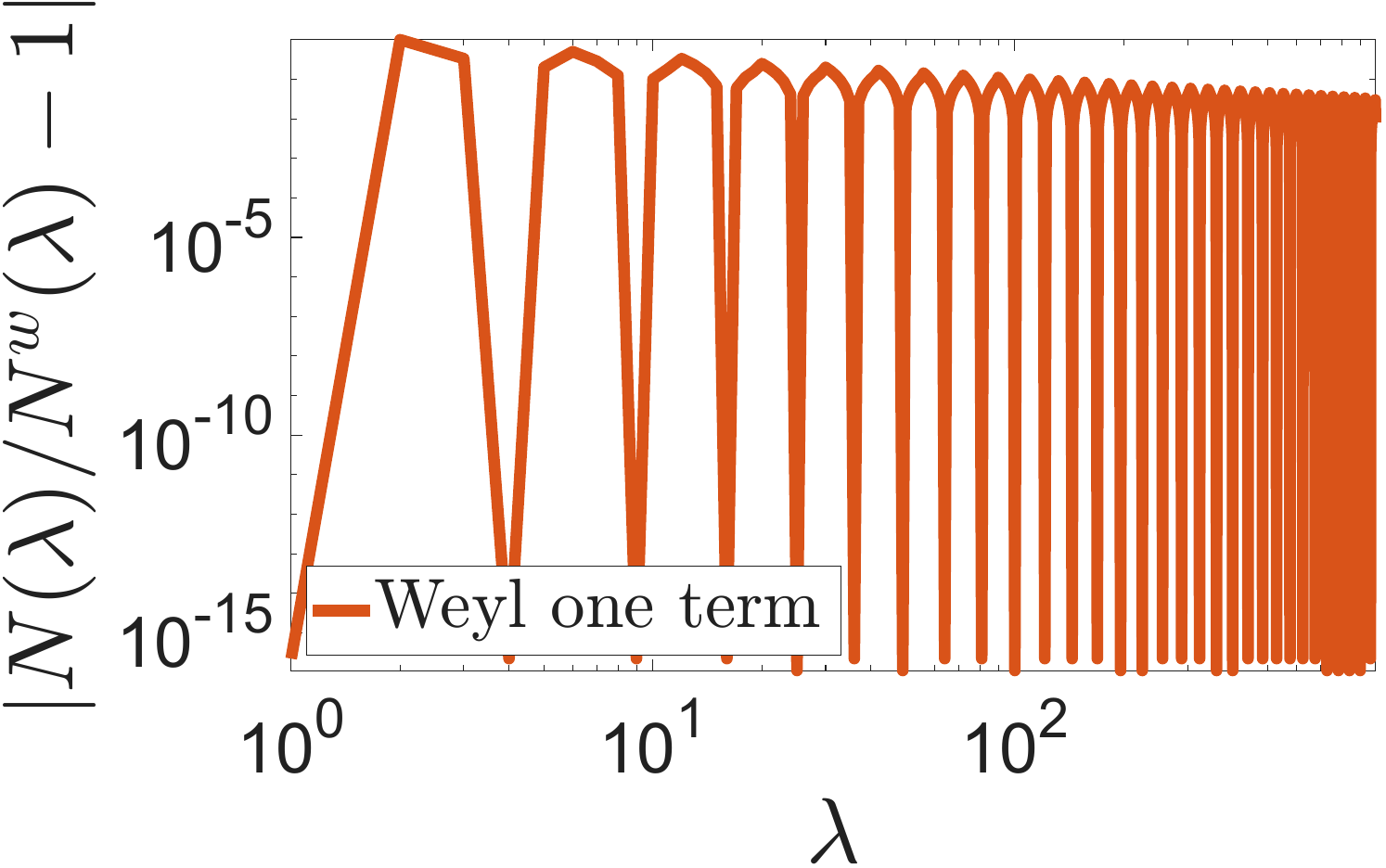}}
\\\addlinespace[3pt]

\begin{tabular}{@{}c@{}}
$\lambda = 2070.00000$\\
$\operatorname{res}(\lambda)=9.77\times 10^{-7}$\\[3pt]
\includegraphics[height=0.10\textwidth]{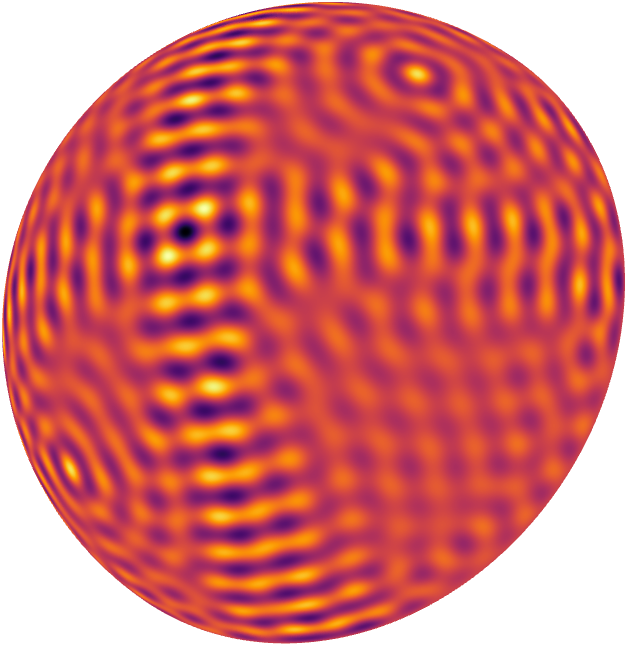}
\end{tabular}
&
\raisebox{-0.5\height}{\includegraphics[width=0.32\textwidth]{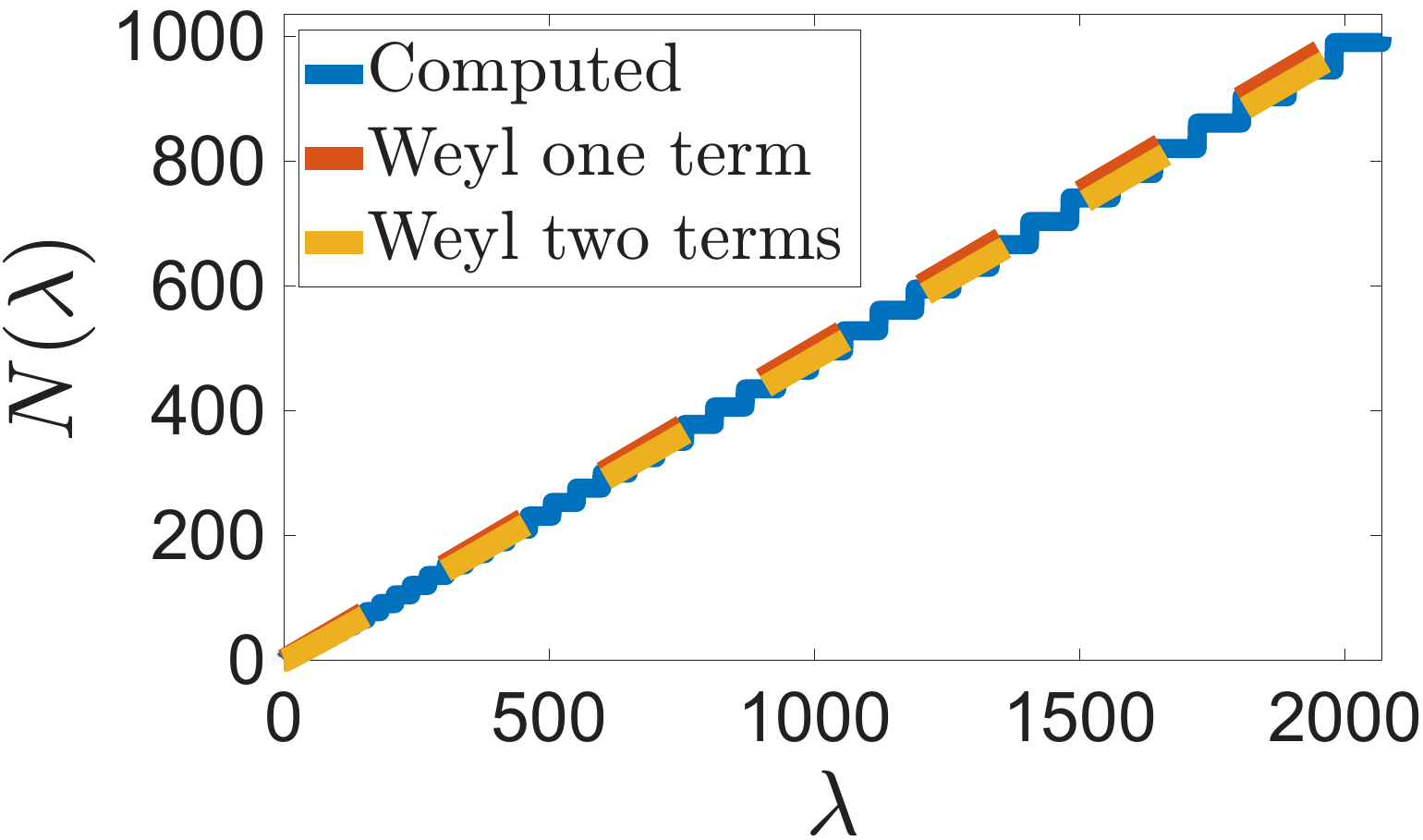}}
&
\raisebox{-0.5\height}{\includegraphics[width=0.32\textwidth]{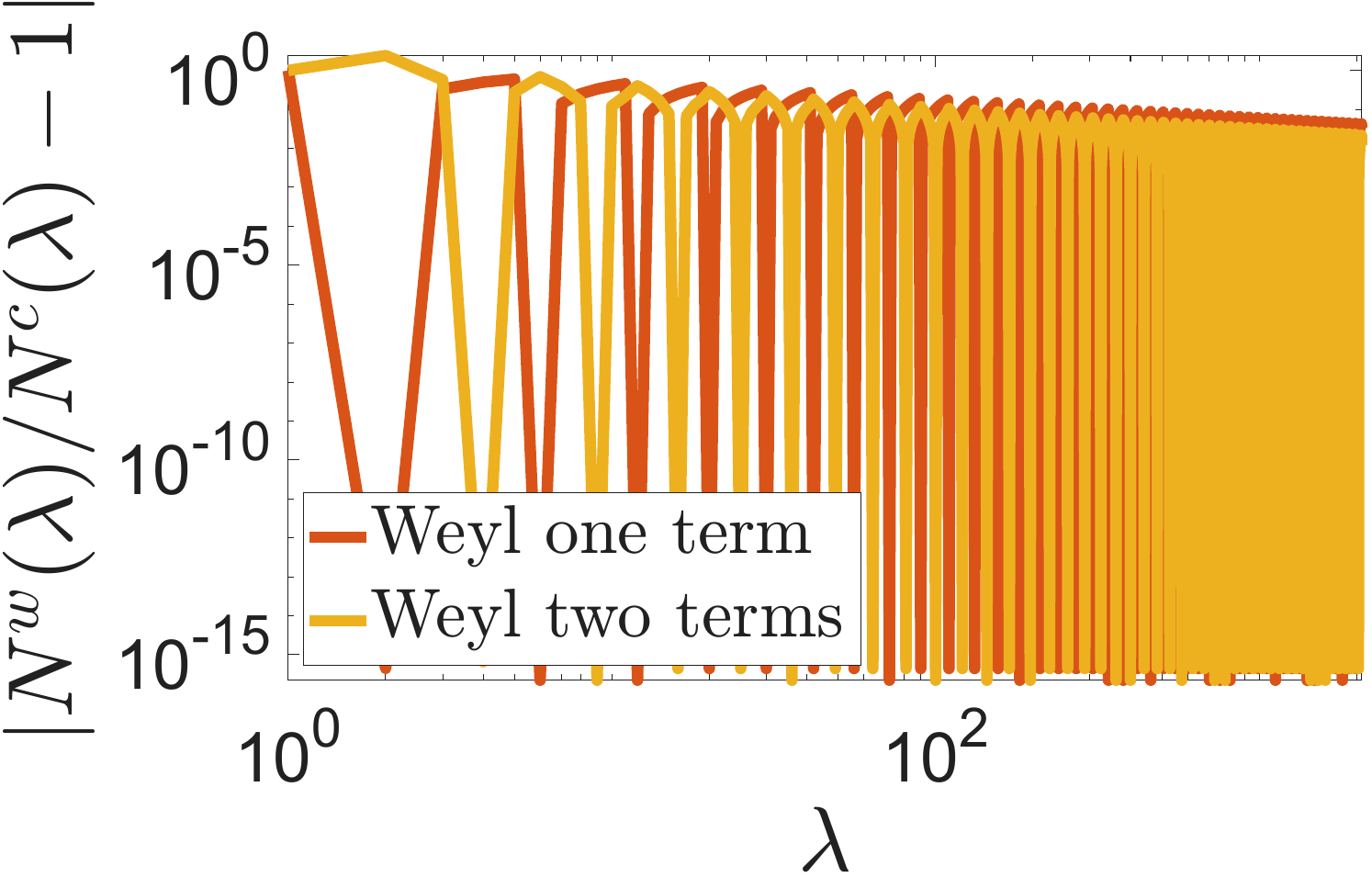}}
\\\addlinespace[3pt]

\begin{tabular}{@{}c@{}}
$\lambda = 682.963454$\\
$\operatorname{res}(\lambda)=1.19\times 10^{-7}$\\[3pt]
\includegraphics[height=0.10\textwidth]{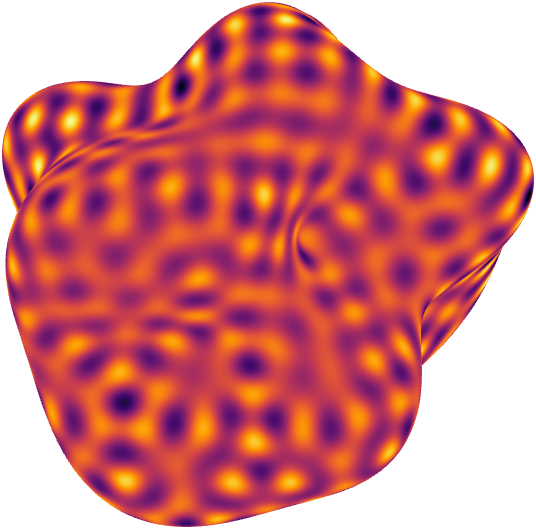}
\end{tabular}
&
\raisebox{-0.5\height}{\includegraphics[width=0.32\textwidth]{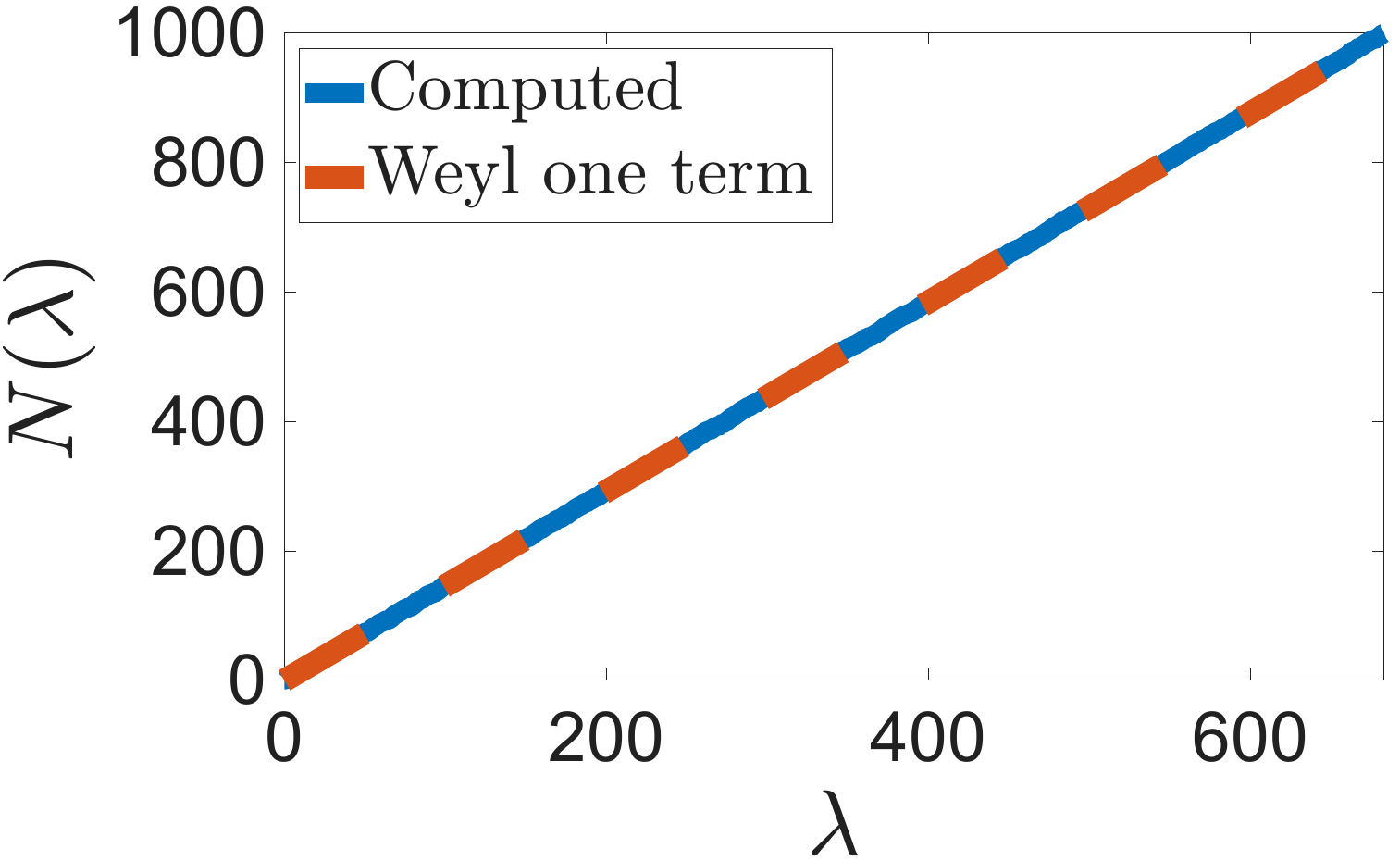}}
&
\raisebox{-0.5\height}{\includegraphics[width=0.32\textwidth]{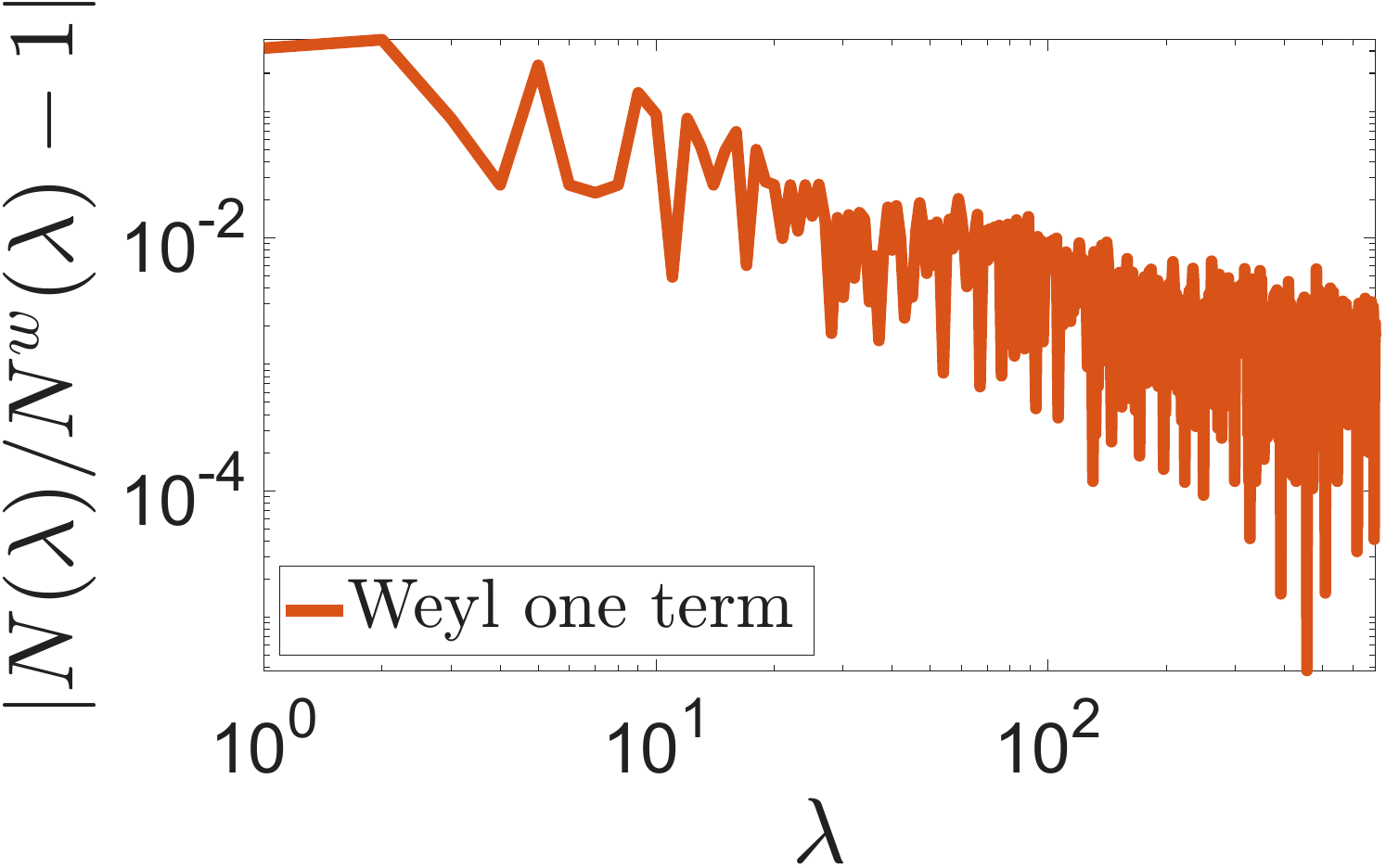}}
\\\addlinespace[3pt]

\begin{tabular}{@{}c@{}}
$\lambda = 158.52046$\\
$\operatorname{res}(\lambda)=2.43\times 10^{-7}$\\[3pt]
\includegraphics[height=0.10\textwidth]{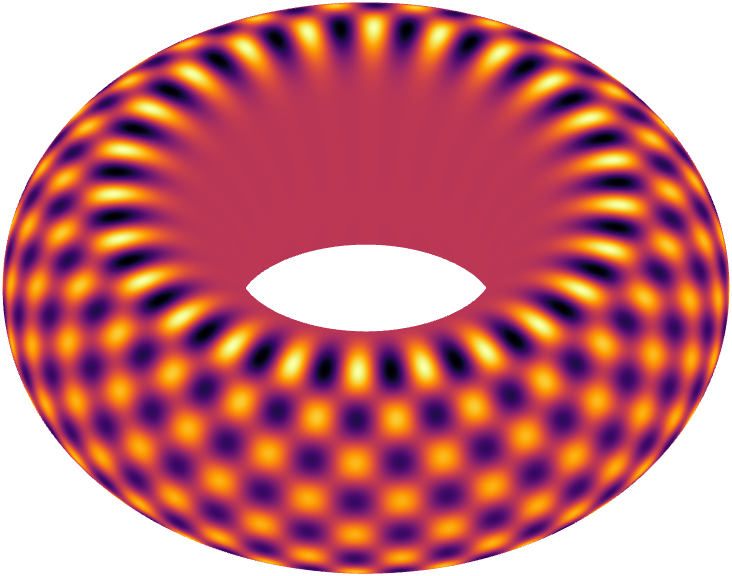}
\end{tabular}
&
\raisebox{-0.5\height}{\includegraphics[width=0.32\textwidth]{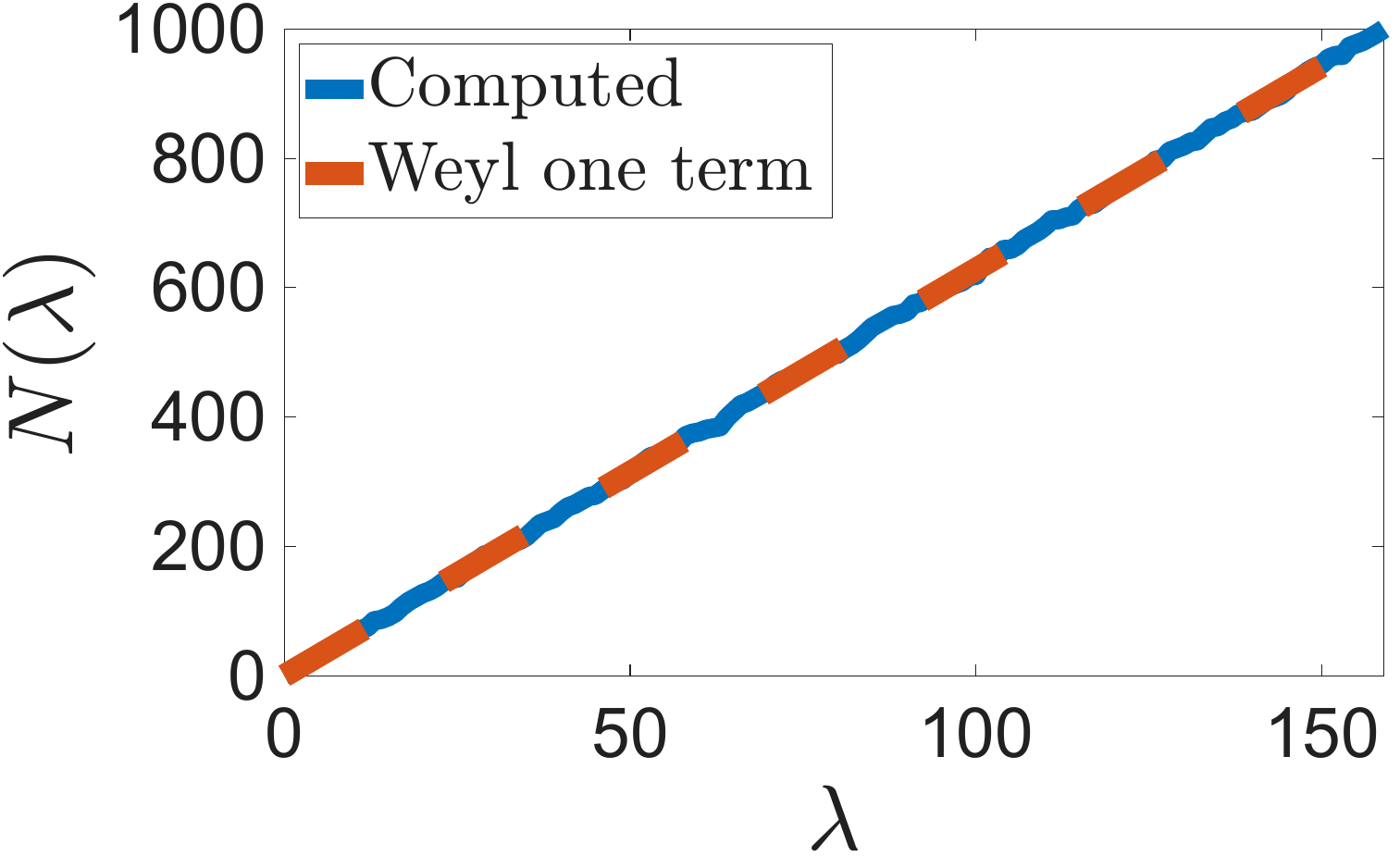}}
&
\raisebox{-0.5\height}{\includegraphics[width=0.32\textwidth]{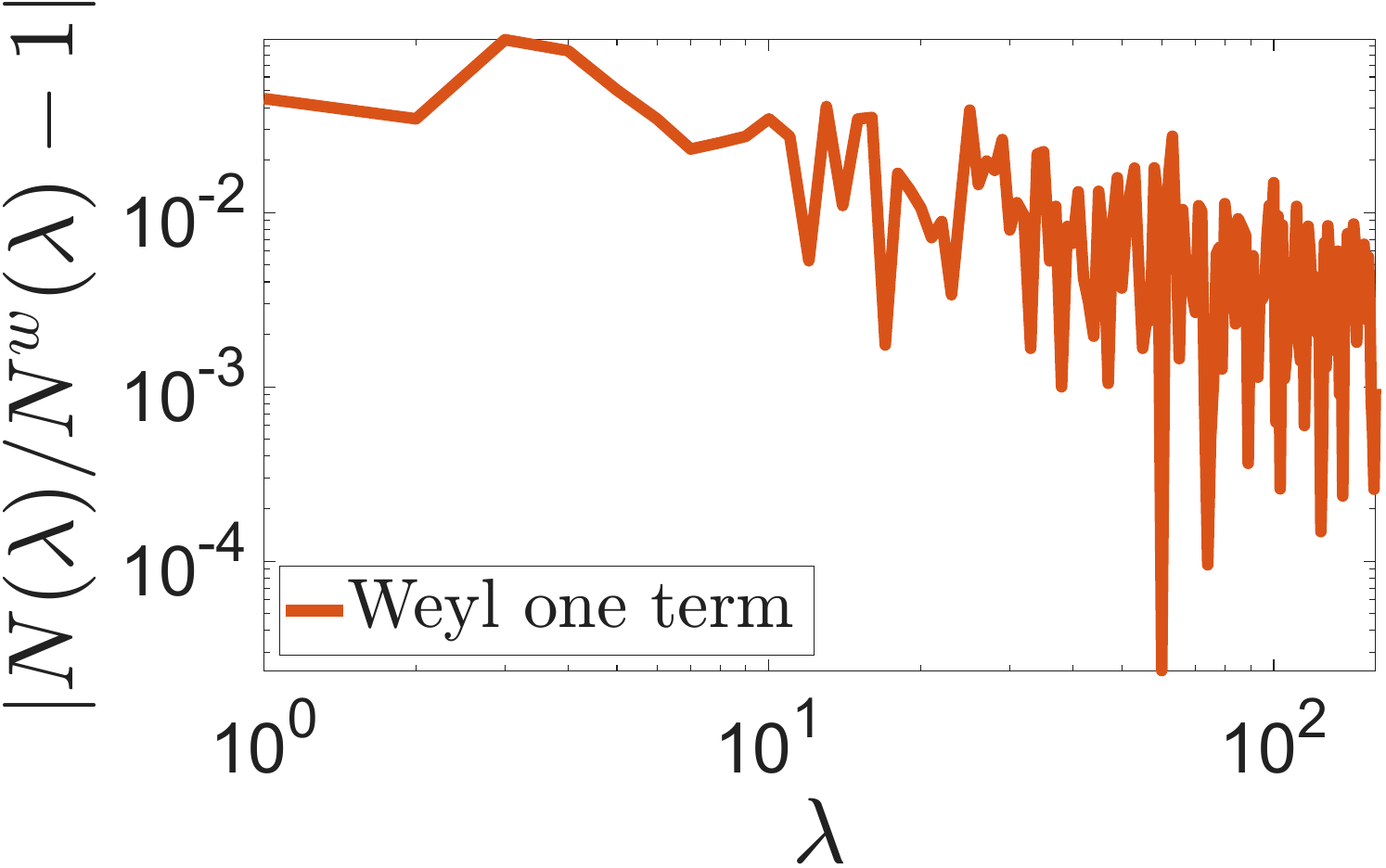}}
\\\addlinespace[3pt]

\begin{tabular}{@{}c@{}}
$\lambda = 2081.0685704$\\
$\operatorname{res}(\lambda)=3.28\times 10^{-9}$\\[3pt]
\includegraphics[height=0.10\textwidth]{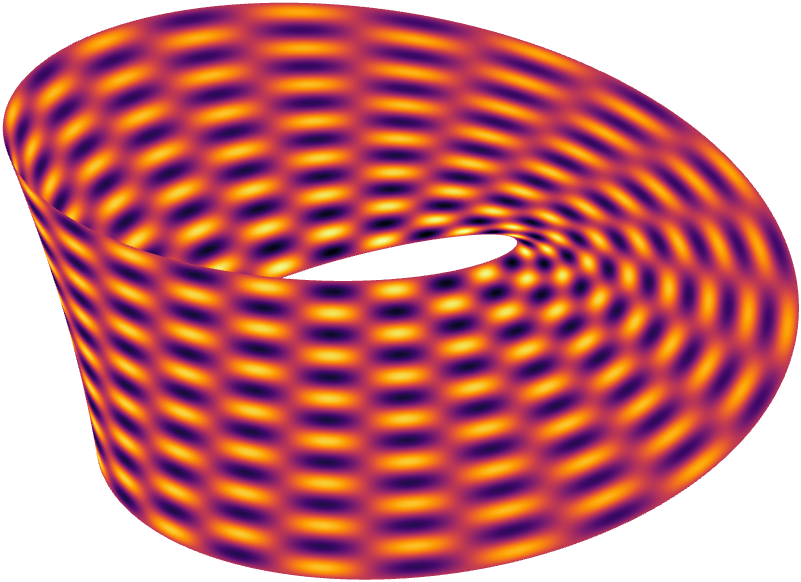}
\end{tabular}
&
\raisebox{-0.5\height}{\includegraphics[width=0.32\textwidth]{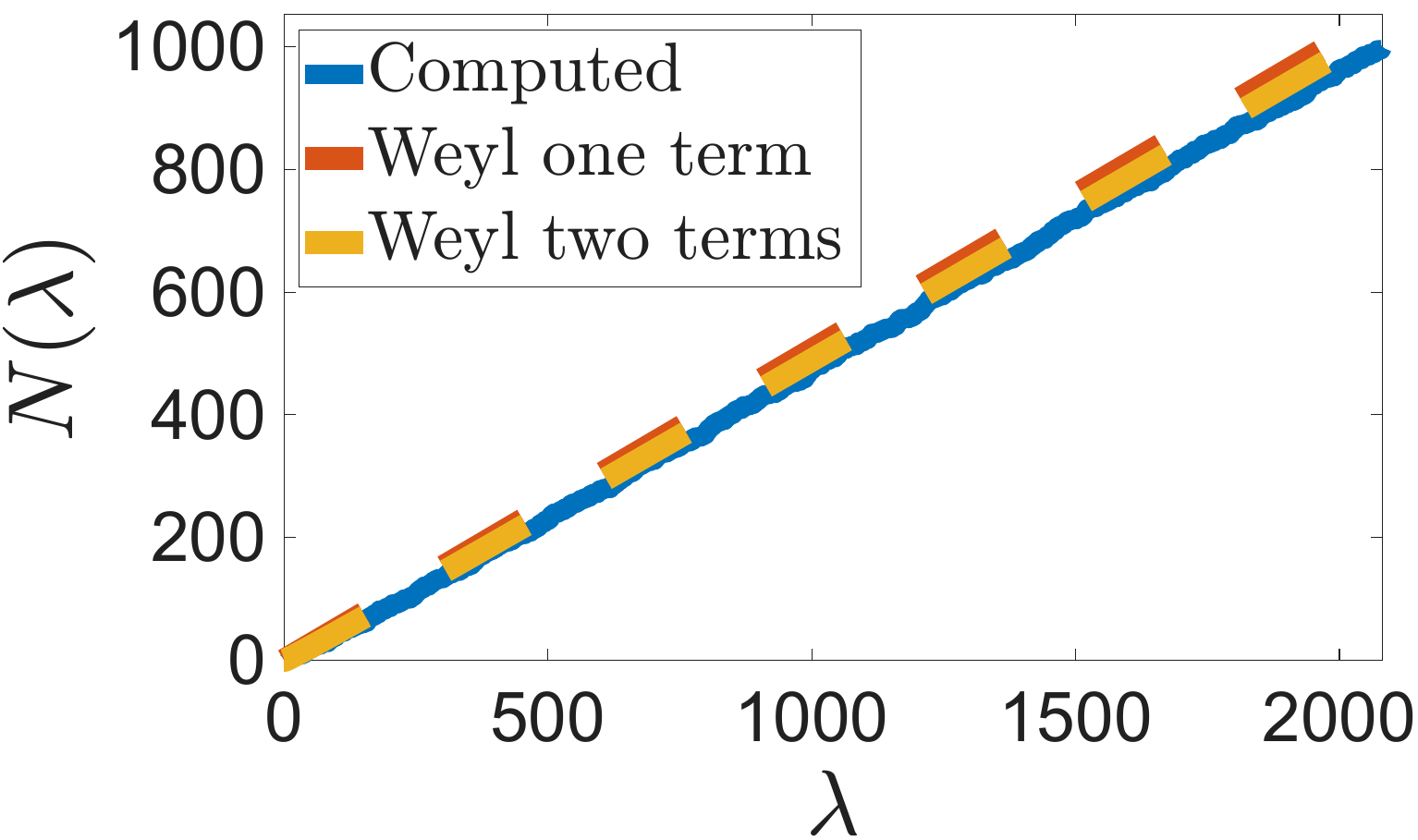}}
&
\raisebox{-0.5\height}{\includegraphics[width=0.32\textwidth]{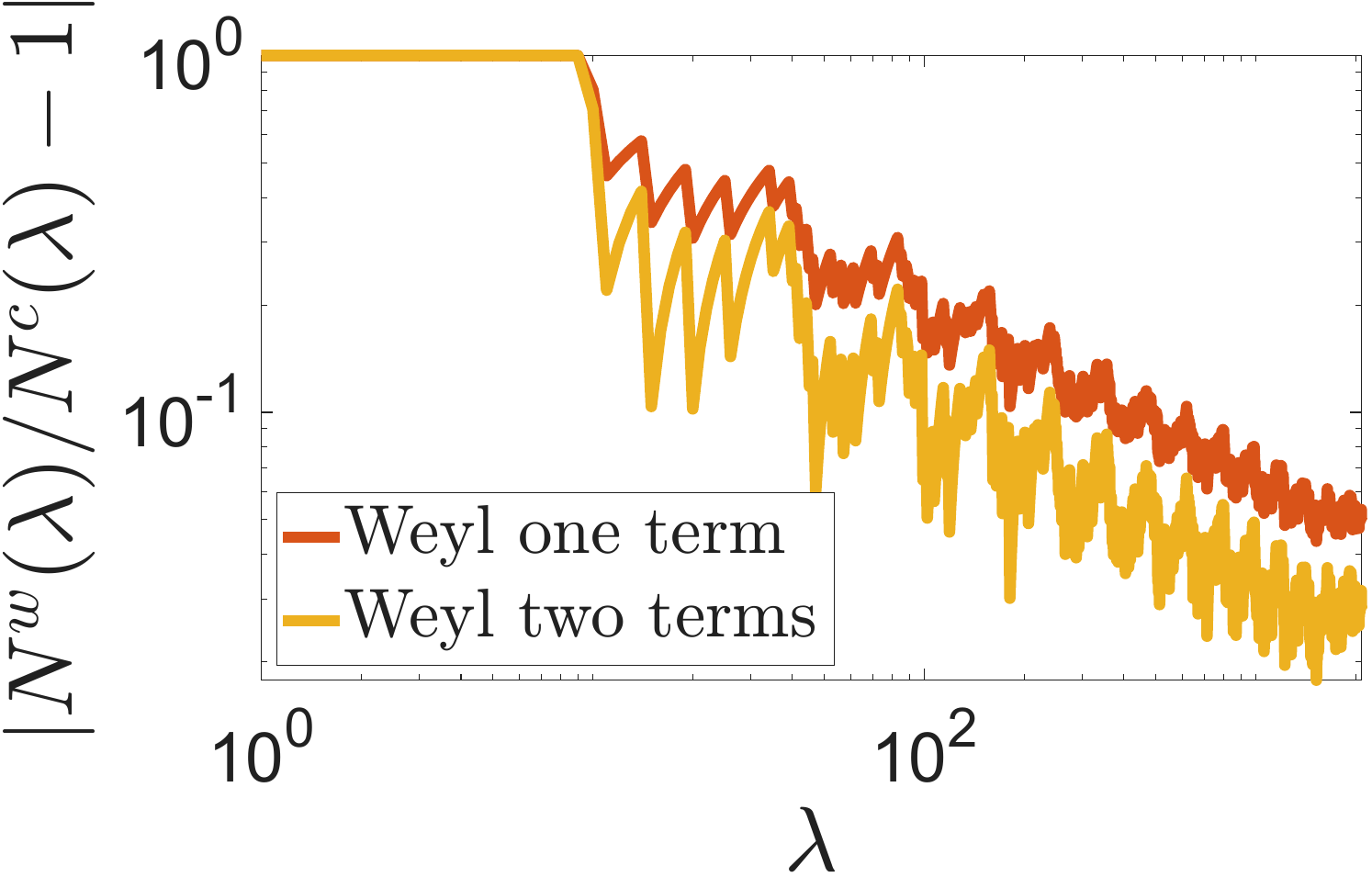}}
\\\addlinespace[3pt]

\begin{tabular}{@{}c@{}}
$\lambda = 887.819$\\
$\operatorname{res}(\lambda)=4.70\times 10^{-4}$\\[3pt]
\includegraphics[height=0.10\textwidth]{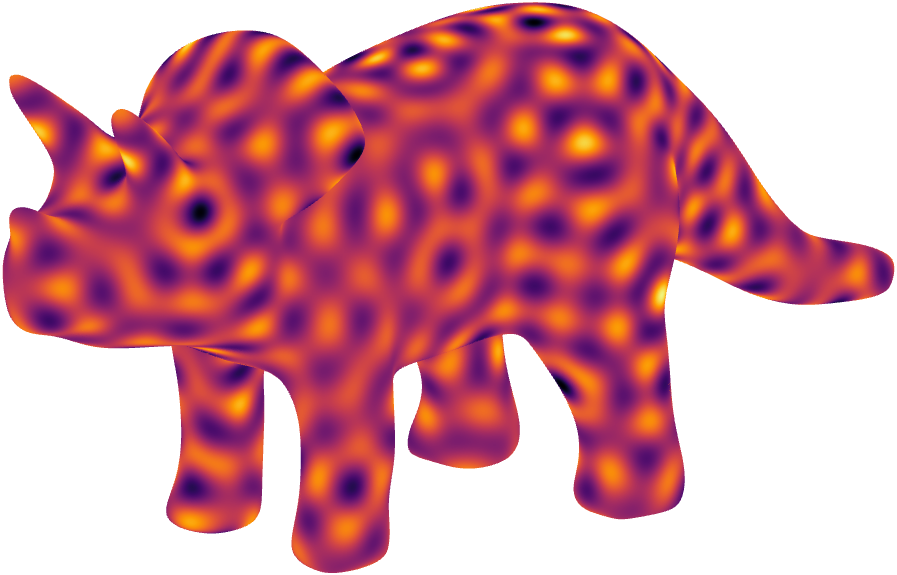}
\end{tabular}
&
\raisebox{-0.5\height}{\includegraphics[width=0.32\textwidth]{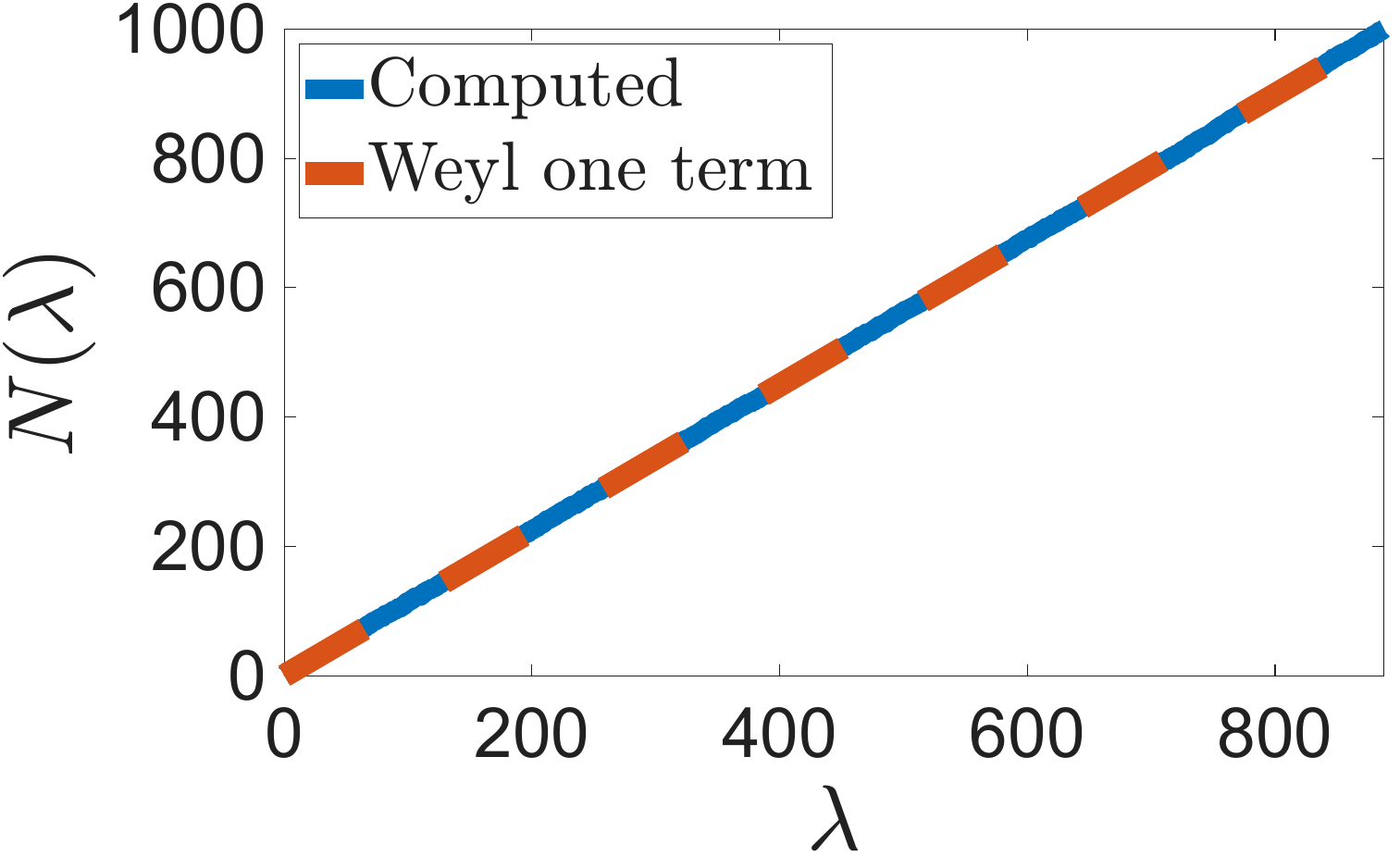}}
&
\raisebox{-0.5\height}{\includegraphics[width=0.32\textwidth]{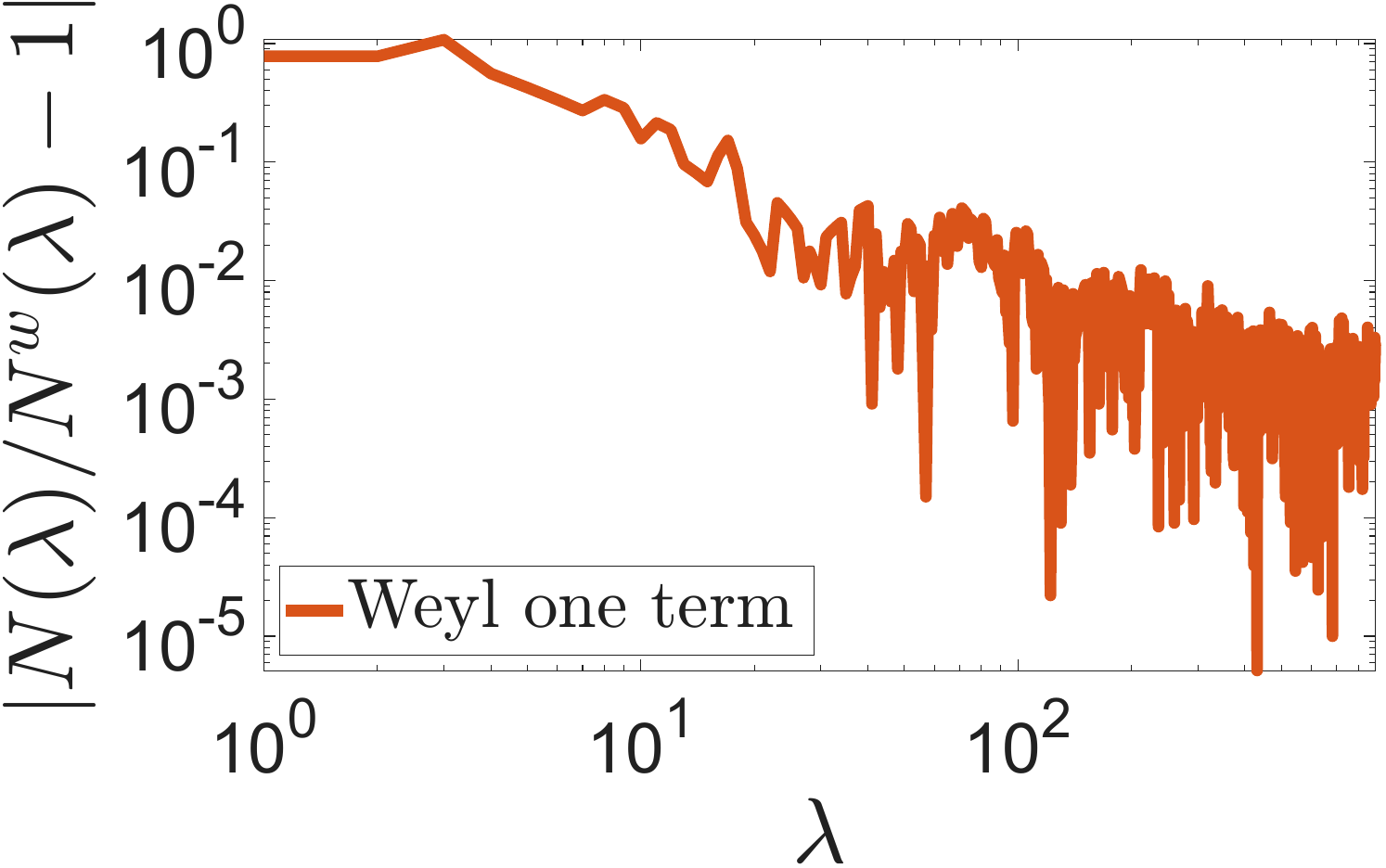}}
\\\addlinespace[3pt]

\bottomrule
\end{tabular}
\caption{Comparison of computed Laplace--Beltrami eigenvalues with Weyl asymptotics, together with the corresponding high-frequency eigenfunctions.}
\label{fig:laplacebeltrami_combined}
\end{figure}

\begin{figure}[t]
\centering
\setlength{\tabcolsep}{6pt}
\begin{tabular}{>{\centering\arraybackslash}m{0.25\textwidth}
                >{\centering\arraybackslash}m{0.32\textwidth}
                >{\centering\arraybackslash}m{0.32\textwidth}}
\toprule
\textbf{1000th Efun} & \textbf{Count} & \textbf{Ratio} \\
\midrule

\begin{tabular}{@{}c@{}}
$\lambda = 684.310343$\\
$\operatorname{res}(\lambda)=1.20\times 10^{-7}$\\[3pt]
\includegraphics[height=0.10\textwidth]{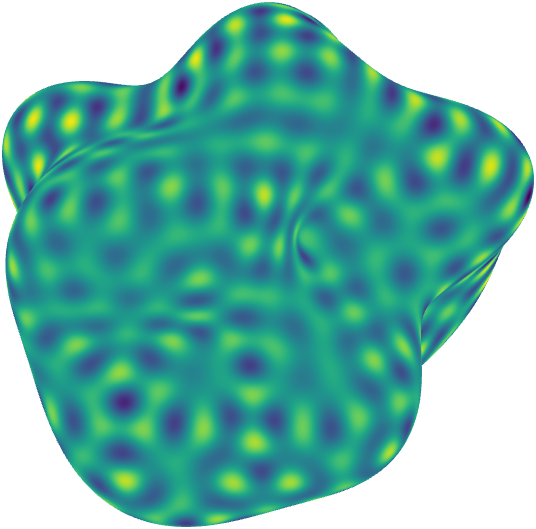}
\end{tabular}
&
\raisebox{-0.5\height}{\includegraphics[width=0.32\textwidth]{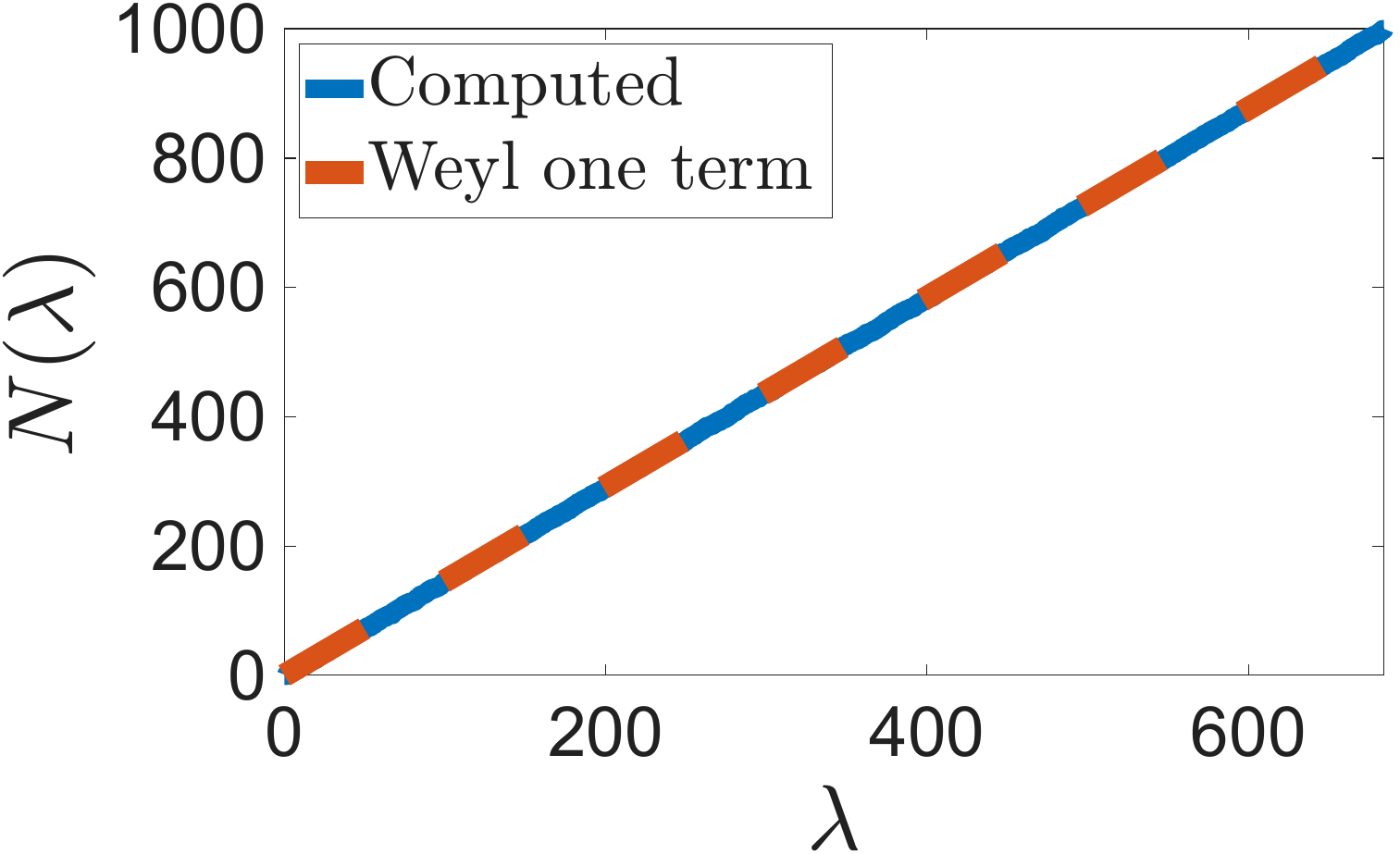}}
&
\raisebox{-0.5\height}{\includegraphics[width=0.32\textwidth]{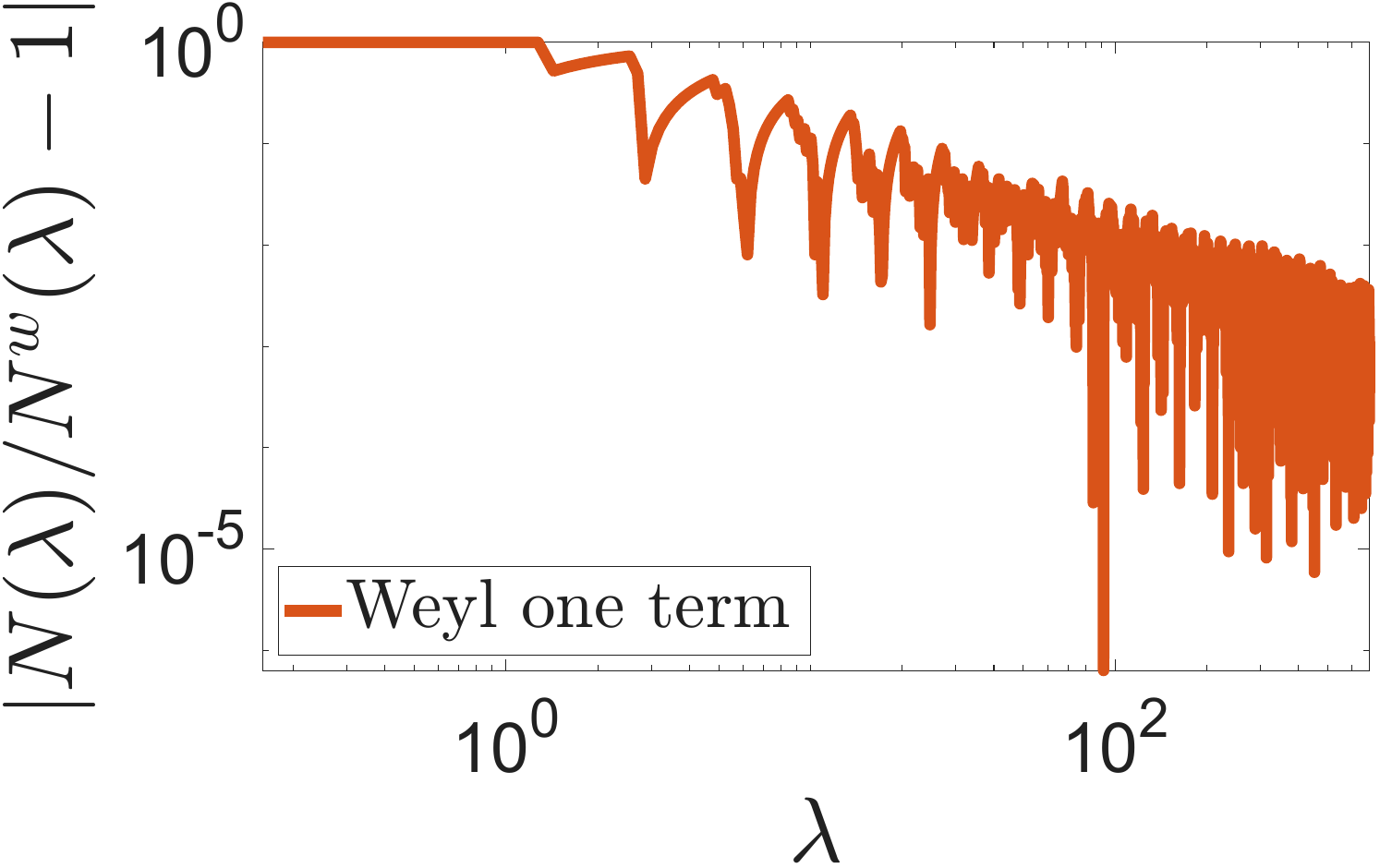}}
\\\addlinespace[3pt]

\begin{tabular}{@{}c@{}}
$\lambda = 165.3611933$\\
$\operatorname{res}(\lambda)=2.65\times 10^{-8}$\\[3pt]
\includegraphics[height=0.10\textwidth]{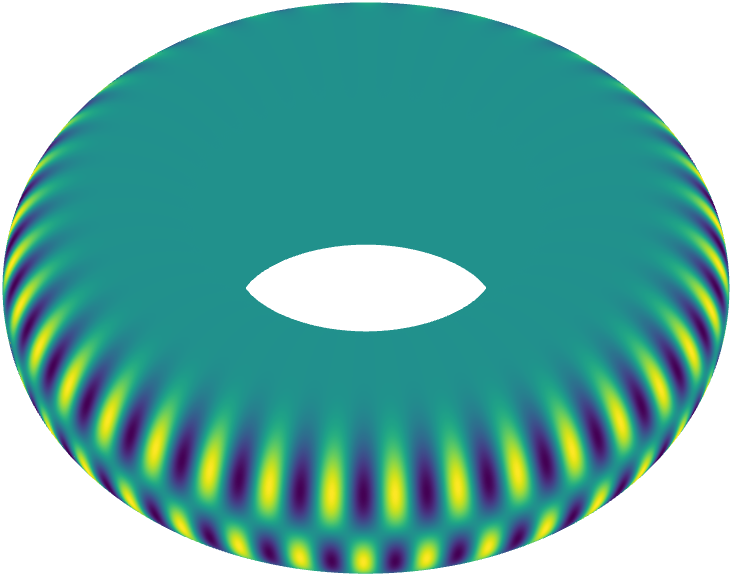}
\end{tabular}
&
\raisebox{-0.5\height}{\includegraphics[width=0.32\textwidth]{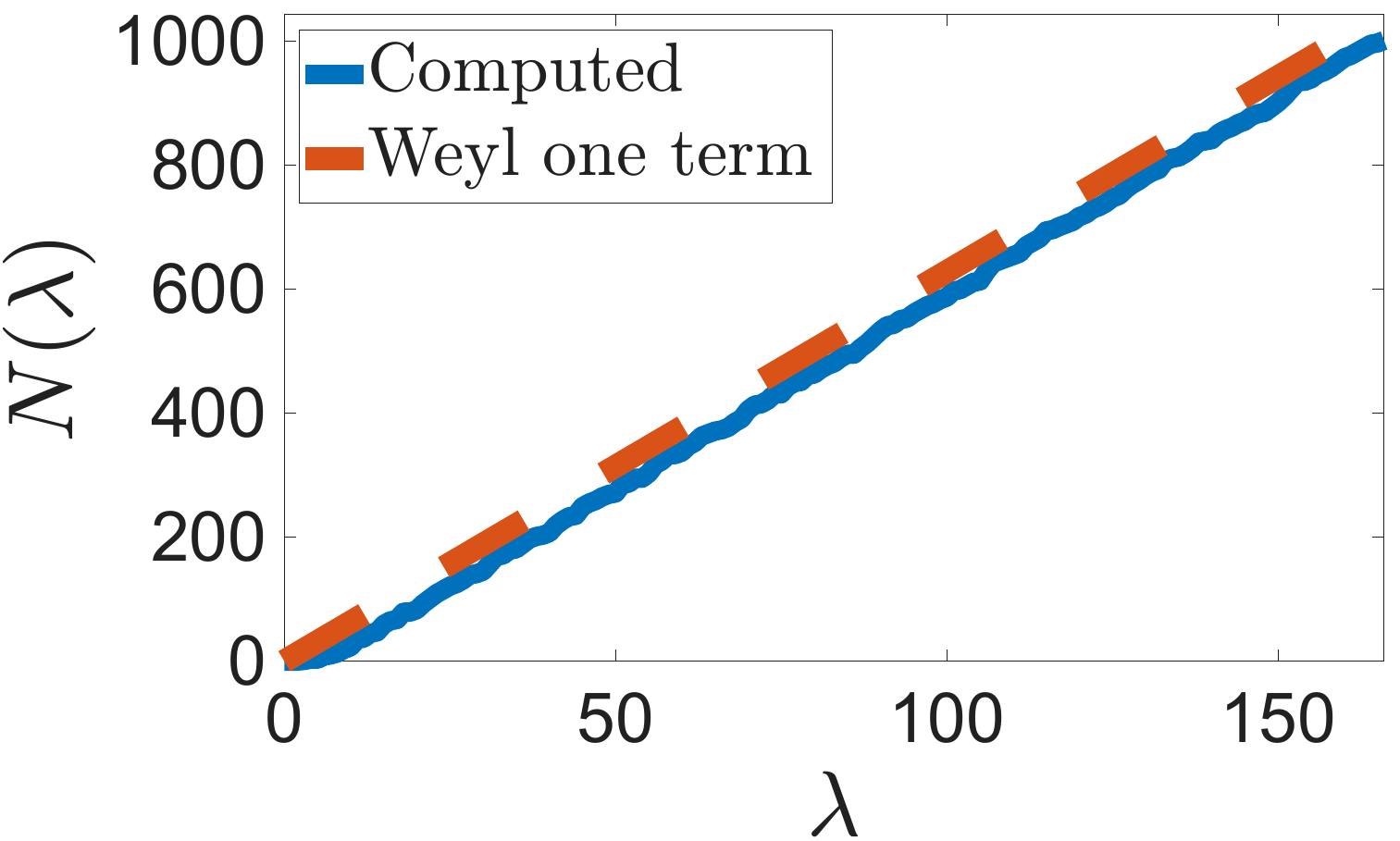}}
&
\raisebox{-0.5\height}{\includegraphics[width=0.32\textwidth]{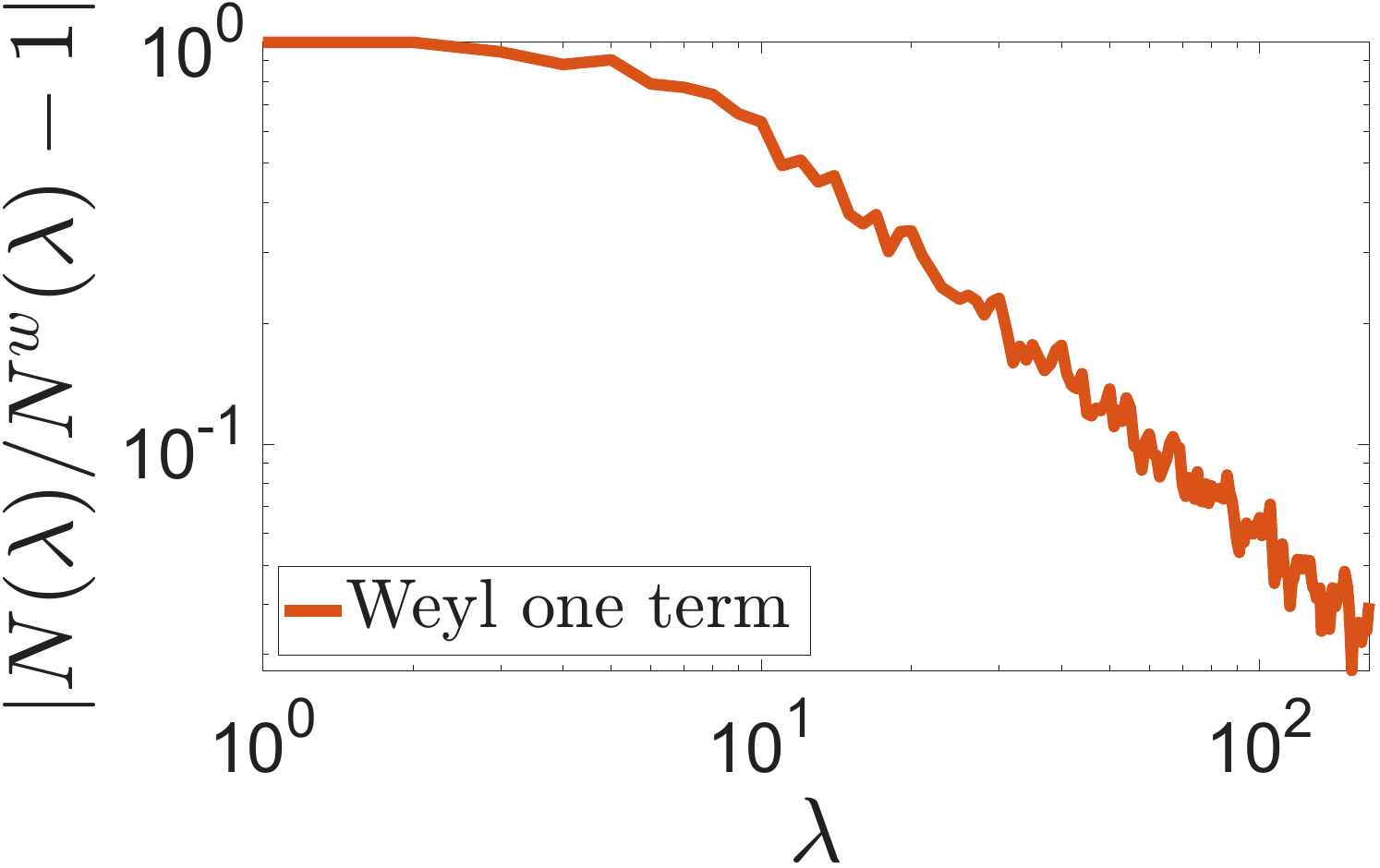}}
\\\addlinespace[3pt]

\begin{tabular}{@{}c@{}}
$\lambda = 2068.0469$\\
$\operatorname{res}(\lambda)=1.84\times 10^{-5}$\\[3pt]
\includegraphics[height=0.10\textwidth]{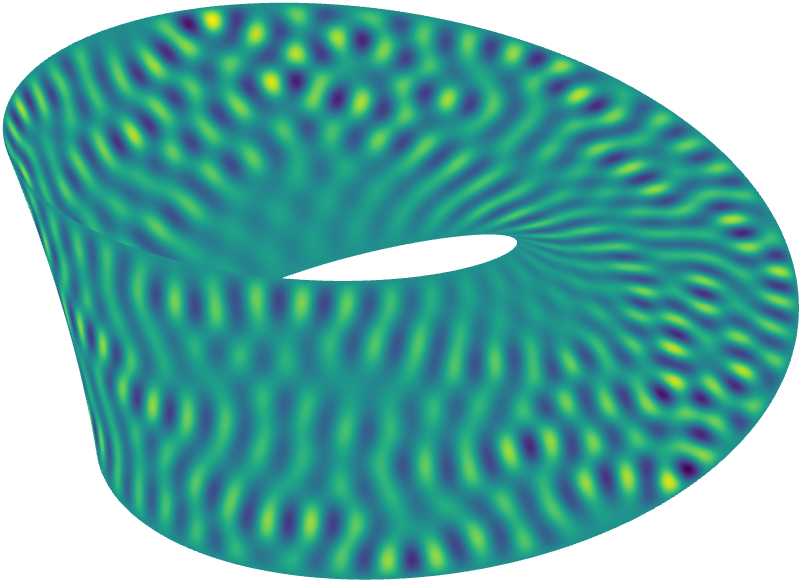}
\end{tabular}
&
\raisebox{-0.5\height}{\includegraphics[width=0.32\textwidth]{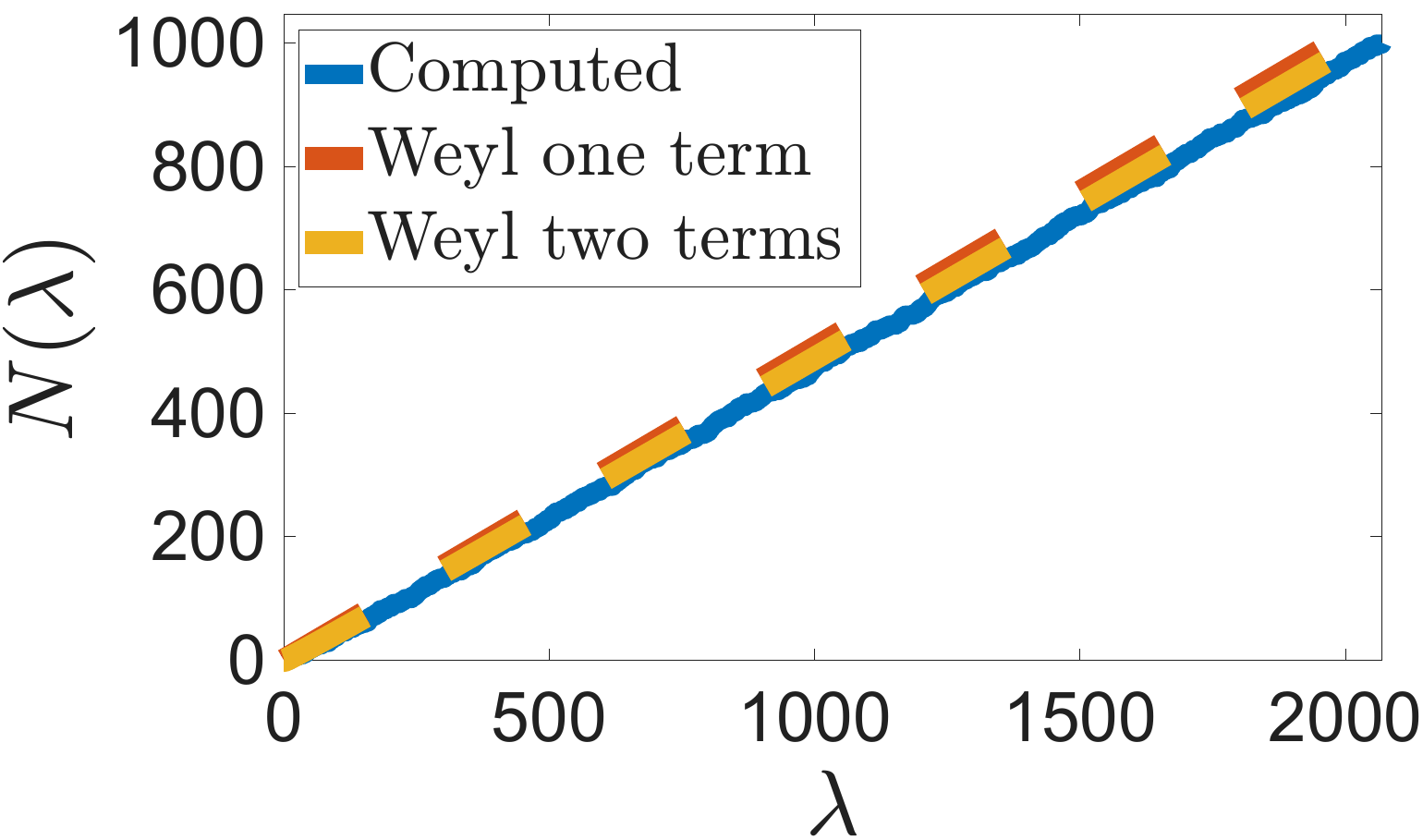}}
&
\raisebox{-0.5\height}{\includegraphics[width=0.32\textwidth]{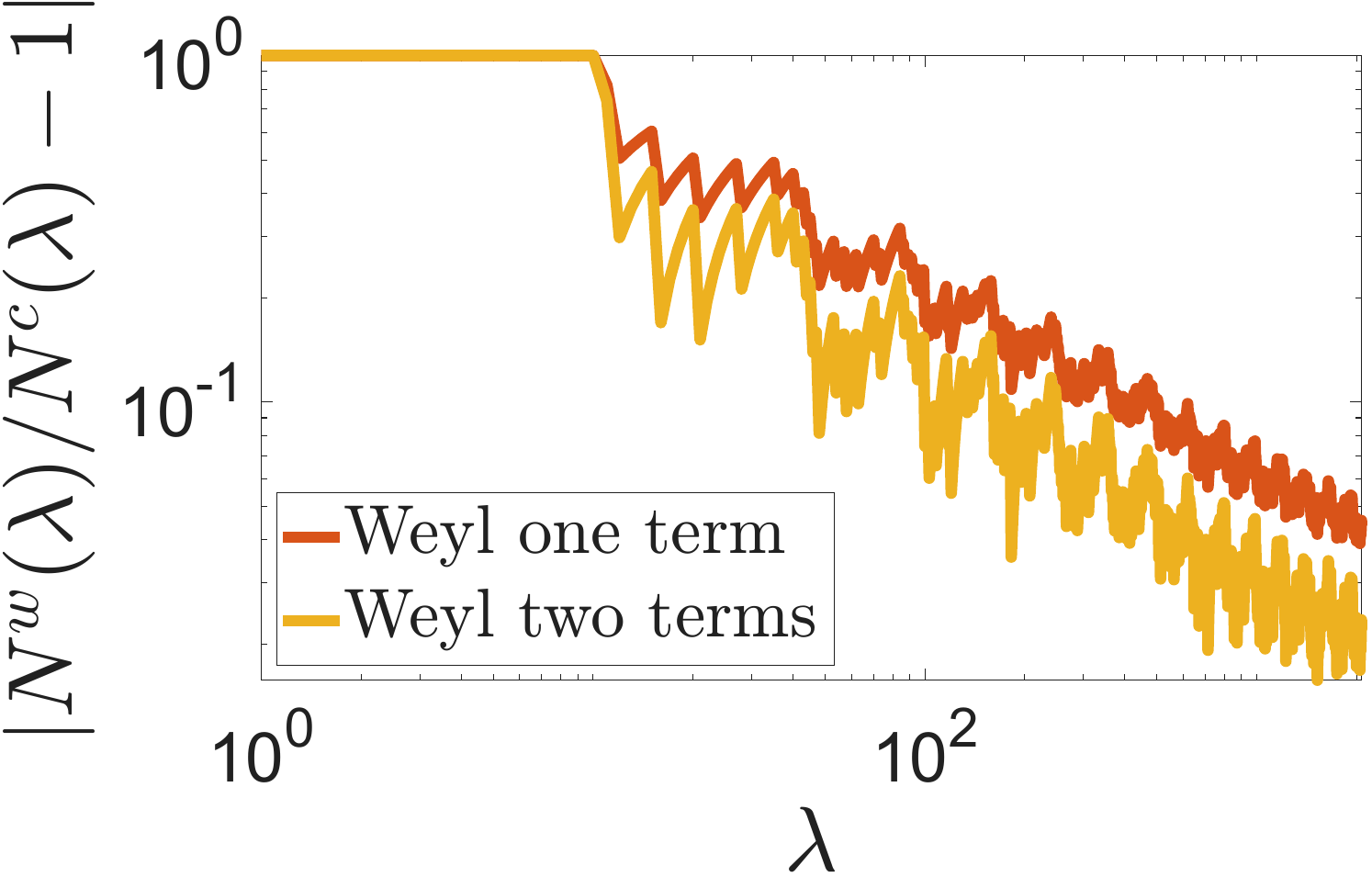}}
\\\addlinespace[3pt]

\begin{tabular}{@{}c@{}}
$\lambda = 889.02$\\
$\operatorname{res}(\lambda)=6.18\times 10^{-4}$\\[3pt]
\includegraphics[height=0.10\textwidth]{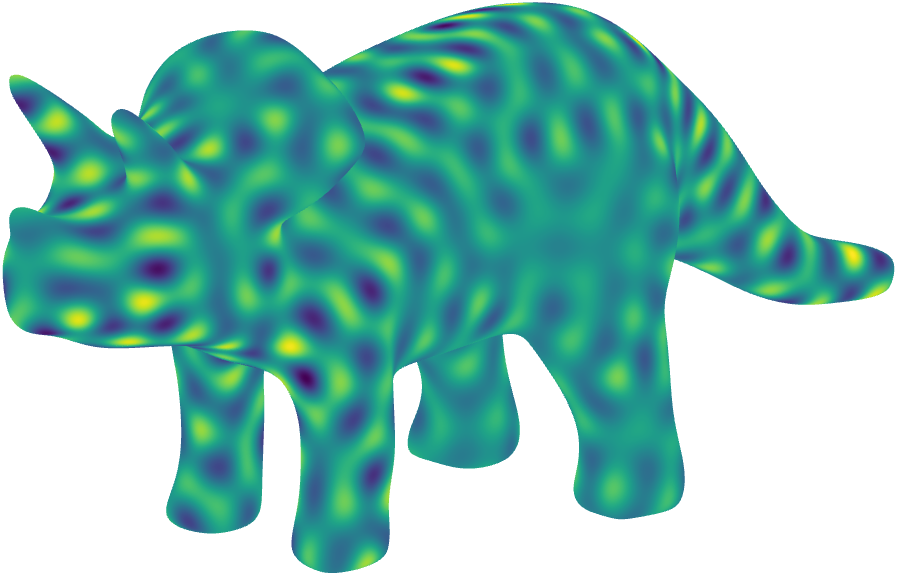}
\end{tabular}
&
\raisebox{-0.5\height}{\includegraphics[width=0.32\textwidth]{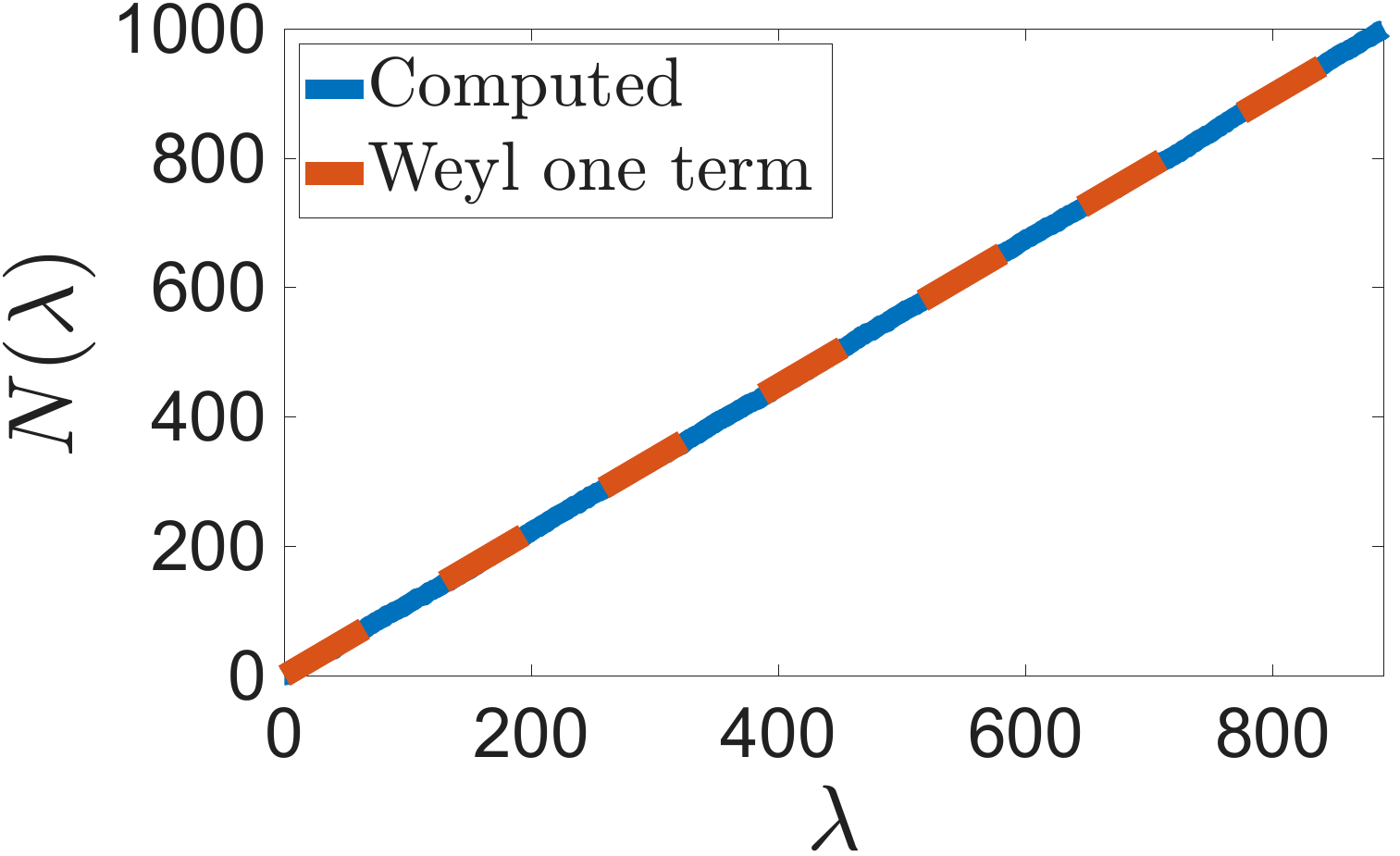}}
&
\raisebox{-0.5\height}{\includegraphics[width=0.32\textwidth]{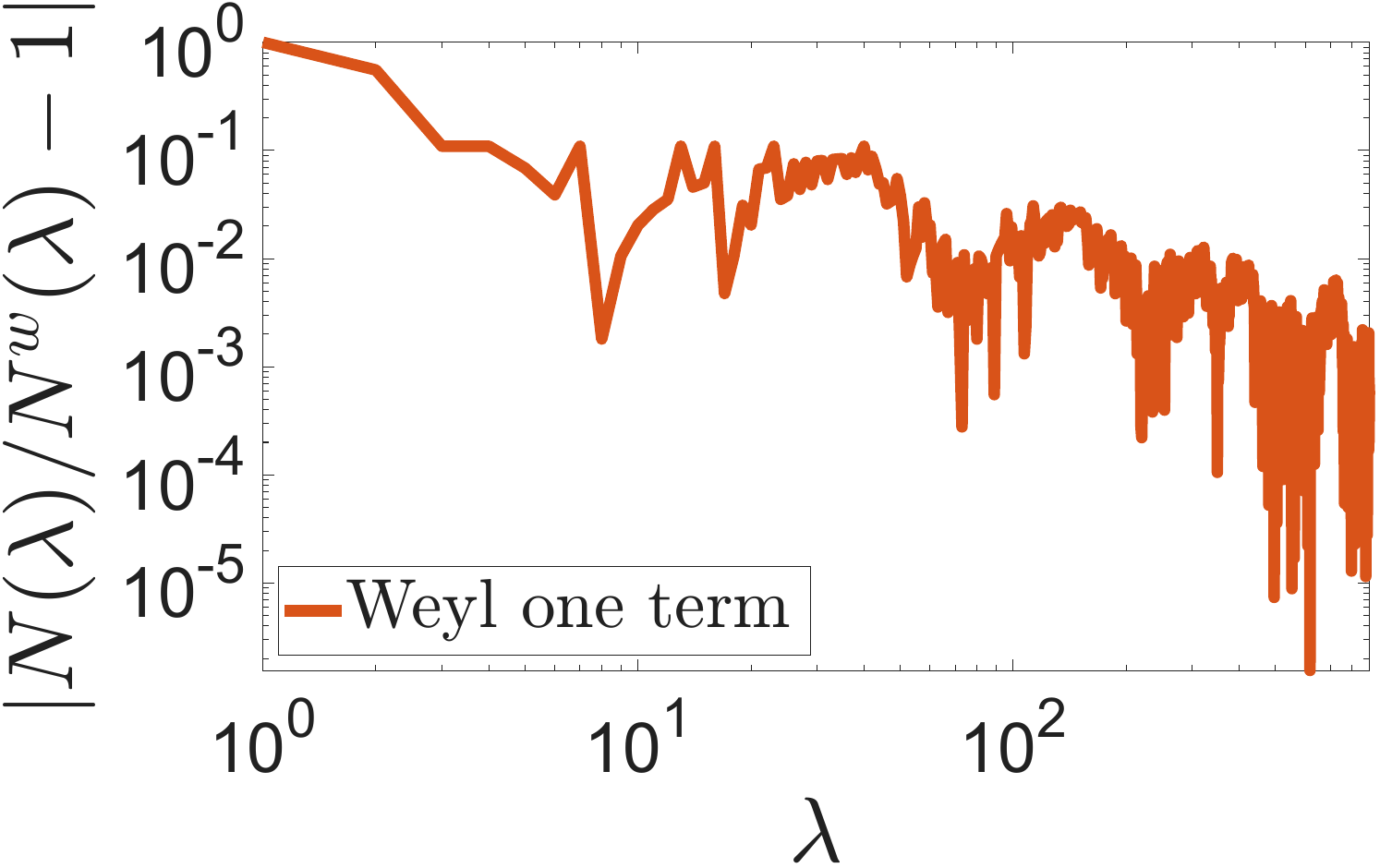}}
\\\addlinespace[3pt]

\bottomrule
\end{tabular}
\caption{Comparison of computed Schr\"odinger eigenvalues with Weyl asymptotics, together with the corresponding high-frequency eigenfunctions.}
\label{fig:schrodinger123}
\end{figure}

{For surfaces with boundary, in each case we use Dirichlet boundary conditions.} As before, the Laplace--Beltrami operator and Schr\"odinger operator with harmonic oscillator potential on each of these surfaces has nonnegative, real eigenvalues, so we can systematically compute the eigenvalues by taking elliptical contours centered on the real axis; by doing so, we compute the first $1000$ eigenvalues, eigenfunctions and residuals of the Laplace--Beltrami and Schr\"odinger operators. Since both operators are self-adjoint, the residuals bound the distance to the true eigenvalue of the infinite-dimensional operator (see \cref{errctrl}). {The parameters are chosen adaptively to ensure small residuals and $m$ is (over)estimated based on Weyl's law; for example, for the triceratops we used the parameters described in \cref{sec:triceratops}.}
We then ``zoom in'' on the $1000\text{th}$ such eigenvalues for each surface/operator pair, taking $m_0=m=1$ and refining the surface to produce small residuals. The eigenvalues, eigenfunctions and residuals are shown in \cref{fig:laplacebeltrami_combined,fig:schrodinger123}; for eigenfunctions with multiplicity, we display the eigenfunction with the largest residual to show the worst-case {problem}. Note that for the sphere and hemisphere we require a larger number of quadrature nodes than usual due to the high multiplicity. We do not show the results for the Schr\"odinger operator on the sphere and hemisphere, as harmonic oscillator potential is constant on these surfaces. 

We also compare the eigenvalue counting function $N(\lambda)$ as computed by our algorithm to the asymptotic expansion $N^w(\lambda)$ given by Weyl's law \cref{weylone}; for open surfaces, we compare to both \cref{weylone} and \cref{weyltwo}. We see that there is strong agreement and as $\lambda$ increases, the relative error between the two quantities decreases, as the asymptotic expansion becomes more accurate for large $\lambda$. The error decreases faster for open surfaces when we include the second term in Weyl's law for the Laplace--Beltrami operator, as expected, and also for the Schr\"odinger operator on the M\"obius strip. For the sphere and hemisphere, the plotted Weyl approximation and the discrete count coincide, leading to the observed machine-precision dips.

\subsection{Nonlinear differential eigenvalue problems}

Contour-integral-based methods are also effective for nonlinear eigenvalue problems (NLEPs). For example, recent work \cite{colbrook2025avoiding} has generalized Beyn's method \cite{beyn}. We consider the eigenvalue problem 
\begin{equation}\setlength\abovedisplayskip{6pt}\setlength\belowdisplayskip{6pt}\label{nlep}
T(\lambda)f=0,\quad T(\lambda):\mathcal{D}(T)\rightarrow\mathcal{H},\quad f\in\mathcal{D}(T),
\end{equation}
where $T$ may depend nonlinearly on $\lambda$ and $\mathcal{D}(T)\subset \mathcal{H}$ is dense.
The spectrum of $T$, $\spec(T)$, is the set of points $\lambda\in\mathbb{C}$ such that $T(\lambda):\mathcal{D}(T)\rightarrow\mathcal{H}$ is not boundedly invertible. Fix an open subset $\Omega\subset\mathbb{C}$. We make the following technical assumptions:
\begin{itemize}
	\item for each $\lambda\in\Omega$, $T(\lambda)$ is a Fredholm operator;
	\item for each $u\in\mathcal{D}(T)$, $\lambda\mapsto T(\lambda)u$ is holomorphic;
	\item there exists at least one point in $\Omega$ not in the spectrum $\spec(T)$.
\end{itemize}
Under these conditions the spectrum is discrete in $\Omega$~\cite[Thm.~1.3.1]{mennicken2003non}. {For example, the nonlinear eigenvalue problem \cref{eqn:nlep_example} with Dirichlet boundary conditions on the M\"obius strip satisfies these conditions: $T(\lambda)$ is Fredholm for every $\lambda$ by elliptic PDE theory \cite[Sec.~5.1]{taylor1996partial}, $\lambda\mapsto T(\lambda)u$ is holomorphic in $\lambda$ for any $u\in\mathcal{D}(\mathcal{T})$, and by considering a relevant variational principle it can be shown that the NLEP only has real eigenvalues \cite[Sec.~3]{keldysh} and so any elliptical or circular domain $\Omega$ satisfies the third condition.}

The key theoretical tool in applying contour integration methods is Keldysh's theorem \cite{keldysh}. If $m$ is the sum of the algebraic multiplicities of the eigenvalues inside $\Omega$, then
\begin{equation*}
\setlength\abovedisplayskip{0pt}\setlength\belowdisplayskip{6pt}
T(z)^{-1} = A(zI-J)^{-1}B^*+R(z)\quad \forall z\in\Omega\backslash \spec(T),
\end{equation*}
where $J\in\mathbb{C}^{m\times m}$ is the Jordan canonical form for the eigenvalues inside $\Omega$, $A$ and $B$ are the normalized canonical Jordan chains and $R(z)$ is a holomorphic remainder.

By Cauchy's theorem, for any holomorphic function $f$ we have that
\begin{equation}\setlength\abovedisplayskip{6pt}\setlength\belowdisplayskip{6pt}
\label{eqn:cauchykeldysh}
\cfrac{1}{2\pi i}\int_{\partial\Omega}f(z)T(z)^{-1}\dd z=Af(J)B^*.
\end{equation}
Using these theoretical results, we develop a solver for nonlinear surface eigenvalue problems, \texttt{SurfBeyn}, based on \cite{colbrook2025avoiding}. As for \texttt{SurfFEAST}, the key idea is to relate the eigenvalues of \eqref{nlep} inside $\Omega$ to a finite-dimensional eigenvalue problem using contour integrals. We define 
\begin{equation*}\setlength\abovedisplayskip{6pt}\setlength\belowdisplayskip{6pt}
    A_0=\cfrac{1}{2\pi i}\int_{\partial\Omega}T(z)^{-1}\dd z,\quad\text{and}\quad A_1=\cfrac{1}{2\pi i}\int_{\partial\Omega}zT(z)^{-1}\dd z.
\end{equation*}
{As in \cref{sec:surfFEAST_linear}, we approximate $A_0$ and $A_1$ by applying them to test functions and using a quadrature rule. If we then take a (truncated, if we do not know $m$ exactly) SVD of our approximation $\hat{A}_0\approx \mathcal{U}\Sigma\mathcal{V}^*$, then the generalized eigenvalue problem $\mathcal{U}^*\hat{A}_1\mathcal{V}X=\Sigma X\Lambda$ gives the eigenvalues of $T(\lambda)$ inside $\Omega$ up to error arising from the linear solves, quadrature rule and oversampling $m$.}\footnote{{This follows by Keldysh's theorem. Assume that we have $A_0$ and $A_1$ exactly, and that we know $m$. Then by \cref{eqn:cauchykeldysh}, $A_0= AB^*$ and $A_1= AJB^*$ for quasimatrices $A$ and $B$. Hence the eigenvalue problem $\mathcal{U}^*{A}_1\mathcal{V}X=\Sigma X\Lambda$ can be written as $(\mathcal{U}^*AJB^*\mathcal{V})X=(\mathcal{U}^*AB^*\mathcal{V})X\Lambda$ and so has the same eigenvalues as $J$.}} The algorithm is summarized in \cref{surfBEYN}. {While the residuals $\|T(\lambda)f\|$ do not bound the distance to the spectrum as in \eqref{eqn:residual_defn}, we can use pseudospectral inclusions to verify computations; the pseudospectrum for nonlinear eigenvalue problems is defined as in \cref{eqn:pspec_defn} but with $(z\mathcal{I}-\mathcal{L})^{-1}$ replaced by $T(z)^{-1}$, and analogous results for its construction in terms of residuals and convergence to the spectrum as $\epsilon\downarrow 0$ hold \cite{colbrook2025avoiding}.}

\begin{algorithm}[t]
\caption{\texttt{SurfBeyn}: Computing eigenvalues of a nonlinear surface differential operator.}
\label{surfBEYN}
\textbf{Input:} Nonlinear operator $T$, quadrature weights $\{w_j\}_{j=1}^\ell$ and nodes $\{z_j\}_{j=1}^\ell$, test functions $F=[f_1\cdots f_{m+m_0}]$, $\epsilon_{\mathrm{res}},\epsilon_{\mathrm{tol}}>0$.\!\!\!\!\!\!\!\!\!\!\!\!\!\!\!\!\!\!\!\!\!\!\!\!\!\!\!\!\!\!
\begin{algorithmic}[1]
\REPEAT
\STATE{Solve $T(z_k)G_k=F$ for $k=1,\ldots,\ell$ with relevant boundary conditions.}
\STATE{Compute $\hat{A}_0= (2\pi i)^{-1}\sum_{k=1}^\ell w_k\hat{G}_k$ and $\hat{A}_1= (2\pi i)^{-1}\sum_{k=1}^\ell w_kz_k\hat{G}_k$.}
\STATE{Compute a truncated SVD of $\hat{A}_0\approx \mathcal{U}\Sigma\mathcal{V}^*$, excluding singular values $<\epsilon_{\mathrm{tol}}$.}
\STATE{Form $B=\mathcal{U}^*\hat{A}_1\mathcal{V}\Sigma^{-1}$ and solve $BX=X\Lambda$ for $X$ and $\Lambda$.}
\STATE{Set $U=\mathcal{U}X$, normalize the columns of $U$, and update $F=U$.}
\UNTIL{$\|T(\lambda_j)u_j\|\leq \epsilon_{\mathrm{res}}\max\{1,|\lambda_j|\}$ for each computed eigenpair $(\lambda_j,u_j)$.}
\end{algorithmic}
\textbf{Output:} Eigenvalues $\Lambda$, eigenfunctions $\mathcal{U} X$.
\end{algorithm}

\subsubsection{Example}

We consider an example inspired by \cite{nlep_coll}, arising from the partial delay differential equation (PDDE) with Dirichlet boundary conditions:
\begin{align*}\setlength\abovedisplayskip{6pt}\setlength\belowdisplayskip{6pt}
u_t(\mathbf{x},t)&=\Delta u(\mathbf{x},t)+a(\mathbf{x})u(\mathbf{x},t)+b(\mathbf{x})u(\mathbf{x},t-1), \quad &\mathbf{x}\in M, t\geq 0,\\
u(\mathbf{x},t)&=0, \quad &\mathbf{x}\in\partial M, t\geq 0.
\end{align*} 
Taking a Laplace transform with respect to the time variable, we obtain a NLEP of the form \cref{nlep} with
\begin{equation}\setlength\abovedisplayskip{0pt}\setlength\belowdisplayskip{6pt}\label{eqn:nlep_example}
T(\lambda)=\Delta+(a(\mathbf{x})-\lambda+e^{-\lambda}b(\mathbf{x}))\mathcal{I}.
\end{equation}
We consider this problem on the M{\"o}bius strip, with $a(\mathbf{x})=\sin^2(x)\sin^2(y)\sin^2(z)$ and $b(\mathbf{x})=\sin(x+y+z)+1.5$, and apply Dirichlet boundary conditions. {The M\"obius strip mesh used has order $p=16$ and $210$ mesh elements.}

\begin{figure}
    \centering
    \begin{tabular}{ c  c  c }
    Eigenfunction $1$&Eigenfunction $10$&Eigenfunction $100$\\
    \includegraphics[width=0.25\textwidth]{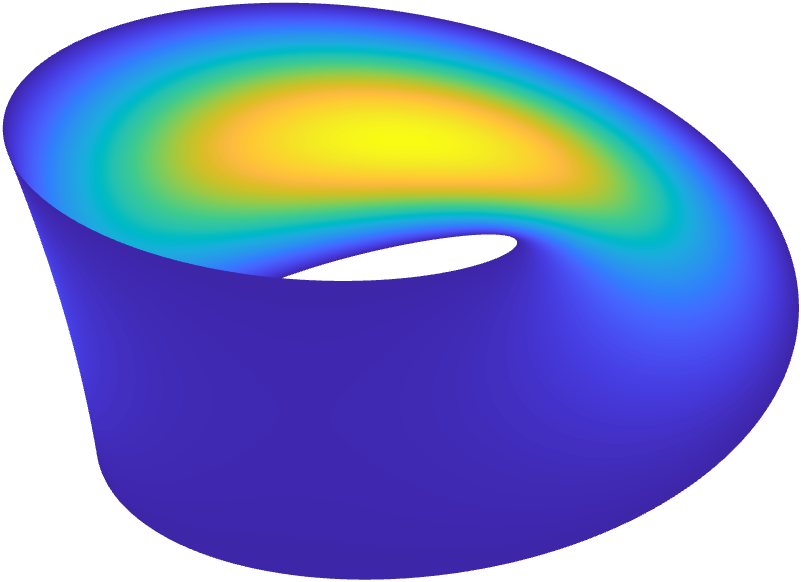}
    &\includegraphics[width=0.25\textwidth]{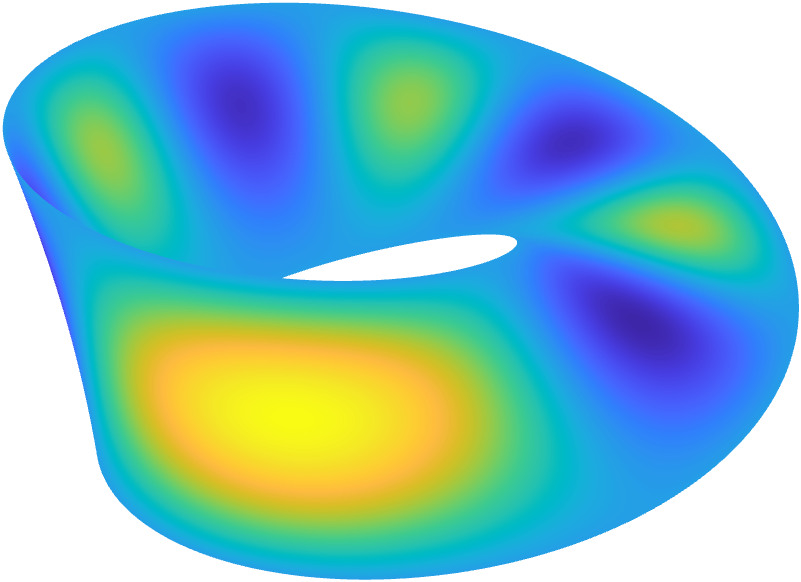}&
    \includegraphics[width=0.25\textwidth]{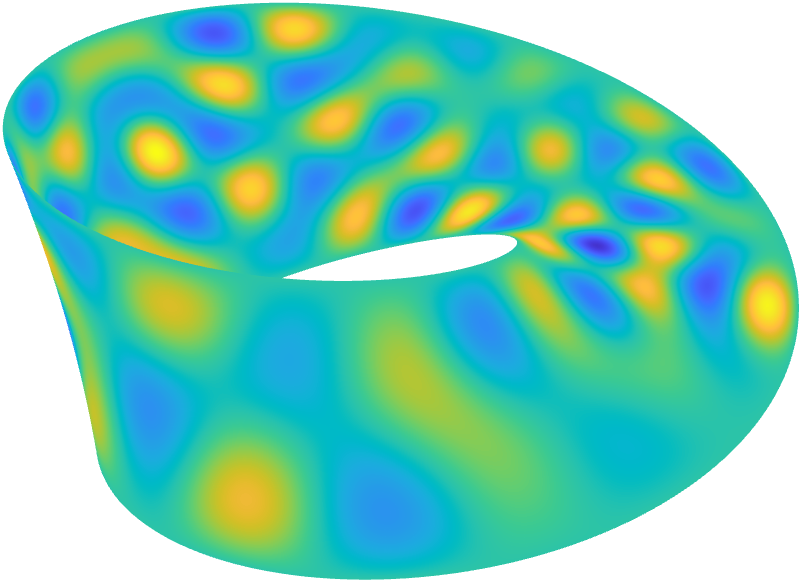}
\end{tabular}
    \label{fig:nlep_efuns}
    \caption{Eigenfunctions of the NLEP \cref{eqn:nlep_example} {on the M\"obius strip mesh with $p=16$ and $210$ mesh elements using $l=50$ quadrature nodes}.}
\end{figure}

As before, we first use elliptical contours to localize the eigenvalues, followed by small, isolated, circular contours to compute them accurately. The eigenvalues are shown in \cref{fig:nlep} along with the corresponding residuals, which allow for interpretations in terms of pseudospectral inclusions and convergence to the true spectrum as in \cref{errctrl}. Residuals remain small, with the expected drop in accuracy for larger eigenvalues as observed in \cref{sec:triceratops}. \Cref{fig:nlep_efuns} displays the 1st, 10th, and 100th eigenfunctions, showing the expected increase in oscillations with eigenvalue magnitude.

\begin{figure}
    \centering
    \includegraphics[width=0.46\linewidth]{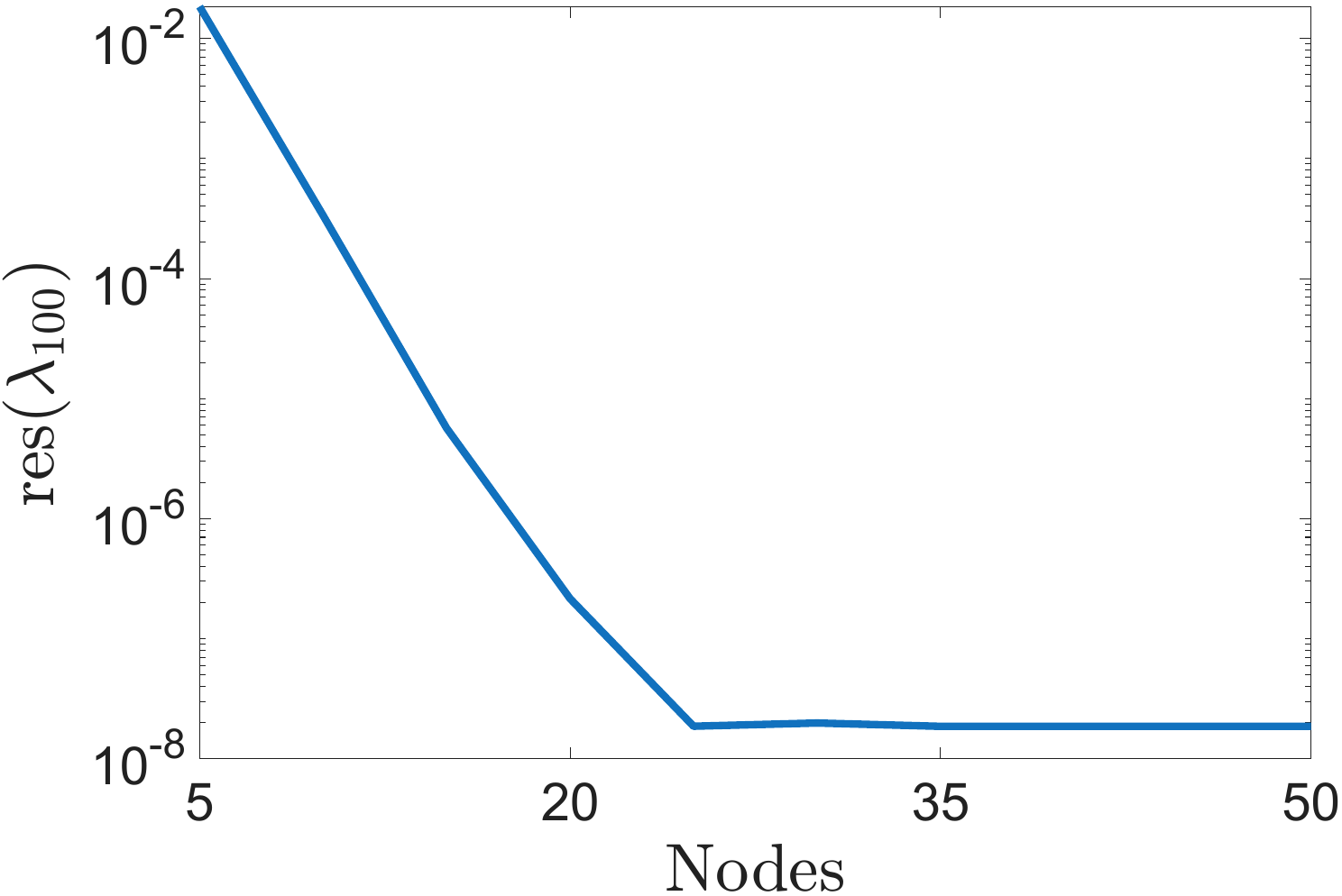}\hfill
        \includegraphics[width=0.46\linewidth]{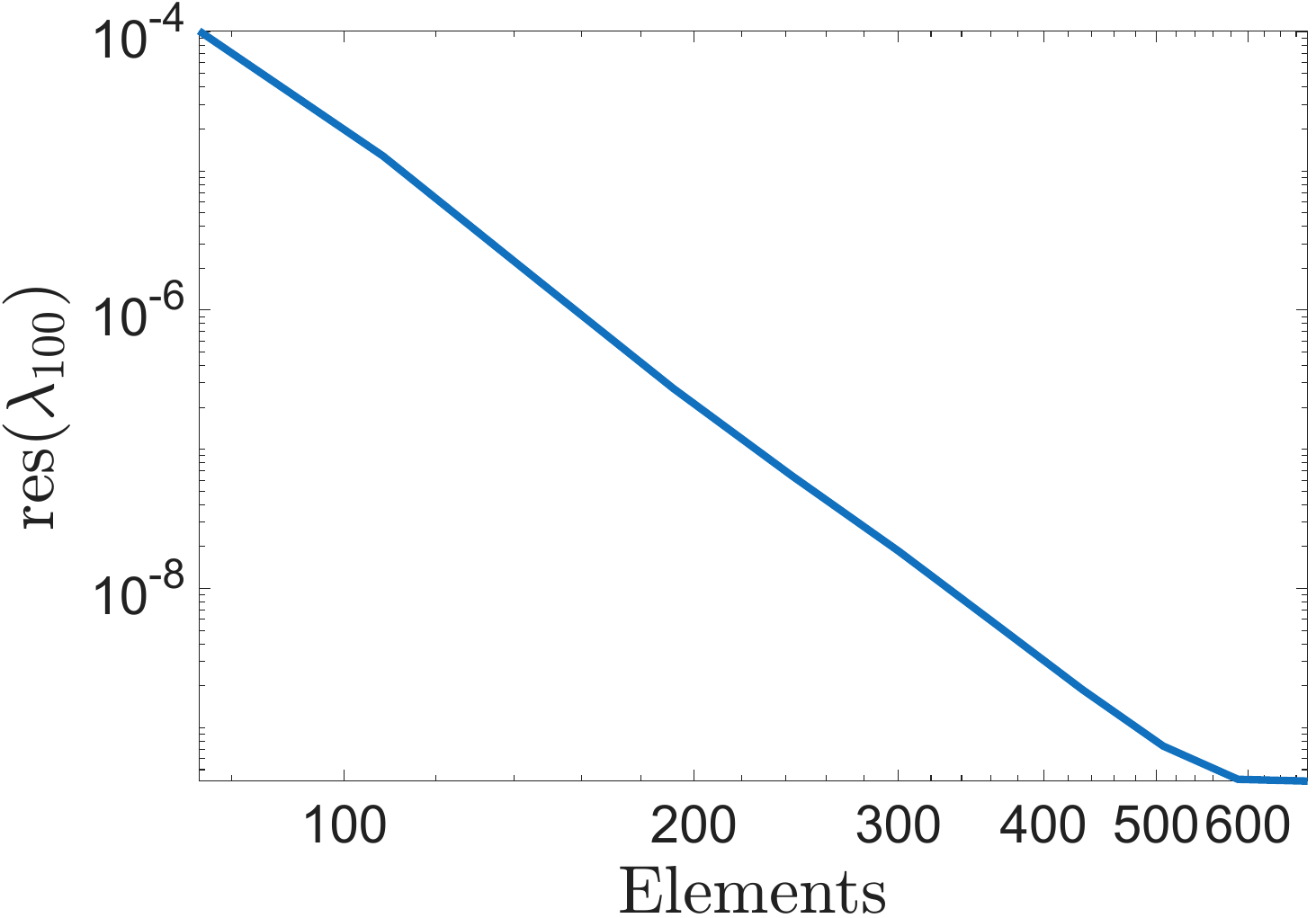}
    \caption{Left: The residual of the $100\text{th}$ eigenvalue of the NLEP \cref{eqn:nlep_example} on a M\"obius strip solved using different numbers of quadrature nodes. Right: The residual of the $100\text{th}$ eigenvalue of the same problem solved using varying numbers of mesh elements.}
    \label{fig:nlep_convergence}
\end{figure}

We also analyze convergence of the algorithm. First we vary the number of quadrature nodes while keeping the mesh fixed; {on the left of \cref{fig:nlep_convergence} we see} that while more nodes are needed for convergence than for \texttt{SurfFEAST} {as in \cref{fig:schr_res_order}}, still only around 25 nodes are required (here we have no symmetry consideration to exploit). To explore how refining the mesh affects convergence, we vary the number of mesh elements rather than the order of the mesh as in \cref{sec:triceratops} (recall \cref{fig:mesh_refinement}), using $50$ nodes throughout. We see {on the right of \cref{fig:nlep_convergence}} that the residual varies linearly (on the log-log plot) with the number of elements, suggesting polynomial convergence.

\begin{table}[t]
    \centering
    \begin{tabular}{|c|c|c|}
    \hline
    Number & Eigenvalue & Residual \\
    \hline
    1 & $-1.359074970$ & $6.10\times 10^{-11}$\\
    2 & $-1.59561416$ & $1.03\times 10^{-9}$\\
    3 & $-1.85101635$ & $6.69\times 10^{-11}$\\
    4 & $-2.089752672$ & $5.95\times 10^{-11}$\\
    5 & $-2.291412643$ & $5.94\times 10^{-11}$\\
    10 & $-2.835581355$ & $9.27\times 10^{-11}$\\
    20 & $-3.55705229$ & $8.89\times 10^{-11}$\\
    50 & $-4.31588030$ & $3.96\times 10^{-10}$\\
    100 & $-5.0004860$ & $1.40\times 10^{-8}$\\
    \hline
    \end{tabular}
    \caption{Eigenvalues and corresponding residuals {(for pseudospectral inclusions)} for the NLEP \cref{eqn:nlep_example}.}
    \label{fig:nlep}
\end{table}

\section{Spectral measures on surfaces}\label{sec:surfspecmeas}

Many infinite-dimensional linear operators in physics and related fields have spectra that are not purely discrete, so their eigenfunctions do not form a basis of the Hilbert space. In particular, an operator $\mathcal{L}$ on $\mathcal{H}$ may have non-empty \textit{continuous spectrum}. However, the spectral theorem for self-adjoint operators provides a resolution: projections onto eigenspaces are replaced by a projection-valued measure $\mathcal{E}$~\cite[Thm.~VIII.6]{reed1972methods}, which assigns an orthogonal projection to each Borel-measurable set such that
\[\setlength\abovedisplayskip{6pt}\setlength\belowdisplayskip{6pt}
f=\left(\int_\mathbb{R} \dd\mathcal{E}(\lambda)\right)f,\quad \forall f\in \mathcal{H}, \qquad\text{and}\qquad \mathcal{L}f=\left(\int_\mathbb{R} \lambda\dd\mathcal{E}(\lambda)\right)f, \quad \forall f\in\mathcal{D}(\mathcal{L}).
\]
Thus $\mathcal{E}$ decomposes $\mathcal{H}$ and diagonalizes $\mathcal{L}$. The spectral measure of $\mathcal{L}$ with respect to $f\in\mathcal{H}$ is a scalar measure defined by $\mu_f(S)=\langle\mathcal{E}(S)f,f\rangle$, where $S\subset\mathbb{R}$ is a Borel-measurable set. Such measures are central in applications ranging from quantum mechanics \cite{reed1972methods,hall2013quantum} and plasma physics \cite{wilkening2015spectral} to stochastic processes \cite{kallianpur1971spectral,girardin2003semigroup} and signal processing \cite{stoica2005spectral}. A similar spectral theorem holds for normal operators.


{For a self-adjoint matrix, the spectral measure will consist of Dirac delta functions located at the eigenvalues of the matrix along the real line, with weights depending on the function $f$. However, infinite-dimensional operators can have continuous spectrum, as in the left of \cref{fig:torus_spec_meas} and the right of \cref{fig:bump_meas}. }When $\mathcal{L}$ has a non-empty continuous spectrum, we construct smoothed approximations of $\mu_f$. In particular, we compute a smooth function $\mu^\epsilon_f$, with parameter $\epsilon>0$, that converges weakly to $\mu_f$, i.e., for any bounded, continuous function $\phi$~\cite[Ch.~1]{billingsley2013convergence},
\[
\int_\mathbb{R} \phi(\lambda)\mu^\epsilon_f(\lambda)\dd \lambda\rightarrow \int_\mathbb{R}\phi(\lambda)\dd\mu_f(\lambda), \qquad\text{as}\qquad\epsilon\downarrow 0.
\]
To achieve this, we again use the resolvent $(\mathcal{L}-z\mathcal{I})^{-1}$. {In particular, an argument} using Stone's formula \cite{stone1932linear} gives a smooth approximation of $\mu_f$ in terms of a convolution with the Poisson kernel:
\begin{equation}\setlength\abovedisplayskip{6pt}\setlength\belowdisplayskip{6pt}
\label{stones_formula}
\mu_f^{\epsilon}(\lambda){\coloneqq}-\cfrac{1}{\pi}\im\left(\langle(\mathcal{L}-(\lambda-i\epsilon)\mathcal{I})^{-1}f,f\rangle\right)=\cfrac{1}{\pi}\int_{\mathbb{R}}\cfrac{\epsilon \dd\mu_f(y)}{\epsilon^2+(\lambda-y)^2}
\end{equation}
The function $\mu_f^{\epsilon}$ converges weakly to the spectral measure $\mu_f$ as $\epsilon\downarrow 0$. Moreover, if $\mu_f$ has sufficient local regularity around the point $\lambda$, then $\mu_f^{\epsilon}(\lambda)$ converges to the Radon--Nikodym derivative $\rho_f(\lambda)$ as $\epsilon\downarrow 0$. However, for small $\epsilon$, solving the shifted linear systems to evaluate the resolvent becomes computationally expensive, as we need a large discretization size (see \eqref{eqn:erb_inverseproblem}). Moreover, the approximation error of \eqref{stones_formula} in $\epsilon$ is $\mathcal{O}(\epsilon\log(1/\epsilon))$. These two properties limit the practicality of the Poisson kernel.

\subsection{High-order smoothing kernels}

Instead of the Poisson kernel, using higher order kernels in the above analysis improves convergence; in particular, since {we can use} larger $\epsilon$ for a specified accuracy, we may use smaller discretization sizes and hence have faster solves. We define
\[\setlength\abovedisplayskip{6pt}\setlength\belowdisplayskip{6pt}
K(x)=\cfrac{1}{2\pi i}\sum_{j=1}^m\left[\cfrac{\alpha_j}{x-a_j}-\cfrac{\overline{\alpha_j}}{x-\overline{a_j}}\right], \quad \text{and}\quad K_{\epsilon}(x)=K(x/\epsilon)/\epsilon,
\]
where the $\{a_j\}_{j=1}^m$ are distinct points in the upper half plane. Following \cite{specsolve}, we use equispaced poles $a_j={2j}/({m+1})-1+i$. We choose the $\alpha_j$ to satisfy the Vandermonde system
\begin{equation}\setlength\abovedisplayskip{6pt}\setlength\belowdisplayskip{6pt}
\label{vandermonde}
\begin{pmatrix}
1&\dots&1\\
a_1&\dots&a_m\\
\vdots&\ddots&\vdots\\
a_1^{m-1}&\dots&a_m^{m-1}
\end{pmatrix}
\begin{pmatrix}
\alpha_1\\
\alpha_2\\
\vdots\\
\alpha_m
\end{pmatrix}=
\begin{pmatrix}
1\\
0\\
\vdots\\
0
\end{pmatrix}.
\end{equation}
This choice ensures that the kernel is an $m$th order kernel. Namely, $K$ integrates to one, its first to $(m-1)$th moments vanish and $K(x)=\mathcal{O}(1/(1+|x|)^{m+1})$ as $|x|\rightarrow\infty$. Since $K$ is a \textit{rational} kernel, we can obtain a new approximation for the spectral measure directly through the resolvent as
\[\setlength\abovedisplayskip{6pt}\setlength\belowdisplayskip{6pt}
\mu_f^{\epsilon}(\lambda)=[K_{\epsilon}\ast \mu_f](\lambda)=-\cfrac{1}{\pi}\im\left(\sum_{j=1}^m\alpha_j\langle(\mathcal{L}-(\lambda-\epsilon a_j)\mathcal{I})^{-1}f,f\rangle\right).
\]
This approximation generalizes \eqref{stones_formula} and can be thought of as an $m$th order version of Stone's formula. Again, this approximation converges weakly to $\mu_f$. However, now with sufficient regularity conditions on $\mu_f$ it converges to $\rho_f$, both pointwise and in the $L^p$ sense, with a rate of convergence $O(\epsilon^m\log(1/\epsilon))$: a significant increase \cite{specsolve}.

In terms of practical implementation, the only requirements are solving the shifted linear systems $(\mathcal{L}-(\lambda-\epsilon a_j)\mathcal{I})g_j=f$ for $j=1,\ldots,m$ and taking inner products of functions on surfaces, both of which have been discussed previously and can be done within the framework of \texttt{surfacefun}. This leads to \cref{surfSpec} for computing spectral measures of operators on surfaces, which we call \texttt{SurfMeas}.

\begin{algorithm}[t]
\caption{\texttt{SurfMeas}: Computing spectral measures of self-adjoint differential operators on surfaces.}
\label{surfSpec}
\begin{algorithmic}[1]
\STATE{\textbf{Input:} $\mathcal{L}:\mathcal{D}(\mathcal{L})\rightarrow\mathcal{H}$ self-adjoint, $f\in\mathcal{H}$, $\lambda\in\mathbb{R}$, $a_1,\ldots,a_m\in\{z\in\mathbb{C}:\im(z)>0\}$ and $\epsilon>0$.}
\STATE{Solve the Vandermonde system \cref{vandermonde} for the residues $\alpha_j$, $1\leq j\leq m$.}
\STATE{Solve $(\mathcal{L}-(\lambda-\epsilon a_j)\mathcal{I})g_j=f$ for $1\leq j\leq m$ for solutions $\hat{g}_j$.}
\STATE{Compute $\mu_f^{\epsilon}(\lambda)=-\im\left(\sum_{j=1}^m\alpha_j\langle \hat{g}_j,f\rangle\right)/\pi$.}
\STATE{\textbf{Output:} $\mu_f^{\epsilon}(\lambda)$.}
\end{algorithmic}
\end{algorithm}

\subsection{Examples}

We consider two examples. The first is on a bounded surface, whereas the second is on an unbounded surface.

\subsubsection{Internal waves}

A common source of continuous spectra of linear operators on compact surfaces is linear pencil problems. We compute the spectra of the linear pencil $\mathcal{A}-\lambda\mathcal{B}$ for $\mathcal{A}$ and $\mathcal{B}$ self-adjoint and $\mathcal{B}$ positive and invertible. We study this through the operator $\mathcal{B}^{-1}\mathcal{A}$. We define the Hilbert space $\mathcal{H}_{\mathcal{B}}$ by equipping $\mathcal{D}(\mathcal{B}^{1/2})$ with the norm $\|f\|_{\mathcal{B}}=\sqrt{\langle\mathcal{B}^{1/2}f,\mathcal{B}^{1/2}f\rangle}$ and define a symmetric closed operator 
$
\mathcal{L}=\overline{\mathcal{B}^{-1}\mathcal{A}|_{\mathcal{D}(\mathcal{A})\cap\mathcal{D}(\mathcal{B}^{1/2})}}$.
Under suitable assumptions on $\mathcal{A}$ and $\mathcal{B}$, $\mathcal{L}$ is self-adjoint \cite{colbrook2022specsolve}, and so we apply the approach of \cref{surfSpec} to compute its spectral measures. For simplicity, we assume that $f\in\mathcal{D}(\mathcal{A})\cap\mathcal{D}(\mathcal{B}^{1/2})$ which can be a strict (but dense) subspace of $\mathcal{H}_{\mathcal{B}}$. Our expression for the convolution now becomes 
\[\setlength\abovedisplayskip{6pt}\setlength\belowdisplayskip{6pt}
\mu_f^{\epsilon}(\lambda)=[K_{\epsilon}\ast\mu](x)=-\cfrac{1}{\pi}\im\left(\sum_{j=1}^m\alpha_j\left\langle(\mathcal{A}-(\lambda-\epsilon a_j)\mathcal{B})^{-1}\mathcal{B}f,\mathcal{B}f\right\rangle\right),
\]
so we must solve $(\mathcal{A}-(\lambda-\epsilon a_j)\mathcal{B})g_j=\mathcal{B}f$ for $j=1,\ldots,m$.

\begin{figure}[t]
\centering
\raisebox{-1\height}{\includegraphics[width=0.45\textwidth,trim={0mm 0mm 0mm 0mm},clip]{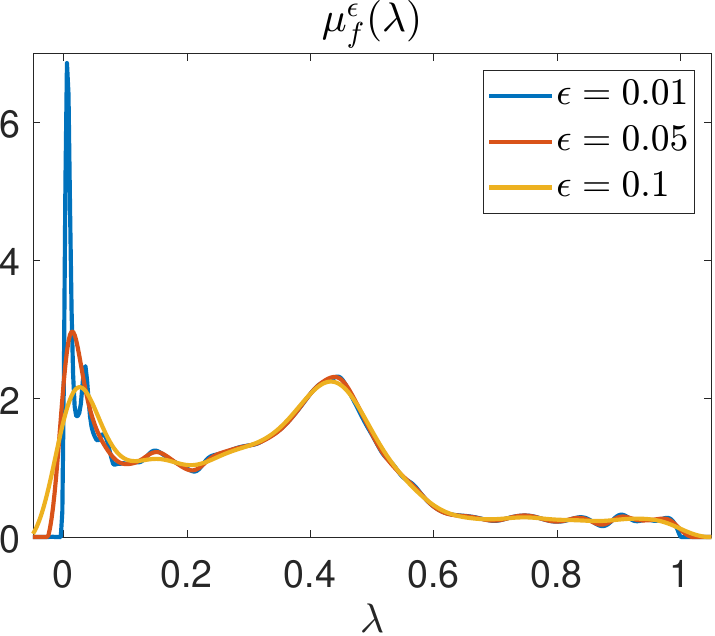}}
\hfill
\raisebox{-1\height}{\includegraphics[width=0.52\textwidth,trim={0mm 0mm 0mm 0mm},clip]{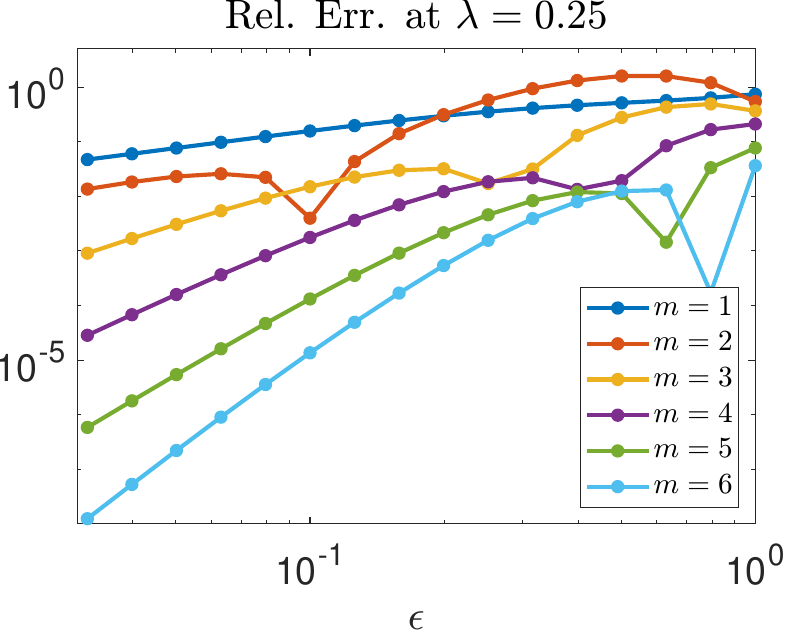}}
\caption{Left: Smoothed spectral measures for different choices of $\epsilon$ computed using the sixth-order kernel. As $\epsilon\downarrow 0$ these approximations converge weakly to the spectral measures. Right: Convergence to the Radon--Nikodym derivative at $\lambda=0.25$. The errors are computed by comparing to a value computed at $\epsilon=0.01$. Note the vastly improved convergence rates for the high-order kernels, with the $m$th order kernel having $m$th order convergence, up to logarithmic factors.}
\label{fig:torus_spec_meas}
\end{figure}

We consider an example inspired by internal waves \cite{dyatlov2024mathematics} and take
\[\setlength\abovedisplayskip{6pt}\setlength\belowdisplayskip{6pt}
\mathcal{A}=-\partial_{yy}+\cos^2(10x),\quad\mathcal{B}=-\Delta +\mathcal{I}
\]
on the torus. {Here $\mathcal{H}_{\mathcal{B}}=H^1(M)$, the Sobolev space of square-integrable functions with square-integrable first (weak) derivatives.} The spectra of internal waves are a significant topic in oceanography and the study of rotating fluids \cite{maas2005wave,sibgatullin2019internal}. In this example, the spectrum of $\mathcal{L}$ is $[0,1]$, but the spectral type is unknown \cite{ralston1973stationary}. As our test function, we take
\[\setlength\abovedisplayskip{6pt}\setlength\belowdisplayskip{6pt}
f(x,y,z)=\sin(2(x-y+z)^2).
\]
\cref{fig:torus_spec_meas} (left) shows the smoothed spectral measures computed using the sixth order kernel for various $\epsilon$. As $\epsilon\downarrow 0$, the approximations appear to converge to a density function, suggesting the spectrum of this operator is absolutely continuous. \cref{fig:torus_spec_meas} (right) shows the convergence to the Radon--Nikodym derivative at the point $\lambda=0.25$ and the improved convergence rates of the higher order kernels.

\begin{figure}
    \centering
    \setlength{\tabcolsep}{6pt}
    \begin{tabular}{ c  c  c }
    $u_f^{\epsilon},\;\epsilon=0.1$&$u_f^{\epsilon},\;\epsilon=0.05$&$u_f^{\epsilon},\;\epsilon=0.01$\\
    \includegraphics[width=0.3\textwidth]{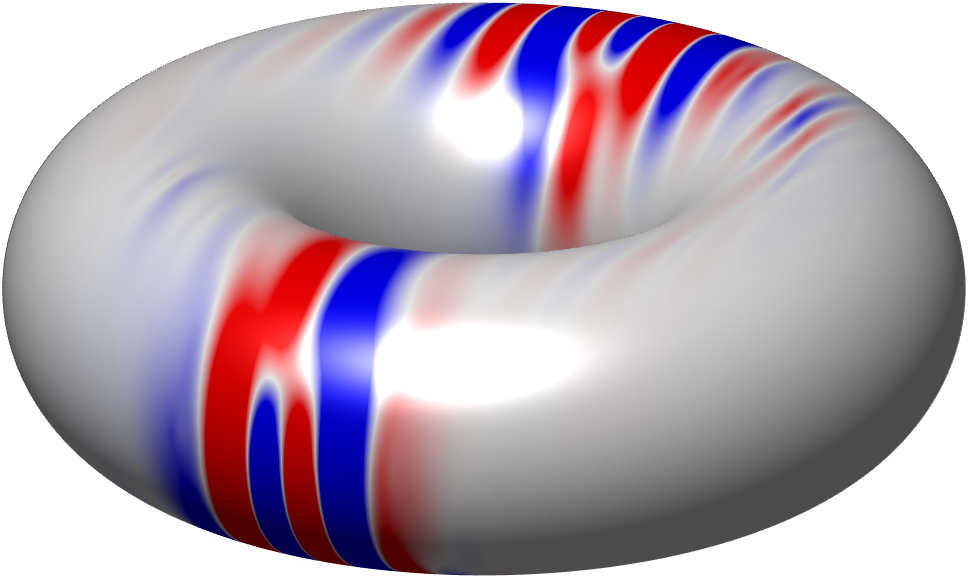}
    &\includegraphics[width=0.3\textwidth]{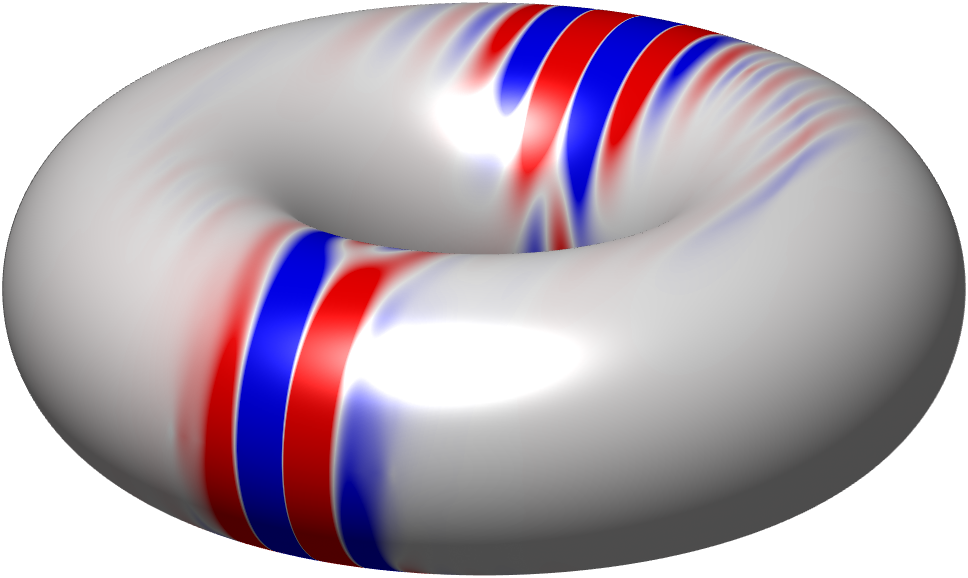}&
    \includegraphics[width=0.3\textwidth]{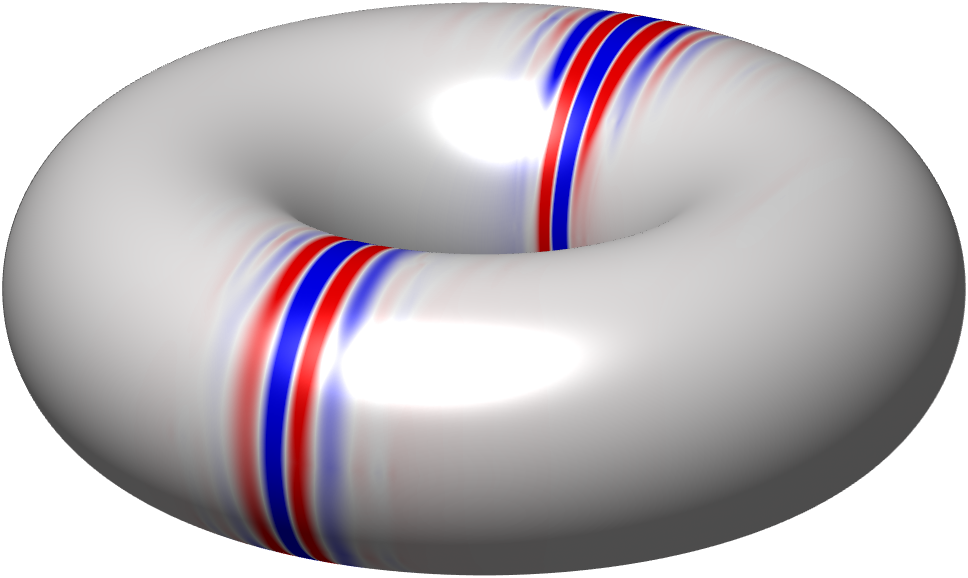}
\end{tabular}
    \label{fig:torus_efuns}
    \caption{The functions in \eqref{surf_gen_efuns} for $\lambda=1$ and various $\epsilon$. As $\epsilon\downarrow 0$, we observe increased oscillations and localization corresponding to a singular generalized eigenfunction. {The mesh has order $p=24$ and $768$ mesh elements.}}
\end{figure}

We can also consider the functions
\begin{equation}\setlength\abovedisplayskip{6pt}\setlength\belowdisplayskip{6pt}
\label{surf_gen_efuns}
u_f^{\epsilon}(x,y,z)=-\frac{1}{\pi}\im\left(\sum_{j=1}^m\alpha_jg_j(x,y,z)\right),
\end{equation}
which correspond to the smoothed spectral projection density applied to \(f\). Under suitable conditions, these functions converge in the sense of distributions to generalized eigenfunctions of the operator \cite{colbrook2025computing}. \cref{fig:torus_efuns} shows these functions for $\lambda=1$ and various $\epsilon$, using a torus mesh with $768$ elements and order $24$. As $\epsilon\downarrow 0$, we observe increasing oscillations and localization of these functions around a narrow band of the torus. In particular, while approximations converge in a weak sense, they do not converge in $\mathcal{H}_{\mathcal{B}}$.

\subsubsection{Schr\"odinger operator on unbounded surface}

\begin{figure}[t]
\centering
\raisebox{-1\height}{\includegraphics[width=0.52\textwidth,trim={0mm 0mm 0mm 0mm},clip]{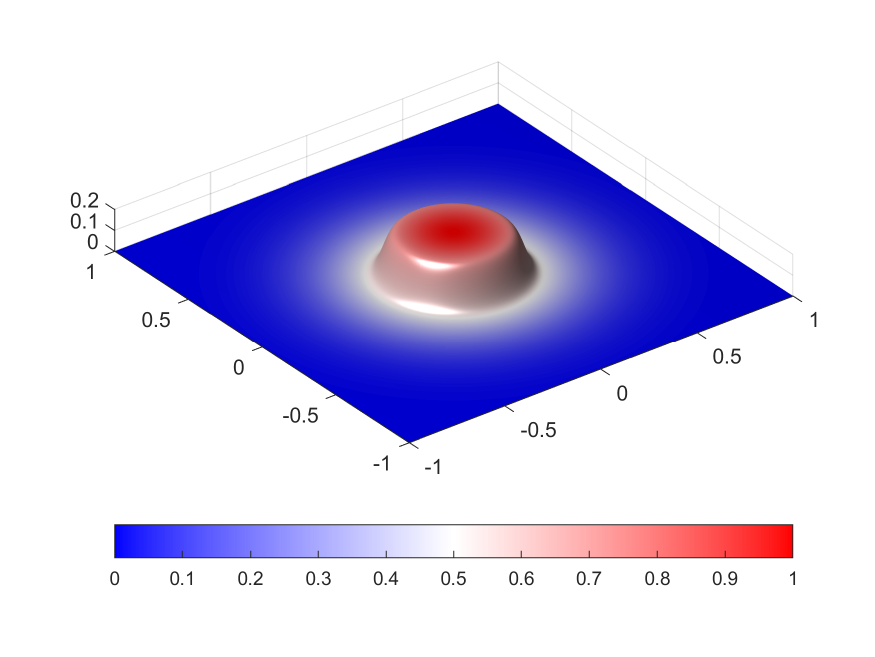}}
\hfill
\raisebox{-1\height}{\includegraphics[width=0.45\textwidth,trim={0mm 0mm 0mm 0mm},clip]{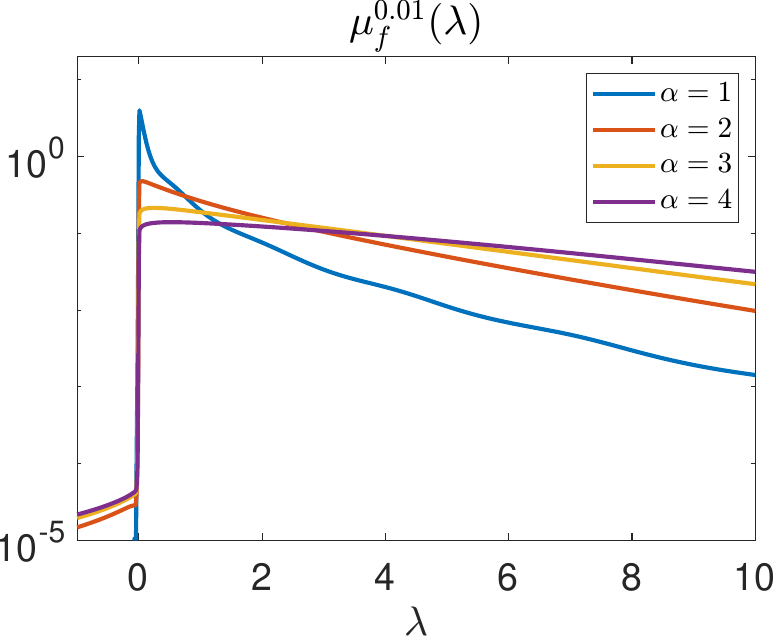}}
\caption{Left: Potential $V(x,y,z)=\exp(-5(x^2+y^2))$ on the infinite surface. Right: Spectral measures for different $\alpha$. Larger $\alpha$ leads to more decay in physical space but less in spectral space.}
\label{fig:bump_meas}
\end{figure}

{Next, we consider a Schr\"odinger operator $\mathcal{L}$ on an unbounded surface $M$. The surface $M = M_\text{int} \cup M_\text{ext}$ is flat outside a compact region $M_\text{int}$ and is shown along with the potential $V(x,y,z)=\chi_{M_\text{int}}(x,y,z)\exp(-5(x^2+y^2))$ in \cref{fig:bump_meas} (left), where $\chi_{M_\text{int}}(x,y,z)$ is the indicator function of the finite region. The infinite-surface Schr\"odinger problem on $M$ may be reformulated into two coupled problems: a bounded-surface, variable-coefficient problem on $M_\text{int}$ and an unbounded, flat, constant-coefficient problem on $M_\text{ext} = M \setminus M_\text{int}$. The latter problem can be solved using potential theory. We construct the surface operator on $M_\text{int}$ using a variant of the method described in \cref{surfacefun} that uses impedance-to-impedance operators to avoid spurious resonances~\cite{Gillman2015}. The boundary integral operators are constructed using the \texttt{chunkIE} package~\cite{chunkie}.}

For \(\alpha>0\), we use the test functions
\[\setlength\abovedisplayskip{6pt}\setlength\belowdisplayskip{6pt}
f(x,y,z)=\frac{c_\alpha \chi_{M_\text{int}}(x,y,z)}{[(x+0.2)^2+(y+0.1)^2+1]^{\alpha}},
\]
where $c_\alpha$ is a normalization constant so that $\|f\|=1$. \cref{fig:bump_meas} (right) shows the spectral measures computed using $\epsilon=0.01$ and the 6th order kernel. The spectrum in this case is the whole positive half-line $[0,+\infty)$. We see that a larger value of $\alpha$ leads to slower decay of the spectral measure, just as in the case of the free Laplacian on $\mathbb{R}^2$. We also observe the excellent localization of the smooth measure to $\mathbb{R}_{\geq 0}$.

\section{Conclusions} The \texttt{SurfSpec} methods developed in this paper provide robust, accurate algorithms for computing spectral properties of differential operators on surfaces; the methods are based on the resolvent and form part of a ``solve-then-discretize'' framework. The algorithms are flexible, come with error control and convergence guarantees and can be applied to a wide variety of surfaces, PDEs, and spectral problems. The power of these algorithms has been demonstrated via examples of Laplace--Beltrami and Schr{\"o}dinger operators on complicated surfaces, such as the triceratops or an unbounded manifold, and internal waves on the torus.

\texttt{SurfSpec} can find application in diverse fields of mathematics. Beyond applications such as spectral shape analysis or data representation, it can be used to inform and prove conjectures in spectral geometry. For example, by combining with numerical optimization algorithms, we might test and formulate conjectures in shape optimization and spectral theory \cite{delfour2011shapes,henrot2018shape}. The framework for computing spectral measures using high-order kernels can be extended to computing the vector-valued functional calculus of operators on manifolds \cite[Sec.~VII.1]{reed1972methods}. This in turn can be used to solve linear evolutionary PDEs \cite{delaubenfels2006existence}; more complicated evolutionary PDEs can be solved by combining this with splitting methods \cite{McLachlan_Quispel_2002}. The contour integral methods developed in this paper could also be utilized for system realization and control for dynamics on surfaces \cite{brennan2023contour}.

\bibliographystyle{jabbrv_siam}
\bibliography{paperreferences}

@inproceedings{reuter2005laplace,
  title={Laplace-spectra as fingerprints for shape matching},
  author={Reuter, Martin and Wolter, Franz-Erich and Peinecke, Niklas},
  booktitle={Proceedings of the 2005 ACM symposium on Solid and Physical Modeling},
  pages={101--106},
  year={2005}
}

@inproceedings{seo2011laplace,
  title={Laplace--{B}eltrami eigenfunction expansion of cortical manifolds},
  author={Seo, Seongho and Chung, Moo K},
  booktitle={2011 IEEE International Symposium on Biomedical Imaging: From Nano to Macro},
  pages={372--375},
  year={2011},
  organization={IEEE}
}

@inproceedings{zhang2007spectral,
  title={Spectral methods for mesh processing and analysis},
  author={Zhang, Hao and van Kaick, Oliver and Dyer, Ramsay},
  booktitle={Proceedings of Eurographics State-of-the-art Report},
  volume={122},
  pages={1865--1894},
  year={2007},
  organization={Eurographics Association Prague}
}

@article{peter2014feast,
  title={{FEAST} as a subspace iteration eigensolver accelerated by approximate spectral projection},
  author={Tang, Peter Ping Tak and Polizzi, Eric},
  journal={SIAM Journal on Matrix Analysis and Applications},
  volume={35},
  number={2},
  pages={354--390},
  year={2014},
  publisher={SIAM}
}

@article{nakatsukasa2020sharp,
  title={Sharp error bounds for {R}itz vectors and approximate singular vectors},
  author={Nakatsukasa, Yuji},
  journal={Mathematics of Computation},
  volume={89},
  number={324},
  pages={1843--1866},
  year={2020}
}

@article{caissard2019laplace,
  title={Laplace--{B}eltrami operator on digital surfaces},
  author={Caissard, Thomas and Coeurjolly, David and Lachaud, Jacques-Olivier and Roussillon, Tristan},
  journal={Journal of Mathematical Imaging and Vision},
  volume={61},
  pages={359--379},
  year={2019},
  publisher={Springer}
}

@article{tan2014spectral,
  title={Spectral {L}aplace--{B}eltrami wavelets with applications in medical images},
  author={Tan, Mingzhen and Qiu, Anqi},
  journal={IEEE Transactions on Medical Imaging},
  volume={34},
  number={5},
  pages={1005--1017},
  year={2014},
  publisher={IEEE}
}

@article{lewin2010spectral,
  title={Spectral pollution and how to avoid it},
  author={Lewin, Mathieu and S{\'e}r{\'e}, {\'E}ric},
  journal={Proceedings of the London mathematical society},
  volume={100},
  number={3},
  pages={864--900},
  year={2010},
  publisher={Oxford University Press}
}

@article{webb2021spectra,
  title={Spectra of {Jacobi} operators via connection coefficient matrices},
  author={Webb, Marcus and Olver, Sheehan},
  journal={Communications in Mathematical Physics},
  volume={382},
  number={2},
  pages={657--707},
  year={2021},
  publisher={Springer}
}

@article{gilles2019continuous,
  title={Continuous analogues of {K}rylov subspace methods for differential operators},
  author={Gilles, Marc Aur{\`e}le and Townsend, Alex},
  journal={SIAM Journal on Numerical Analysis},
  volume={57},
  number={2},
  pages={899--924},
  year={2019},
  publisher={SIAM}
}

@article{olver2009gmres,
  title={{GMRES} for the differentiation operator},
  author={Olver, Sheehan},
  journal={SIAM Journal on Numerical Analysis},
  volume={47},
  number={5},
  pages={3359--3373},
  year={2009},
  publisher={SIAM}
}

@inproceedings{vainikko2004gmres,
  title={{GMRES} and discrete approximation of operators},
  author={Vainikko, Gennadi},
  booktitle={Proc. Est. Acad. Sci. Phys. Math.},
  volume={53},
  pages={124--131},
  year={2004}
  }

@article{colbrook2019infinite,
  title={On the infinite-dimensional {QR} algorithm},
  author={Colbrook, Matthew J and Hansen, Anders C},
  journal={Numerische Mathematik},
  volume={143},
  pages={17--83},
  year={2019},
  publisher={Springer}
}

@phdthesis{webb2017isospectral,
  title={Isospectral algorithms, {T}oeplitz matrices and orthogonal polynomials},
  author={Webb, Marcus},
  school={University of Cambridge},
  year={2017}
}

@article{colbrook2022computing,
  title={Computing semigroups with error control},
  author={Colbrook, Matthew J},
  journal={SIAM Journal on Numerical Analysis},
  volume={60},
  number={1},
  pages={396--422},
  year={2022},
  publisher={SIAM}
}

@inproceedings{niethammer2007global,
  title={Global medical shape analysis using the {L}aplace--{B}eltrami spectrum},
  author={Niethammer, Marc and Reuter, Martin and Wolter, Franz-Erich and Bouix, Sylvain and Peinecke, Niklas and Koo, Min-Seong and Shenton, Martha E},
  booktitle={Medical Image Computing and Computer-Assisted Intervention--MICCAI 2007: 10th International Conference, Brisbane, Australia, October 29-November 2, 2007, Proceedings, Part I 10},
  pages={850--857},
  year={2007},
  organization={Springer}
}

@article{reuter2006laplace,
  title={Laplace--{B}eltrami spectra as ‘{S}hape--{DNA}’ of surfaces and solids},
  author={Reuter, Martin and Wolter, Franz-Erich and Peinecke, Niklas},
  journal={Computer-Aided Design},
  volume={38},
  number={4},
  pages={342--366},
  year={2006},
  publisher={Elsevier}
}

@incollection{glowinski2008computing,
  title={Computing the eigenvalues of the {L}aplace--{B}eltrami operator on the surface of a torus: A numerical approach},
  author={Glowinski, Roland and Sorensen, Danny C},
  booktitle={Partial Differential Equations: Modeling and Numerical Simulation},
  pages={225--232},
  year={2008},
  publisher={Springer}
}

@article{brandman2008level,
  title={A level-set method for computing the eigenvalues of elliptic operators defined on compact hypersurfaces},
  author={Brandman, Jeremy},
  journal={Journal of Scientific Computing},
  volume={37},
  pages={282--315},
  year={2008},
  publisher={Springer}
}

@article{hamidian2019surface,
  title={Surface registration with eigenvalues and eigenvectors},
  author={Hamidian, Hajar and Zhong, Zichun and Fotouhi, Farshad and Hua, Jing},
  journal={IEEE Trans. Vis. Comput. Graph.},
  volume={26},
  number={11},
  pages={3327--3339},
  year={2019},
  publisher={IEEE}
}

@inproceedings{levy2006laplace,
  title={Laplace--{B}eltrami eigenfunctions towards an algorithm that ``understands'' geometry},
  author={L{\'e}vy, Bruno},
  booktitle={IEEE International Conference on Shape Modeling and Applications 2006 (SMI'06) 13},
  year={2006},
  organization={IEEE}
}

@article{alikakos1996critical,
  title={Critical spectrum and stability of interfaces for a class of reaction-diffusion equations},
  author={Alikakos, Nicholas D and Fusco, Giorgio and Stefanopoulos, Vagelis},
  journal={Journal of Differential Equations},
  volume={126},
  number={1},
  eid={106},
  year={1996},
  publisher={Academic Press}
}

@article{harrell1996second,
  title={On the second eigenvalue of the {L}aplace operator penalized by curvature},
  author={Harrell, Evans M, II},
  journal={Differential Geometry and its Applications},
  volume={6},
  number={4},
  pages={397--400},
  year={1996},
  publisher={Elsevier}
}

@article{harrell1998laplace,
  title={On the {L}aplace operator penalized by mean curvature},
  author={Harrell, Evans M, II and Loss, Michael},
  journal={Commun. Math. Phys.},
  volume={195},
  pages={643--650},
  year={1998},
  publisher={Springer}
}

@article{da1981quantum,
  title={Quantum mechanics of a constrained particle},
  author={da Costa, R. C. T.},
  journal={Physical Review A},
  volume={23},
  number={4},
  year={1981},
  publisher={APS}
}

@article{duclos1995curvature,
  title={Curvature-induced bound states in quantum waveguides in two and three dimensions},
  author={Duclos, Pierre and Exner, Pavel},
  journal={Reviews in Mathematical Physics},
  volume={7},
  number={01},
  pages={73--102},
  year={1995},
  publisher={World Scientific}
}

@inproceedings{exner1999optimal,
  title={Optimal eigenvalues for some {L}aplacians and {S}chr{\"o}dinger operators depending on curvature},
  author={Exner, Pavel and Harrell, Evans M, II and Loss, Michael},
  booktitle={Mathematical Results in Quantum Mechanics: QMath7 Conference, Prague, June 22--26, 1998},
  pages={47--58},
  year={1999},
  organization={Springer}
}

@article{venn2024meshfree,
  title={A Meshfree Method for Eigenvalues of Differential Operators on Surfaces, Including Steklov Problems},
  author={Venn, Daniel R and Ruuth, Steven J},
  journal={arXiv preprint arXiv:2410.04336},
  year={2024}
}

@article{nwogbaga2025protein,
  title={Protein diffusion controls how single cells respond to electric fields},
  author={Nwogbaga, Ifunanya and Belliveau, Nathan M and Singh, Amit R and Sun, Daiyue and Mulenga, Natasha and Theriot, Julie A and Camley, Brian A},
  journal={bioRxiv},
  year={2025},
  publisher={Cold Spring Harbor Laboratory}
}

@article{lu2021geometrically,
  title={A geometrically consistent trace finite element method for the {L}aplace--{B}eltrami eigenvalue problem},
  author={Lu, Song and Xu, Xianmin},
  journal={arXiv preprint arXiv:2108.02434},
  year={2021}
}

@article{wu2022surface,
  title={Surface eigenvalues with lattice-based approximation in comparison with analytical solution},
  author={Wu, Yingying and Wu, Tianqi and Yau, Shing-Tung},
  journal={arXiv preprint arXiv:2203.03603},
  year={2022}
}

@article{lee2024simple,
  title={A Simple Embedding Method for the {L}aplace--{B}eltrami Eigenvalue Problem on Implicit Surfaces},
  author={Lee, Young Kyu and Leung, Shingyu},
  journal={Communications on Applied Mathematics and Computation},
  volume={6},
  number={2},
  pages={1189--1216},
  year={2024},
  publisher={Springer}
}

@article{horning2022twice,
  title={Twice is enough for dangerous eigenvalues},
  author={Horning, Andrew and Nakatsukasa, Yuji},
  journal={SIAM Journal on Matrix Analysis and Applications},
  volume={43},
  number={1},
  pages={68--93},
  year={2022},
  publisher={SIAM}
}

@article{belkin2003laplacian,
  title={Laplacian eigenmaps for dimensionality reduction and data representation},
  author={Belkin, Mikhail and Niyogi, Partha},
  journal={Neural Computation},
  volume={15},
  number={6},
  pages={1373--1396},
  year={2003},
  publisher={MIT Press}
}

@article{coifman2006diffusion,
  title={Diffusion maps},
  author={Coifman, Ronald R and Lafon, St{\'e}phane},
  journal={Applied and Computational Harmonic Analysis},
  volume={21},
  number={1},
  pages={5--30},
  year={2006},
  publisher={Elsevier}
}

@article{reuter2009discrete,
  title={Discrete {L}aplace--{B}eltrami operators for shape analysis and segmentation},
  author={Reuter, Martin and Biasotti, Silvia and Giorgi, Daniela and Patan{\`e}, Giuseppe and Spagnuolo, Michela},
  journal={Computers \& Graphics},
  volume={33},
  number={3},
  pages={381--390},
  year={2009},
  publisher={Elsevier}
}

@inproceedings{seo2010heat,
  title={Heat kernel smoothing using {L}aplace--{B}eltrami eigenfunctions},
  author={Seo, Seongho and Chung, Moo K and Vorperian, Houri K},
  booktitle={Medical Image Computing and Computer-Assisted Intervention--MICCAI 2010: 13th International Conference, Beijing, China, September 20-24, 2010, Proceedings, Part III 13},
  pages={505--512},
  year={2010},
  organization={Springer}
}

@article{macdonald2011solving,
  title={Solving eigenvalue problems on curved surfaces using the closest point method},
  author={Macdonald, Colin B and Brandman, Jeremy and Ruuth, Steven J},
  journal={Journal of Computational Physics},
  volume={230},
  number={22},
  pages={7944--7956},
  year={2011},
  publisher={Elsevier}
}

@book {reed1972methods,
    author = {Reed, M. and Simon, B.},
     title = {Methods of Modern Mathematical Physics. {I}: {F}unctional {A}nalysis},
   edition = {Second},
 publisher = {Academic Press},
      year = {1980}
}

@book{billingsley2013convergence,
  title={Convergence of Probability Measures},
  author={Billingsley, P.},
  year={1999},
  publisher={John Wiley \& Sons}
}

@Article{maas2005wave,
  author    = {Maas, Leo R. M.},
  journal   = {International Journal of Bifurcation and Chaos},
  title     = {Wave attractors: {L}inear yet nonlinear},
  year      = {2005},
  number    = {09},
  pages     = {2757--2782},
  volume    = {15},
  publisher = {World Scientific Pub Co Pte Ltd}
}

@Article{sibgatullin2019internal,
  author    = {Sibgatullin, I. N. and Ermanyuk, E. V.},
  journal   = {Journal of Applied Mechanics and Technical Physics},
  title     = {Internal and inertial wave attractors: {A} review},
  year      = {2019},
  number    = {2},
  pages     = {284--302},
  volume    = {60},
  publisher = {Springer}
}

@Article{ralston1973stationary,
  author    = {J. V. Ralston},
  journal   = {Journal of Mathematical Analysis and Applications},
  title     = {On stationary modes in inviscid rotating fluids},
  year      = {1973},
  number    = {2},
  pages     = {366--383},
  volume    = {44},
  publisher = {Elsevier}
}

@book{craioveanu_2001, 
	author={M. Craioveanu and M. Puta and T. M. Rassias},
	title={{Old and New Aspects in Spectral Geometry}},
	publisher={{Springer--Science+Business Media, B.V.}},
	edition={1st},
	year={2001}}

@book{chavel_1984,
	author={Isaac Chavel},
	title={Eigenvalues in Riemannian Geometry},
	publisher={Academic Press},
	edition={1st},
	year={1984}}

@book{berger_2000,
	author={Marcel Berger},	
	title={A Panoramic View of Riemannian Geometry},
	publisher={Springer Verlag},
	edition={1st},
	year={2000}}

@book{kato_1980,
	author={Tosio Kato},	
	title={Perturbation Theory for Linear Operators},
	publisher={Springer Verlag},
	edition={2nd},
	year={1980}}

@book{davies_2007,
	author={E. Brian Davies},
	title={Linear Operators and their Spectra},
	publisher={Cambridge University Press},
	edition={1st},
	year={2007}
}

@article{kallianpur1971spectral,
  title={Spectral theory of stationary H-valued processes},
  author={Kallianpur, G and Mandrekar, V},
  journal={Journal of Multivariate Analysis},
  volume={1},
  number={1},
  pages={1--16},
  year={1971},
  publisher={Elsevier}
}

@article{girardin2003semigroup,
  title={Semigroup stationary processes and spectral representation},
  author={Girardin, Valerie and Senoussi, Rachid},
  journal={Bernoulli},
  volume={9},
  number={5},
  pages={857--876},
  year={2003},
  publisher={Bernoulli Society for Mathematical Statistics and Probability}
}

@book{stoica2005spectral,
  title={Spectral Analysis of Signals},
  author={Stoica, Petre and Moses, Randolph L and others},
  volume={452},
  year={2005},
  publisher={Prentice Hall}
}

@article{contFEAST,
author = {Horning, Andrew and Townsend, Alex},
title = {{FEAST} for Differential Eigenvalue Problems},
journal = {SIAM Journal on Numerical Analysis},
volume = {58},
number = {2},
pages = {1239--1262},
year = {2020}
}

@book{taylor1996partial,
  title={Partial differential equations. 1, Basic theory},
  author={Taylor, Michael Eugene},
  year={1996},
  publisher={Springer}
}

@article{colbrook2023computing,
  title={Computing spectral properties of topological insulators without artificial truncation or supercell approximation},
  author={Colbrook, Matthew J and Horning, Andrew and Thicke, Kyle and Watson, Alexander B},
  journal={IMA Journal of Applied Mathematics},
  volume={88},
  number={1},
  pages={1--42},
  year={2023},
  publisher={Oxford University Press}
}

@article{hale2008new,
  title={New quadrature formulas from conformal maps},
  author={Hale, Nicholas and Trefethen, Lloyd N},
  journal={SIAM Journal on Numerical Analysis},
  volume={46},
  number={2},
  pages={930--948},
  year={2008},
  publisher={SIAM}
}

@article{traprule,
author = {Trefethen, Lloyd N. and Weideman, J. A. C.},
title = {The Exponentially Convergent Trapezoidal Rule},
journal = {SIAM Review},
volume = {56},
number = {3},
pages = {385--458},
year = {2014}
}

@inproceedings{colbrook2022specsolve,
  title={{SpecSolve}: Spectral methods for spectral measures},
  author={Colbrook, Matthew J and Horning, Andrew},
  booktitle={Spectral and High Order Methods for Partial Differential Equations ICOSAHOM 2020+ 1: Selected Papers from the ICOSAHOM Conference, Vienna, Austria, July 12-16, 2021},
  pages={183--195},
  year={2022},
  organization={Springer}
}

@book{delaubenfels2006existence,
  title={Existence Families, Functional Calculi and Evolution Equations},
  author={DeLaubenfels, Ralph},
  year={2006},
  publisher={Springer}
}

@book{mennicken2003non,
  title={Non-Self-Adjoint Boundary Eigenvalue Problems},
  author={Mennicken, Reinhard and M{\"o}ller, Manfred},
  volume={192},
  year={2003},
  series={North-Holland Mathematical Studies},
  publisher={Elsevier}
}

@article{wilkening2015spectral,
  title={A spectral transform method for singular {S}turm--{L}iouville problems with applications to energy diffusion in plasma physics},
  author={Wilkening, Jon and Cerfon, Antoine},
  journal={SIAM Journal on Applied Mathematics},
  volume={75},
  number={2},
  pages={350--392},
  year={2015},
  publisher={SIAM}
}

@article{nigam2025intersection,
  title={At the intersection of Numerical Analysis and Spectral Geometry},
  author={Nigam, Nilima},
  journal={arXiv preprint arXiv:2512.25012},
  year={2025}
}

@article{boulton2012generalized,
  title={Generalized {W}eyl theorems and spectral pollution in the {G}alerkin method},
  author={Boulton, Lyonell and Boussaid, Nabile and Lewin, Mathieu},
  journal={Journal of Spectral Theory},
  volume={2},
  number={4},
  pages={329--354},
  year={2012}
}

@article{boulton2016spectral,
  title={Spectral pollution and eigenvalue bounds},
  author={Boulton, Lyonell},
  journal={Applied Numerical Mathematics},
  volume={99},
  pages={1--23},
  year={2016},
  publisher={Elsevier}
}

@article{nigam2020proof,
  title={A proof via finite elements for {S}chiffer’s conjecture on a regular pentagon},
  author={Nigam, Nilima and Siudeja, Bart{\l}omiej and Young, Benjamin},
  journal={Foundations of Computational Mathematics},
  volume={20},
  number={6},
  pages={1475--1504},
  year={2020},
  publisher={Springer}
}

@article{jakobson2005large,
  title={How large can the first eigenvalue be on a surface of genus two?},
  author={Jakobson, Dmitry and Levitin, Michael and Nadirashvili, Nikolai and Nigam, Nilima and Polterovich, Iosif},
  journal={International Mathematics Research Notices},
  volume={2005},
  number={63},
  pages={3967--3985},
  year={2005},
  publisher={Hindawi Publishing Corporation}
}

@article{colbrook2025computing,
  title={Computing generalized eigenfunctions in rigged {H}ilbert spaces},
  author={Colbrook, Matthew J and Horning, Andrew and Xie, Tianyiwa},
  journal={Pure and Applied Analysis},
  volume={7},
  number={2},
  pages={413--443},
  year={2025},
  publisher={Mathematical Sciences Publishers}
}

@book{stone1932linear,
  title={Linear Transformations in Hilbert Space and their Applications to Analysis},
  author={Stone, Marshall Harvey},
  volume={15},
  series={Colloquium Publications},
  year={1932},
  publisher={American Mathematical Soc.}
}

@book{hall2013quantum,
  title={Quantum Theory for Mathematicians},
  author={Hall, Brian C},
  year={2013},
  publisher={Springer}
}

@article{schrweyl,
 ISSN = {00029947},
 author = {Daniel Ray},
 journal = {Transactions of the American Mathematical Society},
 number = {2},
 pages = {299--321},
 publisher = {American Mathematical Society},
 title = {On Spectra of Second--Order Differential Operators},
 volume = {77},
 year = {1954}
}

@article{colbrook2025avoiding,
  title={Avoiding discretization issues for nonlinear eigenvalue problems},
  author={Colbrook, Matthew J and Townsend, Alex},
  journal={SIAM Journal on Matrix Analysis and Applications},
  volume={46},
  number={1},
  pages={648--675},
  year={2025},
  publisher={SIAM}
}

@article{colbrook2024computation,
  title={On the computation of geometric features of spectra of linear operators on {H}ilbert spaces},
  author={Colbrook, Matthew J},
  journal={Foundations of Computational Mathematics},
  volume={24},
  number={3},
  pages={723--804},
  year={2024},
  publisher={Springer}
}

@book{trefethen2020spectra,
  title={Spectra and Pseudospectra: The Behavior of Nonnormal Matrices and Operators},
  author={Trefethen, Lloyd N and Embree, Mark},
  year={2020},
  publisher={Princeton University Press}
}

@article{weylmore,
	title="{\"U}ber die {R}andwertaufgabe der {S}trahlungstheorie und
	asymptotische {S}pektralgeometrie",
	author={Hermann Weyl},
	journal={J. Reine Angew. Math},
	volume={143},
	pages={177--202},
	year={1913}
}

@article{specsolve,
author = {Colbrook, Matthew J. and Horning, Andrew and Townsend, Alex},
title = {Computing Spectral Measures of Self--Adjoint Operators},
journal = {SIAM Review},
volume = {63},
number = {3},
pages = {489--524},
year = {2021}
}

@article{hemisphere,
 ISSN = {03797570},
 author = {SEUNG-JIN BANG},
 journal = {Chinese Journal of Mathematics},
 number = {4},
 pages = {237--245},
 publisher = {Mathematical Society of the Republic of China},
 title = {EIGENVALUES OF THE {L}APLACIAN ON A GEODESIC BALL IN THE n-SPHERE},
 volume = {15},
 year = {1987}
}

@article{100years,
	author={Victor Ivrii},
	title={{100 years of Weyl's law}},
	journal={Bull. Matt. Sci.},
	volume={6},
	pages={379--452},
	year={2016}
}

@article{OriginalFEAST,
  title = {Density--matrix-based algorithm for solving eigenvalue problems},
  author = {Polizzi, Eric},
  journal = {Phys. Rev. B},
  volume = {79},
  issue = {11},
  eid = {115112},
  numpages = {6},
  year = {2009},
  month = {Mar},
  publisher = {American Physical Society}
}

@article{McLachlan_Quispel_2002, 
title={Splitting methods}, 
volume={11}, 
DOI={10.1017/S0962492902000053}, 
journal={Acta Numerica}, 
author={McLachlan, Robert I. and Quispel, G. Reinout W.}, 
year={2002}, 
pages={341--434}
}

@article{quasiQR,
author = {Townsend, Alex  and Trefethen, Lloyd N. },
title = {Continuous analogues of matrix factorizations},
journal = {Proceedings of the Royal Society A: Mathematical, Physical and Engineering Sciences},
volume = {471},
number = {2173},
year = {2015}
}

@article{randfun,
author = {Filip, Silviu and Javeed, Aurya and Trefethen, Lloyd N.},
title = {Smooth Random Functions, Random {ODE}s, and {G}aussian Processes},
journal = {SIAM Review},
volume = {61},
number = {1},
pages = {185--205},
year = {2019}
}

@article{fortunato2024high,
  title={A high--order fast direct solver for surface {PDE}s},
  author={Fortunato, Daniel},
  journal={SIAM Journal on Scientific Computing},
  volume={46},
  number={4},
  pages={A2582--A2606},
  year={2024},
  publisher={SIAM}
}

@article{dyatlov2024mathematics,
  title={Mathematics of internal waves in a two-dimensional aquarium},
  author={Dyatlov, Semyon and Wang, Jian and Zworski, Maciej},
  journal={Analysis \& PDE},
  volume={18},
  number={1},
  pages={1--92},
  year={2024},
  publisher={Mathematical Sciences Publishers}
}

@misc{surfacefun,
	author = {D. Fortunato},
	title = {https://surfacefun.readthedocs.io},
	year = {2023}
}

@article{beyn,
title = {An integral method for solving nonlinear eigenvalue problems},
journal = {Linear Algebra and its Applications},
volume = {436},
number = {10},
pages = {3839--3863},
year = {2012},
issn = {0024-3795},
author = {Wolf-Jürgen Beyn}
}

@article{keldysh,
  title={The nonlinear eigenvalue problem},
  author={Stefan G{\"u}ttel and Françoise Tisseur},
  journal={Acta Numerica},
  year={2017},
  volume={26},
  pages={1--94}
}

@article{brennan2023contour,
  title={Contour integral methods for nonlinear eigenvalue problems: A systems theoretic approach},
  author={Brennan, Michael C and Embree, Mark and Gugercin, Serkan},
  journal={SIAM Review},
  volume={65},
  number={2},
  pages={439--470},
  year={2023},
  publisher={SIAM}
}

@book{henrot2018shape,
  title={Shape Variation and Optimization},
  author={Henrot, Antoine and Pierre, Michel},
  year={2018},
  publisher={European Mathematical Society}
}

@book{delfour2011shapes,
  title={Shapes and Geometries: Metrics, Analysis, Differential Calculus, and Optimization},
  author={Delfour, Michel C and Zol{\'e}sio, J-P},
  year={2011},
  publisher={SIAM}
}

@article{sci_l2N,
author = {Colbrook, Matthew J. and Hansen, Anders C.},
year = {2022},
month = {11},
title = {{The foundations of spectral computations via the Solvability Complexity Index hierarchy}},
journal = {Journal of the European Mathematical Society}
}

@misc{nlep_coll,
	title={An Updated Set of Nonlinear Eigenvalue Problems},
	author={Higham, Nicholas J. and Negri Porzio, Gian Maria and Tisseur, Francoise},
	year={2019},
	url={https://eprints.maths.manchester.ac.uk/2699/},
	howpublished={MIMS Preprint}
    }

@article{weyl,
	title={{\"U}ber die asymptotische {V}erteilung der {E}igenwerte},
	author={Hermann Weyl},
	pages={110--117},
	year={1911},
	journal={Nachr. Konigl. Ges. Wiss. G\"ottingen}
	}

@book{Frankel2012,
  title = {The Geometry of Physics: {A}n Introduction},
  author = {Frankel, Theodore},
  publisher = {Cambridge University Press},
  year = {2012},
  edition = {3rd}
}

@article{zvonek2025conthutch++,
  title={Cont{H}utch++: Stochastic trace estimation for implicit integral operators},
  author={Zvonek, Jennifer and Horning, Andrew and Townsend, Alex},
  journal={SIAM Journal on Numerical Analysis},
  volume={63},
  number={1},
  pages={334--359},
  year={2025},
  publisher={SIAM}
}

@inproceedings{meyer2021hutch++,
  title={Hutch++: {O}ptimal stochastic trace estimation},
  author={Meyer, Raphael A and Musco, Cameron and Musco, Christopher and Woodruff, David P},
  booktitle={Symposium on Simplicity in Algorithms (SOSA)},
  pages={142--155},
  year={2021},
  organization={SIAM}
}

@article{persson2022improved,
  title={Improved variants of the {H}utch++ algorithm for trace estimation},
  author={Persson, David and Cortinovis, Alice and Kressner, Daniel},
  journal={SIAM Journal on Matrix Analysis and Applications},
  volume={43},
  number={3},
  pages={1162--1185},
  year={2022},
  publisher={SIAM}
}

@article{Gillman2015,
  title = {A Spectrally Accurate Direct Solution Technique for Frequency-Domain Scattering Problems with Variable Media},
  author = {Gillman, Adrianna and Barnett, Alex H. and Martinsson, Per-Gunnar},
  year = 2015,
  month = mar,
  journal = {BIT Numerical Mathematics},
  volume = {55},
  number = {1},
  pages = {141--170},
  doi = {10.1007/s10543-014-0499-8}
}

@software{chunkie,
  author  = {Askham, Travis and Rachh, Manas and O'Neil, Michael and Hoskins, Jeremy and Fortunato, Daniel and Jiang, Shidong and Fryklund, Fredrik and Goodwill, Tristan and Wang, Hai Yang and Zhu, Hai},
  title   = {chunk{IE}: a {MATLAB} integral equation toolbox},
  year    = {2024},
  month   = jun,
  version = {1.0.0},
  url     = {https://github.com/fastalgorithms/chunkie},
  note    = {Documentation: \url{https://chunkie.readthedocs.io/}}
}

\end{document}